**Florentin Smarandache, V. Christianto,**
**Fu Yuhua, R. Khrapko, J. Hutchison**

**Unfolding the Labyrinth:**
**Open Problems in Physics, Mathematics, Astrophysics and**
**Other Areas of Science**

**2006**

Abstract

Throughout this book, we discuss some open problems in various branches of science, including mathematics, theoretical physics, astrophysics, geophysics etc. It is of our hope that some of the problems discussed in this book will find their place either in theoretical exploration or further experiments, while some parts of these problems may be found useful for scholarly stimulation.

The present book is also intended for young physics and mathematics fellows who will perhaps find the unsolved problems described here are at least worth pondering. If this book provides only a few highlights of plausible solutions, it is merely to keep the fun of readers in discovering the answers by themselves.

Bon voyage!

# Unfolding the Labyrinth: Open Problems in Physics, Mathematics, Astrophysics, and Other Areas of Science

F. Smarandache, V. Christianto, Fu Yuhua, R. Khrapko, J. Hutchison

## Contents

# Preface

The reader will find herein a collection of unsolved problems in mathematics and the physical sciences. Theoretical and experimental domains have each been given consideration. The authors have taken a liberal approach in their selection of problems and questions, and have not shied away from what might otherwise be called speculative, in order to enhance the opportunities for scientific discovery.

Progress and development in our knowledge of the structure, form and function of the Universe, in the true sense of the word, its beauty and power, and its timeless presence and mystery, before which even the greatest intellect is awed and humbled, can spring forth only from an unshackled mind combined with a willingness to imagine beyond the boundaries imposed by that ossified authority by which science inevitably becomes, as history teaches us, barren and decrepit.

Revealing the secrets of Nature, so that we truly see '*the sunlit plains extended, and at night the wondrous glory of the everlasting stars*'[*], requires far more than mere technical ability and mechanical dexterity learnt form books and consensus. The dustbin of scientific history is replete with discredited consensus and the grand reputations of erudite reactionaries. Only by boldly asking questions, fearlessly, despite opposition, and searching for answers where most have not looked for want of courage and independence of thought, can one hope to discover for one's self. From nothing else can creativity blossom and grow, and without which the garden of science can only aspire to an overpopulation of weeds.

Stephen J. Crothers
Queensland, Australia
*Progress in Physics Journal,* http://www.ptep-online.com
14th July 2006.

---

[*] A. B. (Banjo) Patterson's 'Clancy of the Overflow'.

# Foreword

"**...T**he central problem is unsolvable: the enumeration, even if only partial**...** I saw the Aleph from all points; I saw the earth in the Aleph and in the earth the Aleph once more and the earth in the Aleph; I saw my face and my viscera;... because my eyes had seen that conjectural and secret object whose name men usurp but which no man has gazed on: *the inconceivable universe.*"---**Aleph, J.L. Borges**

Partly inspired by a well-known paper by Ginzburg [1], the present book discusses various open problems in different areas of Science, including Physics, Mathematics, Geophysics, Astrophysics etc. Therefore this book could be viewed as an extended form of the aforementioned paper of Prof. V. Ginzburg [1]. Nonetheless the writers attempt herein to look deeper into what appear to us as open problems.

Throughout the book the writers describe unsolved problems in various fields of science, with the hope that these problems might perhaps inspire other researchers in their quest of finding new answers. The writers have made their best effort to write the problems here in a refreshing style. This is why the present book is recommended for researchers and graduate students who are looking for potentially new, breakthrough ideas in physics or applied mathematics.

Needless to say, some of the questions posed here will sound a bit weird, if not completely incomprehensible. Some of them also contain things that the reader may not think easy to follow. For instance, a reader might find the extension of 'quark' ideas incomprehensible, because the quarks themselves may not pop-out easily in our daily dose of reality (because of the confinement problem). As Heisenberg once said, more or less: "If quarks exist then we have redefined the word 'exist'." These belong to ideas that perhaps may have a chance to stimulate the neurons inside our brains.

We would like to thank the reviewers of this book, Profs. T. Love and A. Kaivarainen, and also S. Crothers, for their patience in reading the draft version of this book, and for their comments. We are also grateful for valu-able discussions with numerous colleagues from all over the world, for some of the questions in this book were inspired by their comments, in particular Profs. C. Castro, M. Pitkanen, E. Scholz, E. Bakhoum, R.M. Kiehn, Dong Choi, Chen I-wan, D. Rabounski and numerous others. And also special thanks to peer-reviewers for critically reading our papers and suggesting improvements. We also thank Robert Davic for his comments on the Brightsen model.

All in all, hopefully, these unsolved problems could motivate other young researchers in their journey for unfolding the *Labyrinth* of Nature.

FS, VC, FY, RK, JH
August 28th 2006

# 1    Unsolved Problems in Theoretical Physics

*The only way of discovering the limits of the possible is to venture a little way past them into the impossible.--Arthur C. Clarke*

## 1.1    Problems in elementary particles etc.

It is known that Quantum Mechanics is the cornerstone of more recent theories intended to describe the nature of elementary particles, including Quantum Field Theory, Quantum Electrodynamics, Quantum Chromodynamics, and so forth.

But Quantum Mechanics in its present form also suffers from the same limitations as the foundations of logic; therefore it is not surprising that there are difficult paradoxes that astonished physicists for almost eight decades. Some of these paradoxes are:

- Wigner's friend;
- Einstein-Podolski-Rosen paradox;
- Schrödinger's cat paradox.

While numerous attempts have been made throughout the past eight decades to solve all these paradoxes, it seems that only a few of the present theories can solve these paradoxes completely.

As a result, it is therefore not so surprising to find that both Quantum Electrodynamics (QED) and also Quantum Chromodynamics have their own problems. For instance Dirac and Feynman never accepted QED as a complete theory on its own (as Feynman put it: "*It's like sweeping under the rug.*"). This is why Dirac attempted to propose a new theory to replace QED, albeit the result has not been so successful. Recently, there have been some attempts to reconsider Dirac's new theory (1951) in the light of the biquaternions.[2]

Similarly, other big questions in theoretical / particle physics can be described as follows:

(i)        Is there a Dirac æther fluid? [2]

(ii)       Can Dirac's recent theory 1951 solve the infinity problem?

(iii)  Does Dirac's new electron theory 1951 reconcile the quantum mechanical view with the electrodynamical view of the electron? [2]

(iv)  What is the dynamical mechanism behind the Koide mixing matrix of the lepton mass formula? [3][4][5]

(v)  Does the neutrino have mass? [6][7] [8][9]

(vi)  Does the rishon or preon model of elementary particles give better prediction than conventional Quantum Chromodynamics Theory? [10][11][12]

(vii)  Is there a physical explanation of quark confinement?

(viii)  Is there a theoretical link between Quantum Chromodynamics and quantum fluid dynamics?

Harari is a physicist who made one of the earliest attempts to develop a preon [11] model to explain the phenomena appearing in hadrons. Harari proposed the rishon model in order to simplify the quark model of Gell-Mann. The model has two kinds of fundamental particles called "rishon" (which means "primary" in Hebrew).[11] They are T (Third for charge 1/3e or Tohu from "unformed" in Hebrew in Genesis) and V (Vanishes for charge 0 or va-Vohu which means "void" in Hebrew in Genesis).

All leptons and all flavours of quarks are combinations of three rishons. They are as follows: These groups of three rishons have spin ½. They are as follows: TTT=positron; VVV=electron neutrino; TTV, TVT and VTT=three colors of u quarks; TVV, VTV and VVT=three colors of d antiquarks. Each rishon has its antiparticle, therefore: ttt=electron; vvv=anti-electron neutrino; ttv, tvt, vtt=three colors of anti-u quarks; vvt, vtv, tvv=three colors of d quarks.

Furthermore, the search for a neutrino mass has recently become a big industry in recent years. "Today's neutrino detectors, kept deep underground to avoid stray particles on Earth's surface, may contain thousands of tons of fluid. While trillions of neutrinos pass through the fluid every day, only a few dozen are likely to be detected. Scientists have discovered that there are three types of neutrinos, each associated with a different charged particle for which it is named. Thus they are called the electron neutrino, muon neutrino, and tau neutrino. The first type of neutrino to be discovered was the electron neutrino, in 1959. The muon neutrino was discovered in 1962. The tau neutrino has yet to be directly observed. It was inferred from the existence of the tau particle itself, which was discovered in 1978. The tau particle is involved in decay reactions with the same imbalance that Pauli solved for beta decay by postulating the electron neutrino."[14]

As we see, some of these questions are very tough, and it is likely they will trigger new kinds of experiments.

Now from these 'known' questions, we can also ask some new questions for further development of theoretical physics:

(i)     Is it possible to come up with a quantum liquid model of elementary particles? How can it predict the elementary particle masses?

(ii)    Could we find isolated quarks or rishons in Nature?

(iii)   Could we find isolated quarks or rishons in a strong electromagnetic field environment?

(iv)    If Koide's concept of the democratic mixing matrix is proved true, then how could we find fluid a dynamical interpretation of this mixing matrix?

(v)     Could we find a theoretical explanation of quarks / rishons from the viewpoint of multivalued-logic Quantum Mechanics?

(vi)     Could we find a theoretical explanation of quarks / rishons from the viewpoint of Quaternion Quantum Mechanics? If yes, then how could we ascribe physical meaning to a *scalar* in the quaternion field?

(vii)   Is there also quaternion-type symmetry (see Adler's QQM theory, for instance) in neutrino mass?

(viii)  Could we find a theoretical explanation of quarks / rishons from the viewpoint of the Gross-Pitaevskii equation for a rotating Bose-Einstein Condensate? If yes, then how does the Magnus effect affect the rotational dynamics of the quarks?

(ix)    What is the effect of gravitational field on the charges of quarks and rishons?

(x)     Could we alter the charges or masses of quarks? If yes, then how could it be done?

(xi)    Could we transform the quark charges back into the vacuum surrounding it?

(xii)   Could we transform the quark charges into Energy? How could this process be done? Under what conditions?

(xiii)  Does the neutrino mass could transform into an isolated entity?

(xiv)   Could we find signatures of anti-hydrogen (antimatter hydrogen) in astrophysics?

(xv)    Suppose there is a large *anti-hydrogen star* ---similar to neutron star—in the Cosmos. How will it affect the normal star?

(xvi)   Is anti-hydrogen also formed in normal star, like the Sun? If yes, then what is its signature?

(xvii)  Is anti-hydrogen compatible with the ring-model of the electron? If not, why?

(xviii) Is it possible to derive a ring-model of the electron which is consistent with Dirac's model (1951) and also an anti-hydrogen experiment? (http://www.groupkos.com/mtwain/TheElectron.pdf)

(xix)   What is the actual trajectory of a deuterium nucleus in the context of the ring-model of the electron? (ref. http://groups.yahoo.com/group/NuclearStructure/)

(xx)    Could we find a theoretical basis for Quantum Mechanics and Quantum Electrodynamics which automatically includes anti-hydrogen in the theory?

(xxi)   Could we find neutrinos inside the human body?

## 1.2   Problems related to Unmatter [52]-[70]

Some unsolved problems related to unmatter are as follows:
-    Is it possible to make infinitely many combinations of quarks / anti-quarks and leptons / antileptons?

- Unmatter can combine with matter and/or antimatter and the result may be any of these three. Some unmatter could be in the strong force, hence part of hadrons. Could we find signatures of unmatter in hadrons?
- For the containment of antimatter and unmatter would it be possible to use electromagnetic fields (a container whose walls are electromagnetic fields). But is its duration unknown?
-

We describe further these questions in the following sections.

### 1.2.1    Abstract

As shown herein, experiments have detected unmatter: a new kind of matter whose atoms include both nucleons and anti-nucleons, while their life span was very short, no more than $10^{-20}$ sec. Stable states of unmatter can be built on quarks and anti-quarks: applying the unmatter principle here it is obtained a quantum chromodynamics formula is obtained herein that gives many combinations of unmatter built on quarks and anti-quarks.

In the time since the appearance of my articles defining "matter, antimatter, and unmatter" [53,54], and Dr. S. Chubb's pertinent comment [55] on unmatter, there has been new development in the unmatter topic in the sense that experiments verifying unmatter have been performed.

### 1.2.2    Definition of Unmatter

In short, unmatter is formed by matter and antimatter binding together [53,54]. The building blocks (most elementary particles known today) are 6 quarks and 6 leptons; their 12 antiparticles also exist. Then *unmatter* will be formed by at least a building block and at least an antibuilding block which can bind together.

### 1.2.3    Exotic Atom

If in an atom we substitute one or more particles by other particles of the same charge (constituents) we obtain an exotic atom whose particles are held together due to the electric charge. For example, we can substitute for one or more electrons in ordinary atom, by other *negative particles* (say $\pi^-$, anti-Rho meson, $D^-$, $D_s^-$, muon, tau, $\Omega^-$, $\Delta^-$, etc., generally clusters of quarks and antiquarks whose total charge is negative), or the positively charged nucleus replaced by other *positive particles* (say clusters of quarks and antiquarks whose total charge is positive, etc.).

### 1.2.4    Unmatter Atom

It is possible to define unmatter in a more general way, using the exotic atom. The classical unmatter atoms were formed by particles like (a) electrons, protons, and antineutrons, or (b) antielectrons, antiprotons, and neutrons. In a more general definition, an unmatter atom is a system of particles as above, or such that one or more particles are replaces by other particles of the same charge.

Other categories would be (c) a matter atom wherein one or more (but not all) of the electrons and/or protons are replaced by antimatter particles of the same corresponding charges, and (d) an antimatter atom such that one or more (but not all) of the antielectrons and/or antiprotons are replaced by matter particles of the same corresponding charges.

In a more complicated system we can substitute a particle by an unmatter particle and form an unmatter atom.

Of course, not all of these combinations are stable, semi-stable, or quasi-stable, especially when their time to bind might be longer than their lifespan.

### 1.2.5    Examples of an Unmatter Atom

During 1970-1975 numerous purely experimental verifications were obtained proving that "atom-like" systems built on nucleons (protons and neutrons) and anti-nucleons (anti-protons and anti-neutrons) are real. Such "atoms", where nucleon and anti-nucleon are moving in the opposite sides of the same orbit around the common centre of mass, are very unstable, their life span is no more than $10^{-20}$ sec. Then nucleon and anti-nucleon annihilate into gamma-quanta and other light particles (pions), which cannot be connected with one another, see [58,59,60]. The experiments were performed mainly at Brookhaven National Laboratory (USA) and, partially at CERN (Switzerland), where "proton=>anti-proton" and "anti-proton=>neutron" atoms were observed, denoted by $\overline{p}p$ and $\overline{p}n$ respectively, see Fig 1 and 2.

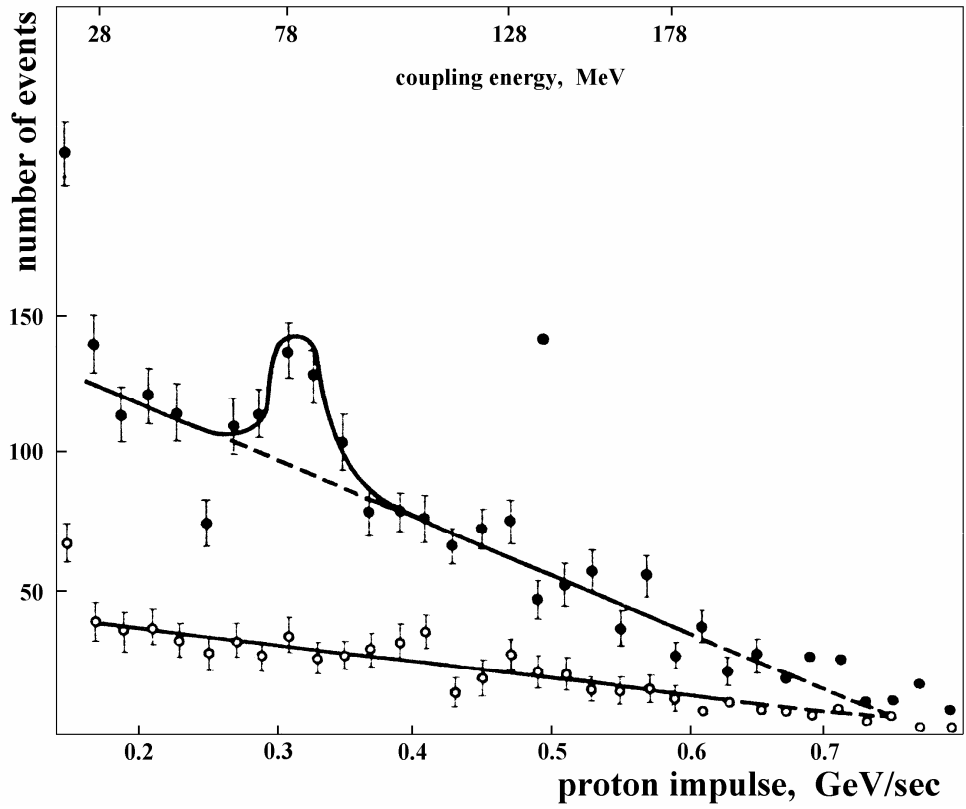

Fig. 1: Spectra of proton impulses in the reaction $\bar{p} + d \rightarrow (\bar{p}n) + p$. The upper arc depicts annihilation of $\bar{p}n$ into an even number of pions, the lower arc --- its annihilation into an odd number of pions. The observed maximum indicates that there is a connected system $\bar{p}n$. Abscissa axis represents the proton impulse in GeV/sec (and the connection energy of the system $\bar{p}n$). Ordinate axis gives the number of events (after [60]).

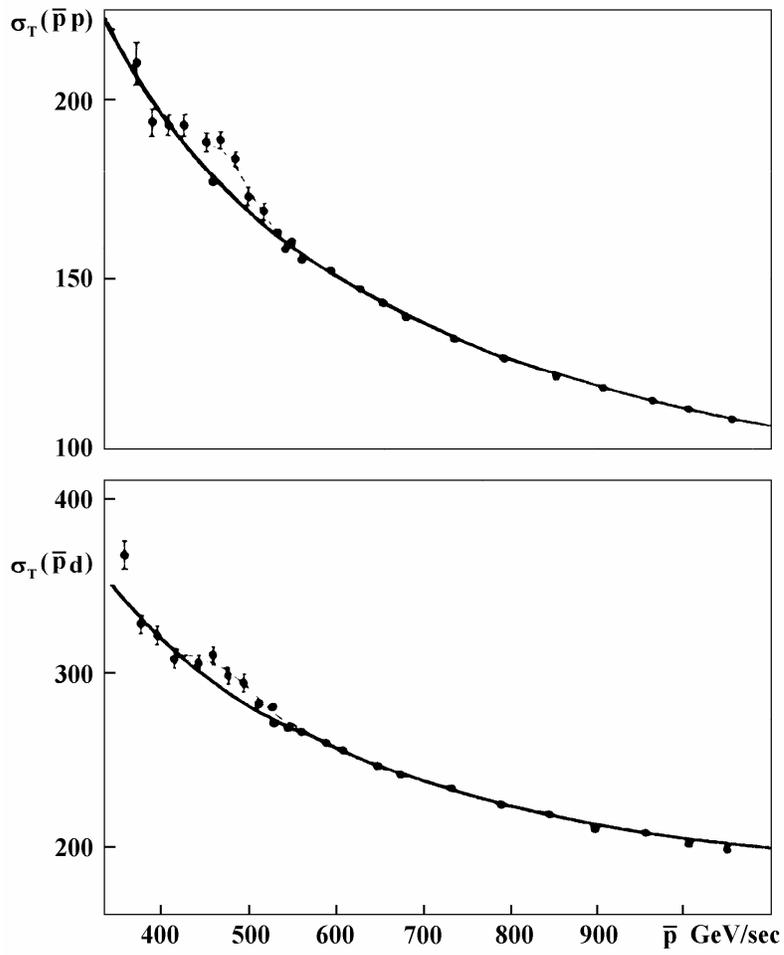

Fig. 2: Probability σ of interaction between $\overline{p}$, $p$ and deuterons $d$ (after from [61]). The presence of a maximum indicates the existence of the resonance state of "nucleon --- anti-nucleon".

After the experiments were completed, the life span of such "atoms" was calculated theoretically in Chapiro's works [61,62,63]. His main idea was that nuclear forces, acting between nucleon and anti-nucleon, can keep them far away from each other, hindering their annihilation. For instance, a proton and anti-proton are located at the opposite side of the same orbit and move around the orbit's centre. If the diameter of their orbit is much larger than the diameter of the "annihilation area", they can be kept from annihilation (see fig. 3). But because the orbit, according to Quantum Mechanics, is an actual cloud spreading far around the average radius, at any radius between the proton and the anti-proton there is a probability that they can meet one another at the annihilation distance. Therefore the nucleon---anti-nucleon system annihilates in any case, as this system is unstable by definition having a life span no more than $10^{-20}$ sec.

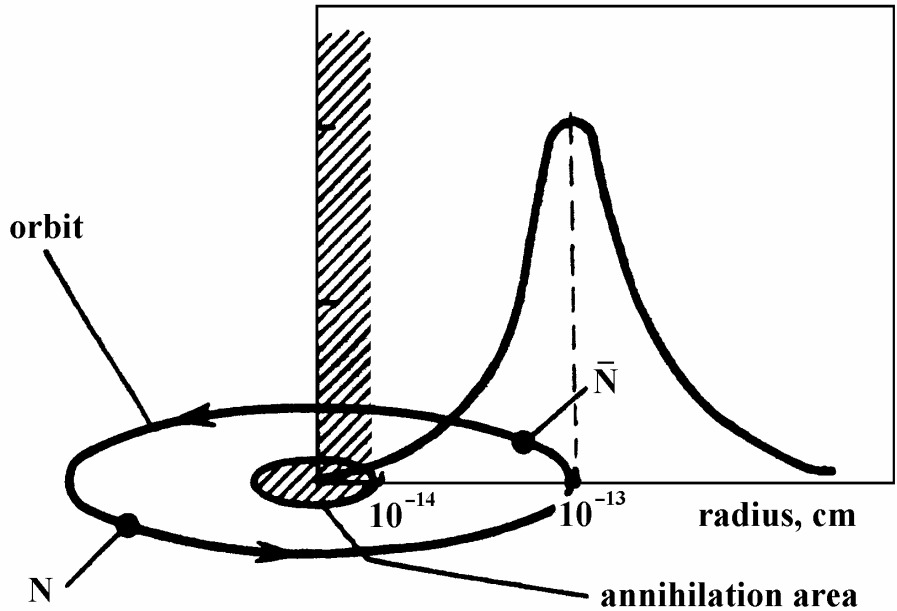

Fig. 3: Annihilation area and the probability arc in a "nucleon --- antinucleon" system (after [63]).

Unfortunately, the researchers limited their investigations to the consideration of $\overline{p}p$ and $\overline{p}n$ nuclei only, for the reason that they, in the absence of a theory, considered $\overline{p}p$ and $\overline{p}n$ "atoms" as only a rare exception, which gives no classes of matter.

Despite Benn Tannenbaum's and Randall J. Scalise's rejections of unmatter and Scalise's personal attack on its author (of unmatter) in a true Ancient Inquisitorial style, under the aegis of MadSci moderator John Link the unmatter does exist, *for example some mesons and antimesons, though for a trifling of a second lifetime, so the pions are unmatter* [which have the composition u^d and ud^, where by u^ we mean anti-up quark, d = down quark, and analogously u = up quark and d^ = anti-down quark, while by ^ means anti], the kaon K$^+$ (us^), K$^-$ (u^s), Phi (ss^), D$^+$ (cd^), D$^0$(cu^), D$_s^+$ (cs^), J/Psi (cc^), B$^-$ (bu^), B$^0$ (db^), B$_s^0$ (sb^), Upsilon (bb^) [where c = charm quark, s = strange quark, b = bottom quark], etc. are unmatter too.

*Also, the pentaquark* ($\Theta^+$), of charge $^+1$, uudds^ (i.e. two quarks up, two quarks down, and one anti-strange quark), at a mass of 1.54 GeV and a narrow width of 22 MeV, is *unmatter*, observed in 2003 at the Jefferson Lab in Newport News, Virginia, in the experiments that involved multi-GeV photons impacting upon a deuterium target. Similar pentaquark evidence was obtained by Takashi Nakano of Osaka University in 2002, by researchers at the ELSA accelerator in Bonn in 1997-1998, and by researchers at ITEP in Moscow in 1986.

Besides $\Theta^+$, evidence has been found in one experiment [56] for other pentaquarks, $\Xi_5^-$(ddssu^) and $\Xi_5^+$(uussd^).

D. S. Carman [57] has reviewed the positive and null evidence for these pentaquarks and their existence is still under investigation.

Let's recall that the *pionium* is formed by a $\pi^+$ and $\pi^-$ mesons, the *positronium* is formed by an antielectron (positron) and an electron in a semistable arrangement, the *protonium* is formed by a proton and an antiproton also semi-stable, the *antiprotonic helium* is formed by an antiproton and electron together with the helium nucleus (semi-stable), and *muonium* is formed by a positive muon and an electron.

Also, the *mesonic atom* is an ordinary atom with one or more of its electrons replaced by negative mesons.

*Strange matter* is ultra-dense matter formed by a big number of strange quarks bounded together with an electron atmosphere (this strange matter is hypothetical).

From the exotic atom, the pionium, positronium, protonium, antiprotonic helium, and muonium are unmatter.

The mesonic atom is unmatter if the electron(s) are replaced by negatively-charged antimesons. Also we can define a mesonic antiatom as an ordinary antiatomic nucleus with one or more of its antielectrons replaced by positively-charged mesons. Hence, this mesonic antiatom is unmatter if the antielectron(s) are replaced by positively-charged mesons. The strange matter can be unmatter if these exists at least an antiquark together with so many quarks in the nucleus. Also, we can define the strange antimatter as formed by a large number of antiquarks bound together with an antielectron cloud around them. Similarly, the strange antimatter can be unmatter if there exists at least one quark together with so many antiquarks in its nucleus.

The bosons and antibosons contribute to the decay of unmatter. There are 13+1 (Higgs boson) known bosons and 14 antibosons at present.

### 1.2.6    Quantum Chromodynamics Formula

In order to save the colourless combinations prevailed in the Theory of Quantum Chromodynamics of quarks and antiquarks in their combinations when binding, we devise the following formula:

$$Q - A \in \pm M3 \qquad\qquad (1)$$

where $M3$ denotes multiples of three, i.e. $\pm M3 = \{3 \cdot k \mid k \in Z\} = \{\dots, -12, -9, -6, -3, 0, 3, 6, 9, 12, \dots\}$, and Q = number of quarks, A = number of antiquarks. But (1) is equivalent to:

$$Q \equiv A \pmod 3 \qquad\qquad (2)$$

(Q is congruent to A modulo 3).

To justify this formula we mention that 3 quarks form a colourless combination, and any multiple of three (**M3**) combination of quarks too, i.e. 6, 9, 12, etc. quarks. In a similar way, 3 antiquarks form a colourless combination, and any multiple of three (**M3**) combination of antiquarks too, i.e. 6, 9, 12, etc. antiquarks. Hence, when we have hybrid combinations of quarks and antiquarks, a quark and an antiquark will annihilate their colours and, therefore, what's left should be a multiple of three number of quarks (in the case when the number of quarks is larger, and the difference in the formula is positive), or a multiple of three number of antiquarks (in the case when the number of antiquarks is larger, and the difference in the formula is negative).

### 1.2.7    Quark-Antiquark Combinations

Let's denote q = quark $\in$ {Up, Down, Top, Bottom, Strange, Charm}, and by a = antiquark $\in$ {Up^, Down^, Top^, Bottom^, Strange^, Charm^}.

Hence, for combinations of n quarks and antiquarks, n $\geq$ 2, colourless prevailing, we have the following possibilities:

- if n = 2, we have: qa (biquark – for example the mesons and antimessons);

- if n = 3, we have qqq, aaa (triquark – for example the baryons and antibaryons);

- if n = 4, we have qqaa (tetraquark);

- if n = 5, we have qqqqa, aaaaq (pentaquark);

- if n = 6, we have qqqaaa, qqqqqq, aaaaaa (hexaquark);

- if n = 7, we have qqqqqaa, qqaaaaa (septiquark);

- if n = 8, we have qqqqaaaa, qqqqqqaa, qqaaaaaa (octoquark);

- if n = 9, we have qqqqqqqqq, qqqqqqaaa, qqqaaaaaa, aaaaaaaaa (nonaquark);

- if n = 10, we have qqqqqaaaaa, qqqqqqqqaa, qqaaaaaaaa (decaquark); etc.

### 1.2.8    Unmatter Combinations

From the above general case we extract the unmatter combinations:

- For combinations of 2 we have: qa (unmatter biquark), [mesons and antimesons]; the number of all possible unmatter combinations will be 6·6 = 36, but not all of them will bind together. It is possible to combine an entity with its mirror opposite and still bind them, such as:,,uu^, dd^, ss^, cc^, bb^ which form mesons. It is possible to combine, unmatter + unmatter = unmatter, as in ud^ + us^ = uud^s^ (of course if they bind together).

- For combinations of 3 (unmatter triquark) we can not form unmatter since the colourless cannot hold.

- For combinations of 4 we have: qqaa (unmatter tetraquark); the number of all possible unmatter combinations will be $6^2 \cdot 6^2 = 1,296$, but not all of them will bind together.

- For combinations of 5 we have: qqqqa, or aaaaq (unmatter pentaquarks); the number of all possible unmatter combinations will be $6^4 \cdot 6 + 6^4 \cdot 6 = 15,552$, but not all of them will bind together.

- For combinations of 6 we have: qqqaaa (unmatter hexaquarks); the number of all possible unmatter combinations will be $6^3 \cdot 6^3 = 46,656$, but not all of them will bind together.

- For combinations of 7 we have: qqqqqaa, qqaaaaa (unmatter septi-quarks); the number of all possible unmatter combinations will be $6^5 \cdot 6^2 + 6^2 \cdot 6^5 = 559,872$, but not all of them will bind together.

- For combinations of 8 we have: qqqqaaaa, qqqqqqqa, qaaaaaaa (unmatter octoquarks); the number of all possible unmatter combinations will be $6^4 \cdot 6^4 + 6^7 \cdot 6^1 + 6^1 \cdot 6^7 = 5,038,848$, but not all of them will bind together.

- For combinations of 9 we have: qqqqqqqaa, qqqaaaaaa (unmatter nonaquarks); the number of all possible unmatter combinations will be $6^6 \cdot 6^3 + 6^3 \cdot 6^6 = 2 \cdot 6^9 = 20,155,392$, but not all of them will bind together.

- For combinations of 10 we have: qqqqqqqqaa, qqqqqaaaaa, qqaaaaaaaa (unmatter decaquarks); the number of all possible unmatter combinations will be $3 \cdot 6^{10} = 181,398,528$, but not all of them will bind together. Etc.

We wonder if it is possible to make infinitely many combinations of quarks / antiquarks and leptons / antileptons. Unmatter can combine with matter and/or antimatter and the result may be any of these three. Some unmatter could be involved in the strong force, and hence a part of hadrons.

**Quantum Chromodynamics Unmatter Formula.**
From formula (2) we derive a particular case in order to characterize the quantum unmatter, and therefore both the quarks and antiquarks should coexist in the same combination:

$$Q \equiv A \pmod 3 \tag{3}$$
$$\text{and } Q \cdot A \neq 0 \text{ (i.e. both Q and A are non-null)},$$

where Q = number of quarks and A = number of antiquarks.

### 1.2.9    Unmatter charge

The charge of unmatter may be positive as in the pentaquark ($\Theta^+$), 0 (as in positronium), or negative as in the anti-Rho meson (u^d) [M. Jordan].

### 1.2.10   Containment

We think that for the containment of antimatter and unmatter it would be possible to use electromagnetic fields (a container whose walls are electromagnetic fields).  But its duration is unknown.

### 1.2.11   Further research

Let's begin with Neutrosophy [70], which is a generalization of dialectics, i.e. not only the opposites are combined but also the neutralities. Why? It is when an idea is propounded, a category of people will accept it, others will reject it, and a third group will ignore it (don't care). But the dynamics between these three categories changes, so somebody accepting it might later reject or ignore it, or someone ignoring it will accept it or reject it, and so on. Similarly for the *dynamics of <A>, <antiA>, <neutA>, where <neutA> means neither <A> nor <antiA>, but in between (neutral).*

Neutrosophy deals not with a kind of di-alectics but a kind of tri-alectics (based on three components: <A>, <antiA>, <neutA>).  Hence unmatter is a kind of neutrality (not referring to the charge) between matter and antimatter, i.e. neither one, nor the other.

In the model of unmatter we may conceive at ungravity, unforce, unenergy, etc. *Ungravity* would be a mixture between gravity and antigravity (for example attracting and rejecting simultaneously or alternatively; or a magnet which changes the + and - poles frequently). *Unforce*. We may consider positive force (in the direction we want), and negative force (repulsive, opposed to the previous). There could be a combination of both positive and negative forces in the same time, or alternating positive and negative, etc.

*Unenergy* would similarly be a combination between positive and negative energies (as the alternating current (a.c.), which periodically reverses its direction in a circuit and whose frequency, f, is independent of the circuit's constants). Would it be possible to construct an alternating-energy generator?

In conclusion: According to the Universal Dialectic, unity is manifests in duality and the duality in unity. "Thus, Unmatter (unity) is experienced as duality (matter vs antimatter). Ungravity (unity) as duality (gravity vs antigravity). Unenergy (unity) as duality (positive energy vs negative energy). and thus also...between duality of being (existence) vs nothingness (antiexistence) must be "unexistence" (or pure unity)." (R. Davic)

## 1.3    Some Unsolved Problems, Questions, and

# Applications of the Brightsen Nucleon Cluster Model

According to the Brightsen Nucleon Cluster Model [1] all nuclides of beta stable isotopes can be described by three fundamental nucleon clusters {NPN,PNP,NP}, with halo clusters (NN,PP,NNN) now experimentally observed. The Brightsen model builds on the early cluster models of the Resonating Group Structure of John Wheeler [2] and the Linus Pauling Close-Packed Spheron Model [3], which predict mathematically that the wave function of a composite nucleus can be viewed quantum mechanically as a combination of partial wave functions that correspond to the multiple ways nucleons (protons, neutrons) can be distributed into close-packed clusters, thus rejecting the standard model Hartree-Fock formalism of average field interactions between independent nucleons in nuclear shells. Presented in this section are a number of unsolved problems, questions, and future experimental pathways based on the Brightsen Nucleon Cluster Model formalism--many additional applications can be gleamed from careful study of the literature cited in the references provided:

1. The Brightsen Model derives the average number of prompt neutrons per fission event for many radioactive isotopes of human importance (U-235, U-233, Pu-239, Pu-241) as well as emission of light charged particles, suggesting that all modes of fission derive from a four step process [4]. Further study of these claims are warranted given the importance of understanding the fission of radioactive isotopes for energy production.

2. The Brightsen Model provides a theoretical pathway for experimentalists to understand the numerous laboratory results of low temperature transformation/low energy reactions, such as the well studied $^{104}$Pd (p, alpha) $^{101}$Rh reaction [5]. Application of the Brightsen Model to low energy fusion reactions as a possible result of interactions between nucleon clusters is of fundamental importance to human energy demands.

3. The Brightsen Model predicts the existence of "unmatter entities" inside nuclei [6], which result from stable and neutral union of matter and antimatter nucleon clusters. As a result, the Brightsen Model predicts that antimatter has corresponding antigravity effects [7]. This prediction can be tested in the future at CERN beginning 2008 using antihydrogen. Once accurate measurements can be made of the gravitational acceleration of antihydrogen, and the results compared with matter hydrogen, if the two forms have opposite acceleration, then a major prediction of the Brightsen

Model will be confirmed (e.g., that antimatter has both anti-gravity effect and anti-mass). If experimentally confirmed, then predictive equations will need to be developed using the Brightsen Model formalism of union of matter and antimatter clusters (e.g., the unsolved mathematical formation of unmatter entities inside nuclei). The importance of this aspect of the Brightsen Model links to the current problem in physics of the missing matter of the universe and possible unification of gravity at relativistic (macroscopic) and quantum (microscopic) states.

4. The Brightsen Model offers a theoretical approach for artificially induced fission of dangerous radioactive nuclei to produce relatively stable elements [5]. In theory, if externally produced electromagnetic radiation can be caused to resonate with the exact magnetic moment of a specific sub-nuclear nucleon cluster (e.g., NPN,PNP,NP nucleon clusters), than an individual nucleon cluster can in theory be excited to a energy such that it is expelled from the nucleus, resulting in transmutation of the parent isotope via fission and/or beta or alpha decay to less radioactive daughter structures. The applications of this process for nuclear energy production are clear and worthy of experimental test.

5. The Brightsen Model predicts that one sub-cluster isodyne [5] of the very stable Helium-4 isotope consists of two weakly stable deuteron [NP] clusters, each with their own distinct energy level, spin, magnetic moment, etc. Experimental tests are needed to confirm this fundamental model prediction. If confirmed, new physics mathematical description of shell structure of isotopes would follow.

6. The Brightsen Model predicts that forces "within" nucleon clusters (NPN,PNP,NP) are stronger that forces "between" such clusters within isotopes, a result of different combinations of the spin doublet and triplet clusters. It is predicted that research here would result in new measurable macroscopic properties of atomic nuclei including new fundamental force interactions.

7. The Brightsen Model predicts that the next "magic number" will be found at N = 172, Z = 106, A = 278 (Seaborgium-278). Experimental confirmation of this prediction would require a revised explanation of magic numbers in isotopes based on nucleon clusters as the fundamental building blocks of shell structure in atomic nuclei, as opposed to independent nucleons in an average field.

8. The Brightsen Model predicts that the large cross section of Boron-10 (as opposed to the small cross section of Boron-11) results from the presence of a stable and independent nucleon cluster structure [PNP], which coexists with two [NP] and one [NPN] clusters that maintain very small cross sections. Thus the vast majority of the cross section dynamics of Boron-10 is predicted by the Brightsen Model to derive from a strongly interacting [PNP] cluster. This four cluster formalism for Boron-10 (e.g., 1 PNP, 2 NP, 1 NPN) also correctly derives the I =3 spin experimentally observed.

# 2   Unsolved Problems in Mathematics

> Imagination, not invention, is the supreme master of art as of life.
> --Joseph Conrad

Most of the following problems come from one of the author's previous book [16]. Other problems come from recent collection of unsolved problems [45].

## 2.1 Maximum number of circles [16]

What is the maximum number of circles of radius 1, at most twice tangentials, which are contained in a circle of radius n? (Gamma 1/1986). This problem was generalized by Mihaly Bencze, who called for the maximum number of circles of radius $\Phi(n)$, at the most twice tangential, which are included into a circle of radius n, where $\varphi$ is function of n (Gamma 3/1986).

Also study a similar problem for circles of radius 1 contained in a given triangle (on Malfatti's problem), similar questions for spheres, cones, cylinders, regular pyramids, etc. More generally: for plane figures contained in a given planar figure, and in space too.

## 2.2 Consecutive sequence [16]

Given the consecutive sequence:
1,12,123,1234,12345,123456,1234567,12345678,123456789,13245678910,1234567891011, 123456789101112, 12345678910111213,…
How many primes are there among these numbers?
In general form, the Consecutive Sequence is considered in an arbitrary numeration base B.

Ref. Student Conference, University of Craiova, Department of Mathematics, April 1979, "Some problems in number theory," F. Smarandache.

## 2.3 Diophantine equation [16]

Conjecture:

Let k$\geq$2 be a positive integer. The Diophantine equation:

$$y = 2\, x_1 x_2 \ldots x_k + 1$$

has an infinite number of solutions in primes. (For example: 571=2*3*5*19+1, 691=2*3*5*23+1, 647=2*17*19+1, when k equals 4 and 3, respectively). (Gamma 2/1986)

## 2.4  Van der Waerden theorem  [16]

Expanding to infinity van der Waerden's theorem: Is it possible to partition N* into an infinity of arbitrary classes such that at least one class contains an arithmetic progression of $\ell$ terms ( $\ell \geq$3)?

Find the maximum $\ell$ having this property.

## 2.5  Differential equation with fractional power  [16]

Let $a \in Q \setminus \{-1, 0, 1\}$. Solve the equation:

$$x a^{\frac{1}{x}} + \frac{1}{x} a^x = 2a$$

Ref.: A generalization of the problem 0:123, *Gazeta Mathematica*, No. 3/1980, p.125.

## 2.6  Representation of odd number by primes [45]

Let k$\geq$3 and 1$\leq$s$<$k, where k and s integers. Then:

(i)     if k is odd, any odd integer can be expressed as the sum of k-s primes (first set) minus a sum of s primes (second set) [so that the primes of the first set are different to the primes of the second set].

(ii)    Is the conjecture true when all the k prime numbers are different?

(iii)   In how many ways can each odd integer be expressed as above?

(iv)    If k is even, any even integer can be expressed as the sum of k-s primes (first set) minus a sum of s primes (second set) [so that the primes of the first set are different to the primes of the second set]. $\rightarrow$ generalized Goldbach conjecture.

(v)     Is the conjecture true when all the k prime numbers are different?

(vi)    In how many ways can each even integer be expressed as above?

Ref.: [45] p. 10.

## 2.7   Magic square problem

A magic square is defined as an array of numbers which yields the same number when we add up all numbers in sequential order, either vertically, horizontally, or diagonally. The problem goes back to an ancient Japanese puzzle, known as the 'magic turtle' problem, which contains of magic square of rank-3 (3x3 square). Today you can solve magic square problems of rank-n with mathematical software like *Maple*.

(i)      Is there a limit of n for a rank-n magic square?

(ii)     Is there a general algorithm to compose any magic square with rank-n?

(iii)    Can we arrange a magic square which yields the same number either vertically, horizontally, or diagonally when we apply multiplication instead of addition?

(iv)     The same question as above, but for division?

(v)      Can we create a magic square of rank-n that consists of fractional numbers instead of integers?

(vi)     Is there any plausible linkage between a magic square and hadronic (quark) charges? Of what rank is it most likely to be link?

## 2.8   Palindromic numbers and iterations [45] p. 49

A number is said to be palindromic if it reads the same backwards and forwards. For example: 121, 1111, 34566543. The Pseudo-Smarandache function $Z(n)$ is defined for any $n \geq 1$ as the smallest integer m, such that n divides evenly into $1+2+\ldots+m$.

There are some palindromic numbers n such that $Z(n)$ is also palindromic $Z(909)=404$, $Z(2222)=1111$. Let $Z^k(n)= Z(Z(Z(\ldots(n)\ldots)))$, where the function Z is executed k times. $Z^0(n)=n$ by convention.

*Unsolved problem:* what is the largest value of m, such that for some n, $Z^k(n)$s is a palindrome for all values of k =0,1,2,....,m?

Conjecture (Ashbacher): there is no largest value of m, such that for some n, $Z^k(n)$s is a palindrome for all values of k =0,1,2,....,m.

### 2.9    Non-Euclidean geometry by giving up the fifth postulate [45] p. 52

- Definition:

A new type of geometry was constructed by F. Smarandache, which is simultaneously in a partially Euclidean and partially non-Euclidean space by replacing the Euclid's fifth postulate (axiom of parallels) with the following five-statement propositions:

a) There is at least one straight line and one point exterior to it in the space, for which only one line passes through the point and does not intersect the initial line [1 parallel];

b) There is at least one straight line and one point exterior to it in the space, for which only a finite number of lines $l_1, l_2, \ldots, l_k$ $(k \geq 2)$ pass through the point and do not intersect the initial line [2 or more (in a finite number) parallels];

c) There is at least one straight line and one point exterior to it in the space, for which any line that passes through the point intersects the initial line [0 parallel];

d) There is at least one straight line and one point exterior to it in the space, for which an infinite number of lines that pass through the point (but not all lines) do not intersect the initial line [an infinite number of parallels, but not all lines passing through the point];

e) There is at least one straight line and one point exterior to it in the space, for which any line that passes through the point does not intersect the initial line [an infinite number of parallels, all lines passing through the point].

- Problem:

(i) Can it be proved that the equation for a circle $x^2 + y^2 = r^2$ holds true in this non-Euclidean geometry? If yes, why?

(ii) How many nontrivial solutions are there provided the above equation holds true?

### 2.10  Smarandache Geometries and Degree of Negation in Geometries

We now present a more general class of geometries extracted from [3].

<u>Definition:</u>

An axiom is said Smarandachely denied if the axiom behaves in at least two different ways within the same space (i.e., validated and invalided, or only invalidated but in multiple distinct ways).

A Smarandache Geometry is a geometry which has at least one Smarandachely denied axiom (1969).

<u>Notations:</u>

Let's note any point, line, plane, space, triangle, etc. in a smarandacheian geometry by s-point, s-line,

s-plane, s-space, s-triangle respectively in order to distinguish them from other geometries.

<u>Applications:</u>

Why these hybrid geometries? Because in reality there does not exist isolated homogeneous spaces, but a mixture of them, interconnected, and each having a different structure.

The Smarandache geometries (SG) are becoming very important now since they combine many spaces into one, because our world is not formed by perfect homogeneous spaces as in pure mathematics, but by non-homogeneous spaces. Also, SG introduce the degree of negation in geometry for the first time [for example an axiom (or theorem, or lemma, or proposition) is denied in 40% of the space and accepted in 60% of the space], that's why they can become revolutionary in science and this thanks to the idea of partially denying and partially accepting of axioms/theorems/lemmas/propositions in a space (making multi-spaces, i.e. a space formed by combination of many different other spaces), similarly as in fuzzy logic (or in neutrosophic logic - the last one is a generalization of the fuzzy logic) the 'degree of truth' (i. e. for example 40% false and 60% true).

Smarandache geometries are starting to have applications in physics and engineering because of dealing with *non-homogeneous spaces.*

In the Euclidean geometry, also called parabolic geometry, the fifth Euclidean postulate that there is only one parallel to a given line passing through an exterior point, is kept or validated.

In the Bolyai-Gauss geometry, called *hyperbolic geometry*, this fifth Euclidean postulate is invalidated in the following way: there are infinitely many lines parallels to a given line passing through an exterior point.

While in the Riemannian geometry, called elliptic geometry, the fifth Euclidean postulate is also invalidated as follows: there is no parallel to a given line passing through an exterior point.

Thus, as a particular case, Euclidean, Bolyai-Gauss, and Riemannian geometries may be united altogether, in the same space, by some Smarandache geometries. These last geometries can be partially Euclidean and partially Non-Euclidean. Howard Iseri [3] constructed a model for this particular Smarandache geometry, where the Euclidean fifth postulate is replaced by different statements within the same space, i.e. one parallel, no parallel, infinitely many parallels but all lines passing through the given point, all lines passing through the given point are parallel.

Linfan Mao [4, 5] showed that SG are generalizations of Pseudo-Manifold Geometries, which in their turn are generalizations of Finsler Geometry, and which in its turn is a generalization of Riemann Geometry.

Let's consider Hilbert's 21 axioms of Euclidean geometry. If we Smarandachely deny one, two, three, and so on, up to 21 axioms respectively, then one gets:

$_{21}C_1 + {}_{21}C_2 + {}_{21}C_3 + _{...} + {}_{21}C_{21} = 2^{21} - 1 = 2,097,151$

Smarandache geometries, however the number is much higher because one axiom can be Smarandachely denied in multiple ways.

Similarly, if one Smarandachely denies the axioms of Projective Geometry, etc.

It seems that Smarandache Geometries are connected with the Theory of Relativity (because they include the Riemannian geometry in a subspace) and with the Parallel Universes (because they combine separate spaces into one space only) too.

A Smarandache manifold is an n-D manifold that supports a smarandacheian geometry.

Examples:

As a particular case one mentions Howard's Models [3] where a Smarandache manifold is a 2-D manifold formed by equilateral triangles such that around a vertex there are 5 (for elliptic), 6 (for Euclidean), and 7 (for hyper-bolic) triangles, two by two having in common a side.  Or, more general, an n-D manifold constructed from n-D submanifolds (which have in common two by two at most one m-D frontier, where m<n) that supports a Smarandache geometry.

<u>A Mode for a particular Smarandache Geometry:</u>

Let's consider an Euclidean plane ($\alpha$) and three non-collinear given points A, B, and C in it.  We define as s-points all usual Euclidean points and s-lines any Euclidean line that passes through one and only one of the points A, B, or C.  Thus the geometry formed is smarandacheian because two axioms are Smarandachely denied:

a) The axiom that through a point exterior to a given line there is only one parallel passing through it is now replaced by two statements: one parallel, and no parallel.

Examples:

Let's take the Euclidean line AB (which is not an s-line according to the definition because passes through two among the three given points A, B, C), and an s-line noted (c) that passes through s-point C and is parallel in the Euclidean sense to AB:

- through any s-point not lying on AB there is one s-parallel to (c).

- through any other s-point lying on the Euclidean line AB, there is no s-parallel to (c).

b)   And the axiom that through any two distinct points there exist one line passing through them is now replaced by: one s-line, and no s-line.

Examples:

Using the same notations:

- through any two distinct s-points not lying on Euclidean lines AB, BC, CA, there is one s-line passing through them;

- through any two distinct s-points lying on AB there is no s-line passing through them.

<u>Miscellanea:</u>

First International Conference on Smarandache Geometries was held, between May 3-5, 2003, at the Griffith University, Queensland, Australia, organized by Dr. J. Allen.

http://at.yorku.ca/cgi-bin/amca-calendar/public/display/conference_info/fabz54

There is a club too on "Smarandache Geometries" at and everybody is welcome.

http://clubs.yahoo.com/clubs/smarandachegeometries

For more information
see
http://www.gallup.unm.edu/~smarandache/geometries.htm

Questions:

Is there a general model for all Smarandache Geometries in such a way
that replacing some parameters one gets any of the desired particular SG?

## 2.11 Non-Archimedean triangle theorem

For the above non-Euclidean geometry by giving up the fifth postulate (2.9), is there a rule similar to the Archimedean triangle theorem?

Prove the new theorem in these multispaces.

## 2.12 The cubic Diophantine equation [16]

(i) The equation

$$x^3 + y^3 + z^3 = 1$$

has as solutions (9,10,-12) and (-6,-8,9). How many other nontrivial integer solutions are there?

(ii) As a generalization, how many solutions has the equation:

$$Ax^3 + By^3 + Cz^3 = D$$

where A,B,C,D are integers and (A,B,C)=1 ? (But $A_{1x1}{}^3 + A_{2x2}{}^3 + ... + A_{nxn}{}^3 = b$, where $A_1$, $A_2$, ..., $A_n$, B are integers and ($A_1$, $A_2$,..., $A_n$)=1?)

Ref.:

## 2.13  Multispaces and applications in physics

a) Multi-space unifies science (and other) fields; actually the whole universe is a multi-space. Our reality is so obviously formed by a union of many different spaces (i.e. a multi-space,
www.gallup.unm.edu/~smarandache/TRANSDIS.TXT).
Unfortunately there is not much theory behind "multi-space" (only some research done about Smarandache Geometries, that are a particular type of multi-space formed as unions of geometrical spaces). So, we can unite nano-scale space with our world scale and with cosmic scale, or we can unify the unorganic nanoscale with organic nanoscale, and so on.
The question is how to develop a multi-space theory? The connection among these heterogeneous spaces could be a problem.

b) Weyl and Kahler geometries are used in quantization somehow, but how should we use the Smarandache geometries
[www.gallup.unm.edu/~smarandache/geometries.htm], that look to be more general, in physics?

c) Can nanochips be used as additional human memories implemented in man's brain? So, like in science fiction, record a whole encyclopedia on an external nanodevice, implement it into a person's brain, and as a conclusion that person (doing a "search" in his artificial memory as done in Google.com) knowns everything from the encyclopedia?

# 3    Unsolved Problems in Astrophysics

> You see things; and you say "Why?"
> But I dream things that never were; and I say "Why not?" --George B. Shaw

In this chapter we discuss some problem in astrophysics and general celestial mechanics. While some of these questions may sound a bit awkward, perhaps they could trigger new ideas.

## 3.1  Some open problems in Celestial Mechanics

### 3.1.1    Photon speed

Is it possible to accelerate a photon (or another particle travelling at say 0.99 c), and thus to achieve a speed greater than c?

### 3.1.2    Flexible bridge

Will it be possible to construct a flexible bridge between two planets, and thus have terrestrial traffic between them?

What about the gravitational field of each planet. What about the gravitational field of each planet (a smooth escape of from gravitational field of one planet and a smooth entry to the other field.). Another difficulty would be the continuous motion of the planets.

Ref. :
- www.star-tech-inc.com/papers/ lse_iaf/LSE_IAF_04_Paper_Final.pdf
- www.theage.com.au/articles/2003/ 12/08/1070732145460.html?from=storyrhs
- www.csci-snc.com/GoingUp.htm

### 3.1.3    Splitting planet

Suppose we are able to dig around and thereby cut our planet into two separated parts.

    i.      In the first case, suppose the two halves are equal. Will these parts attract each other to reform into one planet again, or will they separate from each other?

    ii.     What if one part is significantly greater than the other?

### 3.1.4    Seeing the Moon

Why from the Moon is the Earth seen above, and from the Earth is the Moon seen above too?

(Let's consider a fixed point on the Earth; we are able to see the Moon from this point only when the Moon is above the point, because when the Moon is diametrically opposed it cannot be seen from that Earth point. Similarly when we consider a fixed point on the Moon, from where the Earth is visible from.).

### 3.1.5    Tunnel into the Earth (remember Jules Verne's story)

Let's consider a tunnel from one side to the other side of the Earth, and passing through the centre of the Earth.

    (i)     If one drops an object in the tunnel, will the object stop at the Centre of the Earth or will it oscillate like a pendulum across the centre, up and down, and after a while stop? Will the object then float at the centre?

    (ii)    If an elevator is located in the tunnel and experience free-falling, how much force would be necessary to push it up (especially from the centre of the Earth) to the other side of the Earth's surface? Isn't there any inertial force, from the fall force, that might push the elevator beyond the Earth's centre towards the other side?

    (iii)   Is there Coriolis effect experienced by passenger inside this elevator in story (ii)?

    (iv)   If we let an elevator to free fall inside the tunnel, will a passenger inside the elevator feel the force of gravity or not? (Recall that general relativity assumes *equivalence principle*, but it assumes the motion occurs on or around the surface of gravitational mass, but not inside the rotating mass.)

    (v)    Suppose that one end of the tunnel is located at the bottom of an ocean. Will water flow down into the tunnel only to the centre of the Earth, or will it flow lower nearer to the other side (to some-

how balance about the Earth's centre, the water masses from both sides of the Earth's centre), or will water flood out the first side?

(vi)  Repeat the above three questions for the case when the tunnel runs from one side to the other side of the Earth, but the tunnel doesn't pass through the Earth's centre. Would the midpoint of the tunnel play a similar role of as the Earth centre in the previous three questions?

(vii)  How will the Coriolis force influence this?

## 3.2 Some open problems in Astrophysics

### 3.2.1 Quasar & Pulsars

What is the dynamical origin of quasar and pulsars [17]? Are they ejected from the centres of galaxy(es)?

### 3.2.2 Graviton

Why hasn't the graviton been observed in experiment, supposing it to be the particle carrying gravitational field?

### 3.2.3 Precession of planets other than Mercury

Why didn't general relativity predict precession of planets other than Mercury? Can we predict these effects, and compare with observation?

### 3.2.4 Stargate

-  Provided general relativity is correct, then is it possible to create a real stargate (some theoreticians call it Einstein-Rosen Bridge) machine when one could travel from one place in this galaxy to another galaxy?

-  What are the conditions to keep the Stargate open and connecting the same places, so people could travel back to where they came from without violating causality principle?

### 3.2.5   Knot theory

A recent article described the presence of knots or sponge-like structure between galaxies:

"Surveys of the universe at its largest scales have found that galaxies are arranged into a sponge-like structure, with sheets and filaments of galaxies surrounding nearly empty voids. Places where these sheets and filaments intersect are sometimes called "knots," as they tend to have dense concentrations of galaxies that are merging."

- Provided these knots are real, and then is it possible to describe the structure of the Universe in the form of topological knot theory?
- What are its implications for astrophysics compared to more conventional general relativity predictions?
- Will darkmatter be required within the framework of knot theory?

Ref.: http://www.world-science.net/exclusives/060419_attractorfrm.htm

### 3.2.6   Solar neutrino [15]

- Where do the solar neutrinos come from? Do they have mass?
- Are solar neutrinos produced at the surface of the Sun? Or from its interior dynamics? [18]

### 3.2.7   Double star system

For a given double-star system of equal masses, we could ask some questions:
- Is it possible for this double-star system to rotate around a common circular orbit? If yes why? Is there a similarity to the bipolar electron?
(http://www.groupkos.com/mtwain/BipolarElectron144.jpg)
- If there is another object of equal mass passing near this double-stars system, will these objects form a triple-star system?

- Provided the observation data is given, then how could one calculate the velocities of the double-star system? Are there different possible interpretations?
- Is it possible to find an n-star system, with n <3?

### 3.2.8    Comets

- Where do comets come from? How were they produced in the past?
- Are comets composed of antimatter? Or water?

Ref.: http://matter-antimatter.com/

- Or, are these comets composed of plasma discharge?

"James McCanney, physicist and very active in astronomy, says it's plasma discharge and comets don't have water.  They can't both be right, unless comets have different compositions in their rock, etc., due to their original creation." http://www.jmccanneyscience.com/

### 3.2.9    Gamma ray burst

- Where does gamma ray bursts come from?
- What is the role of gamma ray bursts in the star formation processes?
- Are gamma ray bursts created when antimatter comets collide with stars?

Ref.: http://matter-antimatter.com/

### 3.2.10    Does Tifft's redshift quantization imply quantized distance between galaxies?

It is known to the astronomy society since 1980s that there is anomalous phenomenon called 'redshift quantization' introduced by Tifft. Some recent studies by Guthrie *et al*. seem to support this hypothesis, at least for 250 spiral galaxies. Now the question is:
   -    How to generalise this result for other astrophysics phenomena?
   -    Does Tifft's redshift quantization imply quantized distance between galaxies?

- If yes, then does it mean that quantized distance between galaxies imply formation of quantized vortices at astrophysics scale? (similar to cosmic string [125])
- Is there other astrophysics phenomena supporting the idea of 'quantized vortice' formation in galaxies?

We discuss this issue in this section.

In a recent paper by Moffat [115] it is shown that quantum phion condensate model with Gross-Pitaevskii equation yields an approximate fit to data corresponding to CMB spectrum, and it also yields a modified Newtonian acceleration law which is in good agreement with galaxy rotation curve data. It seems therefore interesting to extend further this hypothesis to explain quantization of redshift, as shown by Tifft *et al.* [116][120][121]. We also argue in other paper that this redshift quantization could be explained as signature of topological quantized vortices, which also agrees with Gross-Pitaevskiian description [117][119].

Furthermore, it is well-known that Gross-Pitaevskii equation could exhibit *topologically* non-trivial vortex solutions [118][119], which also corresponds to quantized vortices:

$$\oint p \cdot dr = N_v \, 2\pi\hbar \tag{1}$$

Therefore an implication of Gross-Pitaevskii equation [118] is that topologically quantized vortex could exhibit in astrophysical scale.

We start with standard definition of Hubble law [116]:

$$z = \frac{\delta\lambda}{\lambda} = \frac{Hr}{c} \tag{1a}$$

Or

$$r = \frac{c}{H} z \tag{1b}$$

Now we suppose that the major parts of redshift data could be explained via Doppler shift effect, therefore [116]:

$$z = \frac{\delta\lambda}{\lambda} = \frac{v}{c} \tag{1c}$$

In order to interpret Tifft's observation of quantized redshift corresponding to quantized velocity 36.6 km/sec and 72.2 km/sec, then we could write from equation (1b) and (1c):

$$\frac{\delta v}{c} = \delta z = \delta(\frac{\delta\lambda}{\lambda}) \tag{1d}$$

In this context we submit the viewpoint that the aforementioned proposition that topologically quantized vortex could exhibit in astrophysical scale (1) has been observed in the form of Tifft's redshift quantization [116][120]:

$$\delta r = \frac{c}{H}\delta z \qquad (2)$$

In other words, we submit the viewpoint that Tifft's observation of quantized redshift implies a quantized distance between galaxies [116][119], which could be expressed in the form:

$$r_n = r_o + n(\delta r) \qquad (3)$$

where n is integer (1,2,3,...) similar to quantum number. Because it can be shown using standard definition of Hubble law that redshift quantization implies quantized distance between galaxies in the same cluster, then one could say that this equation of quantized distance (3) is a result of *topological* quantized vortices (1) in astrophysical scale [119]; and it agrees with Gross-Pitaevskii (quantum phion condensate) description of CMB spectrum [115]. It is perhaps more interesting if we note here, that from (2) then we also get an equivalent expression of (3):

$$\frac{c}{H}z_n = \frac{c}{H}z_o + n(\frac{c}{H}\delta z) \qquad (4)$$

Or

$$z_n = z_o + n(\delta z) \qquad (5)$$

Or

$$z_n = z_o\left[1 + n(\frac{\delta z}{z_o})\right] \qquad (6)$$

In the meantime, it is interesting to note here similarity between equation (6) and (7), as observed in Fundamental Plane clusters and also from various quasars data [120][120a]:

$$z_{iQ} = z_f\left[N - 0.1M_N\right] \qquad (7)$$

Where $z_f$=0.62 is assumed to be a fundamental redshift constant, and N (=1,2,3...), and M is function of N.[6a] Here, the number M seems to play a role similar to second quantum number in quantum physics. [121].

Therefore it seems that we can interpret Tifft's redshift quantized as quantized distance between galaxies. Furthermore, one can rewrite equation (6) as follows:

$$z_n . z_o^{-1} = \left[ 1 + n \left( \frac{\delta z}{z_o} \right) \right] \qquad (6a)$$

And by using (1c), equation (6a) becomes:

$$\frac{d\lambda}{\lambda} . \frac{1}{z_o . n} = \left[ \frac{1}{n} + \left( \frac{dz}{z_o} \right) \right] \qquad (6b)$$

And using differentiation rules for logarithmic equations, yields:

$$\ln \lambda^{\frac{1}{z_o . n}} \approx \frac{k}{n} \ln(z) \qquad (6c)$$

Which implies that there exist neat (logarithmic) relationship between wavelength parameter and redshift quantization, and this relationship goes from microphysics up to macroscale phenomena, as described by Setterfield [121], i.e. it is likely that the *redshift quantization corresponds to the wavelength quantization of Quantum Mechanics.*[121]

It is more interesting if we note here that we can also explain astrophysical quantization using Weyl method in lieu of using generalised Schrödinger equation as Nottale did [118].

For instance, it can be shown that one can obtain Bohr-Sommerfeld type quantization rule from Weyl approach [126, p.12], which for kinetic plus potential energy will take the form:

$$2 \pi N \hbar = \sum_{j=0}^{\infty} \hbar^j S_j(E) \qquad (8)$$

Which can be solved by expressing $E = \sum \hbar^k E_k$ as power series in $\hbar$. [126]. Now Bohr-Sommerfeld quantization rule [126] could be rewritten as follows:

$$\oint p \cdot dr = N_v 2 \pi \hbar = \sum_{j=0}^{\infty} \hbar^j S_j(E) \qquad (9)$$

Or if we consider quantum Hall effect [119a]:

$$\Phi = q \oint A \cdot dr + m \oint \Omega \times r \cdot dr = \iint B \cdot dS = N_v 2 \pi \hbar \qquad (10)$$

then equation (10) can be used instead of equation (9), which yields:

$$\Phi = q \oint A \cdot dr + m \oint \Omega \times r \cdot dr = \iint B \cdot dS = \sum_{j=0}^{\infty} \hbar^j S_j(E) \qquad (11)$$

The above method is known as 'graph kinematic' [127] or Weyl-Moyal's deformation quantization [128].

### 3.2.11 Does Pioneer anomaly imply a modified gravitation theory?

There is a known anomalous observation from Pioneer spacecraft which baffles physicists and general relativists since 1990s [129]-[137]. This anomaly –called Pioneer anomaly— is essentially a slight departure from Newtonian acceleration at the order of ~8.74x10$^{-10}$ m/sec$^2$ which is observed since Pioneer spacecraft entered the orbit of Jupiter. [129][130]

Apart from some other plausible interpretations, there are some 'obvious' questions which apparently deserve further considerations before one considers other 'exotic' theories:

- Does Pioneer anomaly imply a modified gravity theory? [131][132] According to Brownstein & Moffat, the Pioneer anomalous acceleration directed toward the center of the Sun could be written as follows [131]:

$$a_p = -\frac{\delta G(r) M_\oplus}{r^2} \qquad (12)$$

Where

$$\delta G(r) = G_o \alpha(r) \left[ 1 - e^{-\frac{r}{\lambda(r)}} \left( 1 + \frac{r}{\lambda(r)} \right) \right] \qquad (13)$$

- Another remaining question for the above proposition of modified Newtonian acceleration [131][132]: Why was the anomalous acceleration not detected before Pioneer spacecraft reached Jupiter orbit?
- Does Pioneer anomaly exhibit other higher-order gravity effects which could be detected?
- Does Pioneer anomaly imply that Pioneer spacecraft orbit is affected by Jupiter gravitational fields?
- Does Pioneer anomaly imply that there is significant difference between gravitational field of inner planets and Jovian planets, corresponding to 'two-fluid' model of superfluidity?

We could expect that Pioneer anomaly can be explained within five or ten years, and it is likely that it will turn out to be higher-order gravitational effect beyond standard Newtonian acceleration law.

# 4 Unsolved Problems in Geophysics

> There is always an easy solution to every human problem—neat, plausible and wrong. --H. L. Mencken

## 4.1 Introduction

As strange as it may seem, we know very little about the actual dynamics of geophysical systems including our Earth.

In this regards, J. Bowles has argued: "Earth processes are the result of the Newtonian gravitational-forces that accelerate the earth (as with all planets and their moons) into curvilinear - orbital motion. It is Newton's 1st Law: {a body in motion will follow a straight line unless acted upon by an external force.} What this means is that at every infinitesimal moment - the straight line motion that would be the earth's, is changed {by solar gravitational forces} into a curvilinear orbital motion. The stresses induced by this acceleration, and the fact that we're in simultaneous rotation .. are what cause all geo-physical activity .. including the generation of immense currents that eventually strike as earthquake." <jimbow1@mindspring.com>

From this viewpoint, we think that we could ask some interesting questions.

## 4.2 Some new questions

### 4.2.1 Newtonian dynamics

- What is the governing Newtonian dynamical equation of the geophysical processes?
- Can we put this governing dynamical equation into a Cauchy equation form? What is the numerical solution?
- What is the role of the Coriolis force in this dynamical equation?

Ref.: Nash, J., "Le probleme de Cauchy pour les equations differentielles d'un fluide general," *Bull. De la S.M.F.*, tome 90 (1962), p. 487-497. (www.numdam.org)

### 4.2.2     Geophysical processes and solar radiation:

-   Is there any link between global warming and earthquakes?
-   Or is there any linkage between fluctuation of solar radiation and earthquakes?
-   Is there a definitive correlation between fluctuation of solar radiation and Earth climatic changes? (Some articles describing the Numerical climatic models, predict that a change in solar output of only 1 percent per century would alter the Earth's average temperature by between 0.5-1.0 degrees. See Ref.)
-   Is there a mathematical model to describe this theoretical correlation between Solar radiation and Earth climatic changes?
-   How good is the mathematical model (if it exists) to predict new climatic changes like large volcano activities or hurricanes etc.?

Ref.: (i) www.physicalgeography.net/fundamentals/7y.html; (ii) en.wikipedia.org/wiki/Climate_change

### 4.2.3     Geophysical process and Chandler polar wobble:

-   Is there a plausible linkage between Newtonian dynamics of the geophysical processes inside the Earth and Chandler polar wobble?
-   Is this Chandler polar wobble caused by Earth's internal geodynamics, including rotational inertia and quadrupole moments?
-   Or are there external forces that may affect this polar wobble (such as the motion of the Sun, or fluctuation of Solar radiation etc.)?
-   What is the role of the Coriolis force on the dynamics of Chandler polar wobble?
-   Can we describe a new prediction of Chandler wobble, at least for a few months ahead?
-   Is there any plausible linkage between Chandler wobble and some earthquakes, because of Earth instability?

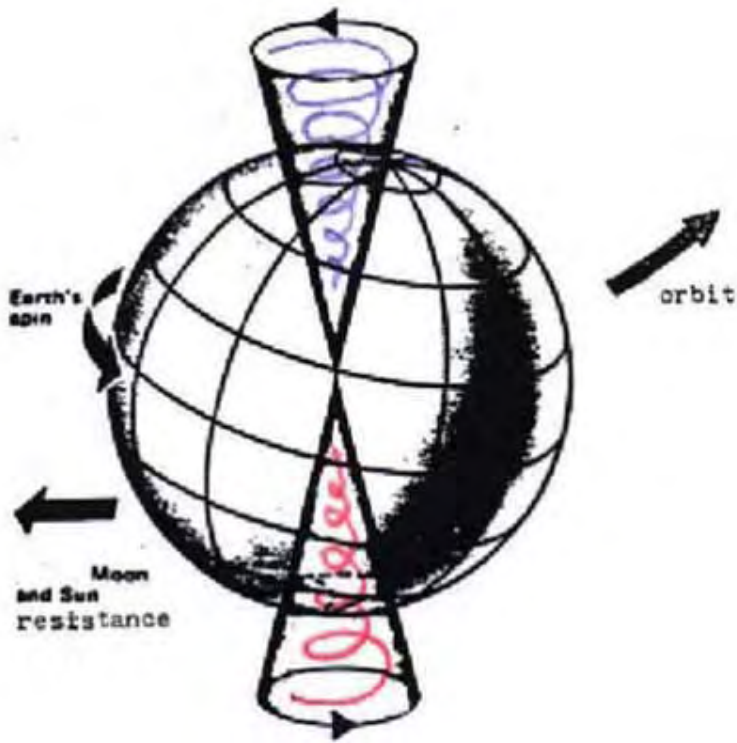

Figure 4. Illustration of Chandler polar wobble

Ref.: (i) http://www.newscientist.com/article/mg18925391.500.html,

(ii) www.tmgnow.com/repository/global/outside_forces.html,

(iii) www.huttoncommentaries.com/Summaries.htm

**Other questions:**

- Provided the Earth can be viewed as an electric-capacitor, then is it possible to stabilise the Chandler polar wobble via introducing an artificial electromagnetic field?
- Story: The Earth is an electrified body, moving in plasma. We who stand on its surface are seldom aware of its electrical properties. That's because we live in balance with the Earth's electric field. Like a high-tension wire, our Earth produces hums and crackles as it responds to surges of power in the electric currents of space. Perhaps the most obvious sparks are the auroras. (see: http://www.gnn.tv/threads/8218/VIOLENT_SOLAR_STORMS_MAKE_EARTHQUAKES_30_TIMES_MORE_LIKELY).

Background theory:

We begin with the proposition that planetary systems, including the Earth, can be modelled as a rotating superconductor, therefore we can use the Gross-Pitaevskii equation.

It is well-known that we can find an Euler-type fluid from the Gross-Pitaevskii equation. First, we could rewrite equation (1) in a more common form [97][98], by neglecting Planck's constant:

$$-i.\partial \psi / \partial t - \Delta \psi + (1 - \kappa|\psi|^2)\psi = 0.$$  (1)

Now, using the Madelung transformation [7, p. 7],

$$\psi = \sqrt{\rho}\exp(i\varphi),$$  (2)

one can write equation (1) in the variables ($\rho, \vec{\upsilon} := 2\nabla\varphi$) [7]:

$$\frac{\partial \rho}{\partial t} + div(\rho\vec{\upsilon}) = 0,$$  (3)

$$\frac{\partial \vec{\upsilon}}{\partial t} + \vec{\upsilon}.\nabla\vec{\upsilon} + \nabla(2\rho) = -\nabla(\frac{|\nabla\rho|^2}{2\rho^2} - \frac{\Delta\rho}{\rho}).$$  (4)

Neglecting the last term on the right hand side, often called 'quantum pressure', this system reduces to the Euler equations for compressible ideal fluids, with speed $\vec{\upsilon}$ and pressure $\rho^2$.[98] Interestingly, Nottale has derived a similar Euler fluid from Schrödinger equation [99], which again seems to support a plausible linkage between GPE and Schrödinger equation via the Madelung transformation.

As an alternative to this known Madelung transform (2), one could use quaternion Madelung transform which leads to a quaternionic Schrödinger equation [100]. First, we write the group elements as [100]:

$$g = e^{iu}e^{j\phi}e^{kv},$$  (5)

where i,j,k are the quaternionic elements satisfying [101]:

$$ii = jj = kk = ijk = -1,$$  (6)

and $\mu, \nu, \phi$ are real angles [100]. Substitution of this form for g into equation (5) yields:

$$(\nabla g) g^{-1} = iu + jv + kw, \qquad (7)$$

where $\mathbf{u,v,w}$ are quaternionic generalization of Madelung variables [100].

Now, instead of one equation of motion for velocity (4), one gets equations for each $\mathbf{u,v,w}$ [100, p.6]:

$$\frac{\partial u}{\partial t} = -\nabla V' - \frac{\nabla}{2}\left(u^2 + v^2 + w^2\right) - \frac{v}{\rho}\nabla.(\rho v) - \frac{w}{\rho}\nabla.(\rho w), \quad (8a)$$

$$\frac{\partial v}{\partial t} = \frac{\nabla}{2}\left(\frac{\nabla.(\rho w)}{\rho}\right) - 2w[V'(x) + \frac{1}{2}\left(u^2 + v^2 + w^2\right)] + \frac{u}{\rho}\nabla.(\rho w), \qquad (8b)$$

$$\frac{\partial w}{\partial t} = \frac{\nabla}{2}\left(\frac{\nabla.(\rho v)}{\rho}\right) - 2v[V'(x) + \frac{1}{2}\left(u^2 + v^2 + w^2\right)] + \frac{u}{\rho}\nabla.(\rho v), \qquad (8c)$$

Alternatively, our quaternion transform could be written in the form of quaternion field [101, p.4]:

$$\Phi = i\Psi_1 + j\Psi_2 + k\Psi_3 + \Psi_4 \qquad (9)$$

which has a neat link to the original form of the Maxwell's equations [101].

A quaternionic description may also be useful for describing the Euler equation of motion for rigid body dynamics, because it could describe rotations directly [102]. We shall discuss its basic method in the following section based on [102]. While we do not introduce new equations in the following sections, the fact that we can use quaternions to describe solid dynamics along with Maxwell's electromagnetic fields lead us to a conjecture that one could expect to stabilize solid dynamics rotation when it is necessary using electromagnetic fields. This is our motivation in this section.

In keeping track of the motion of rigid body, one needs to be able to give its orientation in space at each instant. In accordance with [102] we shall use **x** to denote a vector, $\underline{\mathbf{x}}$ to denote its coordinates in space, and $\underline{\mathbf{X}}$ to denote its coordinates in the body. A vector X in the body will have coordinates in space:

$$\underline{x} = Q(t)\underline{X} \qquad (10)$$

where Q(t) is an orthogonal transformation. For the time evolution of $\underline{\mathbf{X}}$ and $\underline{\mathbf{x}}$, one get [102]:

$$\frac{d\underline{x}}{dt} = \frac{dQ}{dt}\underline{X} + Q\frac{d\underline{X}}{dt} \qquad (11)$$

And because $\underline{\mathbf{X}} = Q^T\underline{\mathbf{x}}$, we can write equation (11) as [102]:

$$\frac{d\underline{x}}{dt} = \frac{dQ}{dt}Q^T\underline{x} + Q\frac{d\underline{X}}{dt}$$

(12)

In the meantime, Newton laws give:

$$\frac{d\underline{m}}{dt} = \underline{t} ,$$

(13)

where $\underline{m} = A\underline{\omega}$, and $\underline{\mathbf{t}}$ is some *externally applied torque* [102].

In the body, this becomes [102]:

$$\frac{d\underline{m}}{dt} = \underline{t} = Q\left(\left[\underline{\Omega}, \underline{M}\right] + \frac{d\underline{M}}{dt}\right),$$

(14)

or, rewritten entirely in terms of body variables [102]:

$$\frac{d\underline{M}}{dt} = \left[\underline{M}, \underline{\Omega}\right] + \underline{T} ,$$

(15)

with $\underline{T} = Q^T\underline{t}$. These are the Euler equations of motion for a rigid body. In the absence of an external torque (such as a tidal wave effect, or a nearby/flyby asteroid), there are two integrals: the energy T and the total angular momentum M=$\underline{\mathbf{M}}.\underline{\mathbf{M}}$, i.e.

$$T = \frac{1}{2}\left(\frac{M_1^2}{T_1} + \frac{M_2^2}{T_2} + \frac{M_3^2}{T_3}\right)$$

(16)

and

$$M^2 = M_1^2 + M_2^2 + M_3^2$$

(17)

What is more interesting here is that in the context of the Ginzburg-Landau equation, one could write Newton's second law in the quaternionic form of equation (13):

$$\frac{d\underline{m}}{dt} = e(E + v \times B) + \nabla(e\varphi + \frac{1}{2}mv^2)$$

(18)

We submit the viewpoint that this equation (18) is quite useful for analysing interactions between rotational rigid body dynamics and electromagnetic torque fields such as the '*Birkeland effect.*' [103]. It is recommended, however, that experimental support be found for this proposition.

Ref:
http://www.theworld.com/~sweetser/quaternions/gravity/unified_field/unified_field.html

### 4.2.4     General relativistic

- Can we translate Newtonian dynamics for geophysical processes into a general relativistic expression?
- What are the new predictions from this new expression? Can we find an experimental test result supporting this prediction?
- With regard to the condensed-matter analogue of general relativity (Volovik et al.), does it mean that we can consider geophysical processes (like Earth's), as composed of condensed-matter superconductivity?
- If yes, then what are implications of this superconductivity model of Earth? Can we find experimental verification of these implications?

# 5    Unsolved Problems in Sorites Quantum Paradox and Smarandache Class of Paradoxes

*We used to think that if we knew one, we knew two, because one and one are two. We are finding that we must learn a great deal more about `and`.—Sir A. Eddington*

There can be generated many paradoxes or quasi-paradoxes that may occur from the combination of quantum and non-quantum worlds in physics. Even the passage from the micro-cosmos to the macro-cosmos, and conversely, can generate unsolved questions or counter-intuitive ideas. We define a quasi-paradox as a statement which has a *prima facie* self-contradictory essence or an explicit contradiction, but which is not completely proven as a paradox.

We present herein four elementary quantum quasi-paradoxes and their corresponding quantum Sorites paradoxes, which form a class of quantum quasi-paradoxes.

## 5.1    Introduction [38]-[51]

According to Dictionary of Mathematics (Borowski & Borwein, 1991), the **paradox** is "an apparently absurd or self-contradictory statement for which there is *prima facie* support, or an explicit contradiction derived from apparently unexceptionable premises". Some paradoxes require the revision of their intuitive conception (Russell's paradox, Cantor's paradox), others depend on the inadmissibility of their description (Grelling's paradox), others show counter-intuitive features of formal theories (Material Implication paradox, Skolem Paradox), others are self-contradictory [Smarandache Paradox: "All is <A> the <Non-A> too!", where <A> is an attribute and <Non-A> its opposite; for example "All is possible the impossible too!". Paradoxes are normally true and false at the same time.

The **Sorites paradoxes** are associated with Eubulides of Miletus (fourth century B.C.) and they state that there is not a clear frontier between visible and invisible matter, the determinate and indeterminate principle, stable and unstable matter, long time living and short time living matter. Generally, between <A> and <Non-A> there is no clear distinction, no exact frontier. Where does <A> really end and <Non-A> begin?  One extends Zadeh's "fuzzy set" concept to the "neutrosophic set" concept.

Let's now introduce the notion of quasi-paradox: A **quasi-paradox** is a statement which has a *prima facia* self-contradictory essence or an explicit contradiction, but which is not completely proven as a paradox.  A quasi-paradox is an *informal* contradictory statement, while a paradox is a *formal* contradictory statement.

Some of the below quantum quasi-paradoxes below can be proven to be real quantum paradoxes.

## 5.2   Quantum Paradox and Quantum Sorites Paradox [38]-[51]

It is the interaction of the quantum world with the "environment", associated with the large-scale world, which is thought to cause wave function collapse. For this reason we do not perceive the quantum behavior of every particle inside Schrödinger's cat; the presence of such an "environment" (the body of the cat) is thought to cause the cat to be seen to be either dead or alive, even though it may be poisoned as a result of a quantum phenomenon.

The following quasi-paradoxes and Sorites paradoxes are based on the antinomies: visible/invisible, determinatet/indeterminate, stable/unstable, long time living/short time living, as well as on the fact that there is not a clear separation between these pairs of antinomies.

5.2.1.1. **Invisible Quasi-Paradox:** Our visible world is composed of a totality of invisible particles

5.2.1.2. **Invisible Sorites Paradox:** There is not a clear frontier between visible matter and invisible matter.

    a)    An invisible particle does not form a visible object, nor do two invisible particles, three invisible particles, etc. However, at some point, the collection of invisible particles becomes large enough to form a visible object, but there is apparently no definite point where this occurs.

    b)    A similar paradox is developed in an opposite direction. It is always possible to remove a particle from an object in such a way that what is left is still a visible object. However, repeating and repeating this process, at some point, the visible object

is decomposed so that the remaining part becomes invisible, but there is no definite point where this occurs.

5.2.2.1. This **Uncertainty Quasi-Paradox**: Bulk matter, which is to some degree subject to the 'determinate principle', is formed by a totality of elementary particles, which satisfy Heisenberg's 'indeterminacy principle'.

5.2.2.2. **Uncertainty Sorites Paradox:** Similarly, there is no a clear frontier between matter subject to the 'determinate principle' and matter subject to the 'indeterminate principle'.

5.2.3.1. **Unstable Quasi-Paradox**: 'Stable' matter is formed by 'unstable' elementary particles (elementary particles decay when free).

5.2.3.2. **Unstable Sorites Paradox:** Similarly, there is not a clear frontier between the 'stable matter' and the 'unstable matter'.

5.2.4.1. **Short-Time-Living Quasi-Paradox**: 'Long-time-living' matter is formed by very 'short-time-living' elementary particles.

5.2.4.2. **Short-Time-Living Sorites Paradox:** Similarly, there is not a clear frontier between the 'long-time-living' matter and the 'short-time-living' matter.

Additional quantum quasi-paradoxes and paradoxes can be designed, all of them forming a class of Smarandache quantum quasi-paradoxes." (Dr. M. Khoshnevisan, Griffith University).

## 5.3   Smarandache Class of Paradoxes [19]-[37]

Standard definition of Smarandache Class of Paradoxes are as follows:
Let <A> be an attribute, and <Non-A> its negation. Then:

Paradox 1. ALL IS <A>, THE <Non-A> TOO.
Examples:
E11: All is possible, the impossible too.
E12: All are present, the absents too.
E13: All is finite, the infinite too.

## 5.4   Paradox

Paradox 1. ALL IS <A>, THE <Non-A> TOO.
Examples:

E11: All is possible, the impossible too.
E12: All are present, the absent too.
E13: All is finite, the infinite too.

Paradox 2. ALL IS <Non-A>, THE <A> TOO.
Examples:
E21: All is impossible, the possible too.
E22: All are absent, the presents too.
E23: All is infinite, the finite too.

Paradox 3. NOTHING IS <A>, NOT EVEN <A>.
Examples:
E31: Nothing is perfect, not even the perfect.
E32: Nothing is absolute, not even the absolute.
E33: Nothing is finite, not even the finite.

Remark: The three kinds of paradoxes are equivalent. They are called: The Smarandache Class of Paradoxes.

## 5.5 Generalization

We can put the above Smarandache Class of Paradoxes in more general statement:

Paradox: ALL (Verb) <A>, THE <Non-A> TOO

(<The Generalized Smarandache Class of Paradoxes>)

Replacing <A> by an attribute, we find a paradox.

Let's analyse the first one (E11):
< All is possible, the impossible too. >

If this sentence is true, then we get that <the impossible is possible too>, which is a contradiction; therefore the sentence is false. (Object Language).

But the sentence may be true, because <All is possible> involves that <the impossible is possible>, i.e.< it's possible to have impossible things>, which is correct. (Meta-Language).

Of course, this method leads to some unsuccessful paradoxes, but the proposed method yields other beautiful results. Consider this pun which reminds you of Einstein:

*All is relative, the (theory of) relativity too!*

(or if you try with Godel's incompleteness theorem, we can put forth the similar argument:

*This statement is unprovable.*

Godel's first theorem too is unprovable within its own numerization method.)

So:

1. The shortest way between two points is the meandering way!

2. The unexplainable is, however, explained by the word: "unexplainable"! (This is another version of the Grelling paradox: Is the word 'heterological', heterological?)

# 6    Origin of Spin: Paradox of the classical Beth experiment

*Creative effort must always call for guessing; and even the best guessing cannot avoid error. –A. Osborn*

*Contribution of R. Khrapko.*

A celebrated Beth's experiment contradicts the angular momentum conservation law in the frame of Maxwell electrodynamics because Beth's birefringent plate experienced a torque without an angular momentum flux in the surrounding space. However, this paradox can be removed by introducing a classical spin tensor.

Questions:

(I)     Is there (classical) Maxwell electromagnetic description of quantum spin?

(II)     Is there plausible linkage between spin, Aharonov effect and hidden-variable interpretation of Quantum Mechanics?

Ref.: - HY. Cui, arXiv.org:physics/0408025, arxiv.org/pdf/hep-th/0110128

## 6.1   Angular momentum of circularly polarized light

It has been known for long time that, on the basis of either the wave theory [1, 2] or the quantum theory (by assigning an angular momentum of $\pm h/2\pi$ to a photon), a circularly polarized light should exert a torque on a doubly refracting plate which changes the state of polarization of the light, or on a medium which (maybe partly) absorbs the light.

R. A. Beth explained [3] that the moment of force or torque exerted on a doubly refracting medium by a light wave passing through it arises from the fact that the dielectric constant $\varepsilon$ is a tensor. Consequently the electric intensity $\mathbf{E}$ is, in general, not parallel to the electric polarization $\mathbf{P}$ in the medium. The torque per unit volume produced by the action of the electric field on the polarization of the medium is

$$\tau / V = \mathbf{P} \times \mathbf{E} \qquad (1)$$

R. Feynman repeated this explanation [4]. We quote him from [4] with insignificant abridgements.

"If we have a beam of light containing a large number of photons all cir-cularly polarized the same way, it will carry angular momentum. If the total energy carried by the beam in a certain time is $W$, then there are $N = 2\pi W / h\omega$ photons. Each one carries the angular momentum $h/2\pi$, so there is a total angular momentum of

$$J_z = Nh / 2\pi = W / \omega .\qquad (2)$$

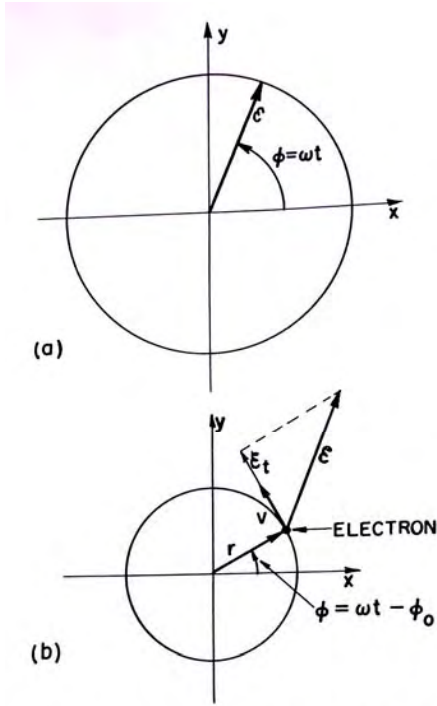

Can we provide classically that light which is right circularly polarized carries an angular momentum and energy in proportion $1/\omega$? Here we have a case where we can go from the quantum things to the classical things. Remember what right circularly polarized light is, classically. It's described by an electric field so that the electric vector $\mathbf{E}$ goes in a circle – as drawn in Fig. 17-5(a). Now suppose that such a light shines on a wall which is going to absorb it – or at least some of it – and consider an atom in the wall according to the classical physics. We'll suppose that the atom is isotropic, so the result is that the electron moves in a circle, as shown in Fig. 17-5(b). The electron is displaced at some displacement $\mathbf{r}$ from its equilibrium position at the origin and goes around with some phase lag with respect to the vector $\mathbf{E}$. The relation between $\mathbf{E}$ and $\mathbf{r}$ might be as shown in Fig. 17-5(b). As time goes on, the electric field rotates and the displacement rotates with the same frequency, so their relative orientation stays the same. Now let's look at the work being done on this electron. The rate that energy is being put into this electron is $v$, its velocity, times the component of $\mathbf{E}$ parallel to the velocity:

$$dW / dt = eE_t v \quad .$$
(3)

**Fig. 17–5. (a) The electric field ε in a circularly polarized light wave. (b) The motion of an electron being driven by the circularly polarized light.**

But look, there is angular momentum being poured into this electron, because there is always a torque about the origin. The torque is $eE_t r$ which must be equal to the rate of change of angular momentum $dJ_z / dt$:

$$dJ_z / dt = eE_t r \quad .$$
(4)

Remembering that $v = \omega r$, we have that

$$dJ_z / dW = 1 / \omega \quad .$$
(5)

Therefore, if we integrate the total angular momentum which is absorbed, it is proportional to the total energy – the constant of proportionality being $1 / \omega$."

Thus *Beth's and Feynman's reasoning prove that a circularly polarized plane wave carries angular momentum whose density is proportional to the energy density*. Unfortunately, the authors did not give an expression for the angular momentum flux density through the field quantities. At the same time, the scientific community refutes *the existence of angular momentum of plane waves.*

Heitler wrote [5]:

"In Maxwell's theory the Poynting vector $\mathbf{E} \times \mathbf{H}$ (divided by $c^2$) is interpreted as the density of momentum of the field. We can then also define an angular momentum relative to a given point $O$ or to a given axis,

$$\mathbf{J} = \int_V \mathbf{r} \times (\mathbf{E} \times \mathbf{H}) dV$$
(6)

where $\mathbf{r}$ is the distance from $O$ and $V$ is the volume of a transverse slice of a beam [$c = 1$ in this paper].

A plane wave traveling in the z-direction and with infinite extension in the xy-directions can have no angular momentum about the z-axis, because $\mathbf{E} \times \mathbf{H}$ is in the z-direction and $[\mathbf{r} \times (\mathbf{E} \times \mathbf{H})]_z = 0$. However, this is no longer the case for a wave with finite extension in the xy-plane. Consider a cylindrical wave with its axis in the z-direction and traveling in this direction. At the wall of the cylinder $r = R$, say, we let the amplitude drop to zero. It can be shown that the wall of such a wave packet gives a finite contribution to $J_z$."

Ohanian wrote [6]:

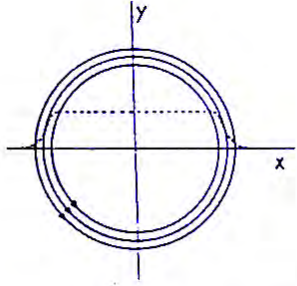

Fig. 1. This pattern of circular flow lines represents the time-average energy flow, or the momentum density, in a circularly polarized electromagnetic wave packet. On a given wave front, say $z = 0$, the fields are assumed to be constant within a circular area and to decrease to zero outside of this area (the dashed line gives the field amplitude as a function of radius). The energy flow has been calculated from an approximate solution of Maxwell's equations. The picture only shows the flow in the transverse directions. The flow in the longitudinal direction is much larger; the net flow is helical.

502    Am. J. Phys., Vol. 54, No. 6, June 1986

"In an infinite plane wave, the $\mathbf{E}$ and $\mathbf{H}$ fields are everywhere perpendicular to the wave vector and the energy flow is everywhere parallel to the wave vector. However, in a wave of finite transverse extent, the $\mathbf{E}$ and $\mathbf{H}$ fields have a component parallel to the wave vector (the field lines are closed loops) and the energy flow has components perpendicular to the wave vector. For instance, Fig. 1 shows the time-average transverse energy flow in a circularly polarized wave propagating in the z-direction; the wave has a finite extent in the x and y directions. The circulating energy flow in the wave implies the existence of angular momentum, whose direction is along the direction of propagation."

Simmonds and Guttman wrote [7]:

"The electric and magnetic field of a cylindrical beam can have a nonzero z-component only within the 'skin' region of the wave. Having z-component within this region implies the possibility of a nonzero z-component of angular momentum within this region. Since the wave is identically zero outside the skin and constant inside the skin region, the skin region is the only one in which the z-component of angular momentum does not vanish.

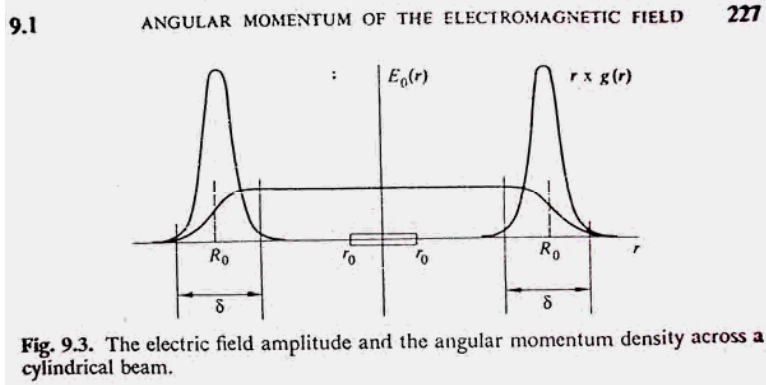

**9.1**          ANGULAR MOMENTUM OF THE ELECTROMAGNETIC FIELD          **227**

**Fig. 9.3.** The electric field amplitude and the angular momentum density **across a** cylindrical beam.

In Fig. 9.3 we plot an acceptable function $E_0(x, y) = E_0(r)$. We have explicitly made $E_0$ constant over a large central region of the wave and confined the variation of the function from this constant value to zero to lie within a 'skin' of thickness $\delta$ which lies a distance $R_0$ from the axis."

A calculation of the angular momentum **J**, according to eqn. (6), requires an explicit expression for the beam. We use the Jackson's expressions [8] with $k = \omega$ here,

$$\mathbf{E} = \exp(i\omega z - i\omega t)[\mathbf{x} + i\mathbf{y} + \frac{1}{\omega}\mathbf{z}(i\partial_x - \partial_y)] E_0(x, y)$$
$$\mathbf{H} = -i\mathbf{E}.$$

(7)

Transform the integrand of $J_z$ from eqn. (6),

$$\Re[x(\mathbf{E} \times \mathbf{H}^*)_y - y(\mathbf{E} \times \mathbf{H}^*)_x]/2$$

$$= \Re[x(E_z H_x^* - E_x H_z^*) - y(E_y H_z^* - E_z H_y^*)]/2$$

$$= \Re[x(E_z i E_x^* - E_x i E_z^*) + y(E_z i E_y^* - E_y i E_z^*)]/2$$

$$= -\Im[x E_z E_x^* + y E_z E_y^*]$$

$$= -\Im[x E_0(i\partial_x - \partial_y)E_0 + y(-iE_0)(i\partial_x - \partial_y)E_0]/\omega$$

$$= -(x\partial_x + y\partial_y)E_0^2/2\omega. \tag{8}$$

Substituting (8) into (6) and integrating by part yields

$$J_z = \int E_0^2 dV / \omega \ . \tag{9}$$

The power $\mathsf{P}$ of the beam is,

$$\mathsf{P} = \int_a (\mathbf{E} \times \mathbf{H})_z \, da = \int_a \Re(E_x H_y^* - E_y H_x^*) da / 2 = \int_a E_0^2 \, da \tag{10}$$

where $a$ is the cross-section area of the beam. If $l$ is a length of the transverse slice of the beam, i.e. $V = la$, the energy of the slice is

$$W = \int E_0^2 dV \tag{11}$$

because $c = 1$. So the relation between the total angular momentum $J_z$ and the total energy $W$,

$$J_z / W = 1 / \omega, \tag{12}$$

is the same in Beth – Feynman paradigm and in the scientific community paradigm. However, the distribution of the angular momentum is different. According to Beth and Feynman, the angular momentum density is proportional to energy density in a beam or in a plane wave, but, according to the community, the angular momentum is located near the wall of the beam and is absent in the plane wave.

In connection with this difference an important question was raised at the V. L. Ginsburg, *Moscow Physical Seminar* in the spring of 1999. The question concerned absorption of a circularly polarized light by a round flat target, which is divided concentrically into an inner disc and a closely fitting outer annulus [8].

If the target absorbs a circularly polarized beam, the annulus absorbs the wall or 'skin' of the beam, which carries the angular momentum, according to the community, and the disc absorbs the body of the beam, which has no angular momentum. Since the Poynting vector is perpendicular to the disc, an infinitesimal force

$$dF^i = T^{ij} da_j$$

(13)

acting on a surface element $da_j$ of the disc is also perpendicular to the disc ($T^{ij}$ is the Maxwell stress tensor). So, the disc does not perceive a torque when the target absorbs a circularly polarized beam. There are no ponderomotive forces, which are capable of twisting the disc. Tangential forces act only on the annulus.

But it is clear that *in reality the disc does perceive a torque from the wave, since the disc gets angular momentum, according to Beth – Feynman.* The disc will be twisted in contradiction with the common paradigm.

Allen and Padgett [9] attempted to explain the torque acting on the disc within the scope of the paradigm. They mentally decomposed the beam into three beams: the inner beam, the annulus beam, and the remainder. They wrote, *"Any form of aperture introduces an intensity gradient, so a field component is induced in the propagation direction and the dilemma is potentially resolved."*

Alas! A small clearance between the inner disc and outer annulus does not aperture a beam and does not induce longitudinal field components. The *imaginary decomposition of a wave is not capable of generating longitudinal field components and, correspondingly, transverse momentum and torque acting on the disc.* Maxwell stress tensor **cannot supply the disc with a torque. According to the Maxwell theory, the disc absorbs energy and experiences normal pressure only.**

Thus the mental experiment shows a weakness of the common paradigm. Does Beth's experiment confirm the formula (6)?

### 6.2  Beth experiment result is a puzzle

The classical Beth's experiment [3] was made 70 years ago. A beam of circularly polarized light exerts a torque on a doubly refracting plate, which changes the state of polarization of the light beam. The apparatus used involves a torsional pendulum with about a ten minute period consisting of a round quartz half-wave plate one inch in diameter (M at Fig. 3 from [3]) suspended with its plane horizontal from a quartz fiber about 25 centimeters long. A circularly polarized light beam (power $P$ = 80 mW, $\lambda$ = 1.2 $\mu$m,

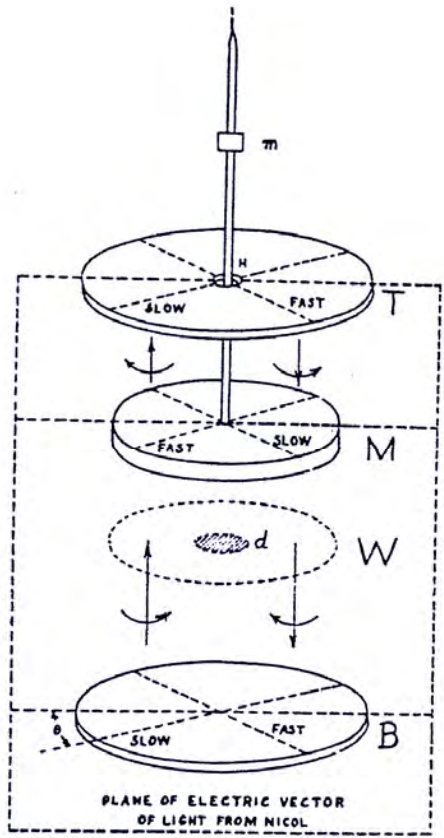

FIG. 3. Wave plate arrangement.

$\omega = 1.6 \cdot 10^{15} \, \mathrm{s}^{-1}$)
travels upwards, passing through the half-wave plate from below upwards. Because the plate reverses the handedness of the circular polarization of the beam, according to (6) and (12), the torque acting on the plate must be

$$\tau = 2P/\omega \qquad (14)$$

However, and this is the main point, in order to redouble the torque, the beam is reflected and passes through the plate a second time on the way back. For this, a fixed quartz quarter-wave plate T is mounted about 4 millimetres above the plate (Fig. 3). The *top* side of the upper plate was coated by evaporation with a reflecting layer of aluminium The rotation of the pendulum is observed by a telescope using the small mirror *m* at Fig.3. As a result, the torque exerting on the half-wave plate is 20 dyne cm. This result is in accordance with the formula

$$\tau = 4P/\omega. \qquad (15)$$

It is evident that the reflected beam cancels the energy flux in the Beth's apparatus, i.e. the Poynting vector $\mathbf{E} \times \mathbf{H} = 0$ in the experiment. Thus, according to equation (6), no angular momentum is contained in the double beam. So, no torque must act on the Beth plate according to eqn. (6). Why then the plate experiences the torque (15)?

We verify our claim $\mathbf{E} \times \mathbf{H} = 0$ here. Let us start from the Jackson beam (7) with $\omega = 1$ for simplicity,

$$\mathbf{E}_1 = \exp(iz - it)[\mathbf{x} + i\mathbf{y} + \mathbf{z}(i\partial_x - \partial_y)]E_0\,,$$

$$\mathbf{H}_1 = \exp(iz - it)[-i\mathbf{x} + \mathbf{y} + \mathbf{z}(\partial_x + i\partial_y)]E_0\,. \quad (16)$$

Changing the sign of $z$ we get the reflected beam. But the quarter-wave plate T changes the handedness of the circularly polarization of the beam. Thus, the sign of $y$ must be changed as well. So, the reflected beam is

$$\mathbf{E}_2 = \exp(-iz - it)[\mathbf{x} - i\mathbf{y} + \mathbf{z}(-i\partial_x - \partial_y)]E_0\,,$$

$$\mathbf{H}_2 = \exp(-iz - it)[-i\mathbf{x} - \mathbf{y} + \mathbf{z}(-\partial_x + i\partial_y)]E_0\,. \quad (17)$$

Adding together expressions (16) and (17) we get the total field

$$E_x = \Re[\exp(iz - it) + \exp(-iz - it)]E_0 = 2E_0 \cos z \cos t\,, \quad (18)$$

$$E_y = \Re[i\exp(iz - it) - i\exp(-iz - it)]E_0 = -2E_0 \sin z \cos t\,, \quad (19)$$

$$E_z = \Re[\exp(iz - it)(i\partial_x - \partial_y) + \exp(-iz - it)(-i\partial_x - \partial_y)]E_0$$

$$= -2(\sin z\partial_x + \cos z\partial_y)E_0 \cos t\,, \quad (20)$$

$$H_x = \Re[-i\exp(iz - it) - i\exp(-iz - it)]E_0 = -2E_0 \cos z \sin t\,, \quad (21)$$

$$H_y = \Re[\exp(iz - it) - \exp(-iz - it)]E_0 = 2E_0 \sin z \sin t\,, \quad (22)$$

$$H_z = \Re[\exp(iz - it)(\partial_x + i\partial_y) + \exp(-iz - it)(-\partial_x + i\partial_y)]E_0$$

$$= 2(\sin z\partial_x + \cos z\partial_y)E_0 \sin t \quad (23)$$

As a result we get

$$\mathbf{E} = 2[(\mathbf{x}\cos z - \mathbf{y}\sin z) - \mathbf{z}(\sin z\partial_x + \cos z\partial_y)]E_0 \cos t \quad (24)$$

$$\mathbf{H} = -2[(\mathbf{x}\cos z - \mathbf{y}\sin z) - \mathbf{z}(\sin z\partial_x + \cos z\partial_y)]E_0 \sin t \quad (25)$$

The $\mathbf{E}$ and $\mathbf{H}$ fields are parallel to each other everywhere. So, the Poynting vector is zero.

### 6.3  An explanation of Beth result

Formula (6) predicts the zero result of the Beth experiment because this formula is incorrect. As Ohanian wrote, the existence of angular momentum (6) is caused by circulating energy flow in the wave. In other words, equation (6) represents an orbital angular momentum of electromagnetic field. It is in accord with the fact that Maxwell's electrodynamics does not know spin. Spin is considered to be a pure quantum phenomenon. Maxwell electrodynamics knows the energy-momentum tensor $T^{\lambda\mu}$ (Maxwell-Minkowski tensor), but it does not know a spin tensor, or rather, *the spin tensor of the modern classical electrodynamics is zero*. We introduce classical spin into the electrodynamics. We introduce a spin tensor $Y^{\lambda\mu\nu}$ [10, 11], i.e. we add a spin term to equation (6):

$$J^{ij} = \int_V 2r^{[i}T^{j]0}dV + \int_V Y^{ij0}dV \, . \qquad (26)$$

The energy flux density, i.e. the Poynting vector $T^{0j} = T^{j0}$, is zero, $T^{j0} = 0$, in Beth's experiment. So, the first term on the right hand side of equation (26), i.e. the orbital term, is zero. However, spin flows from the beam into the Beth plate, and a torque acts on the plate due to a spin term.

The sense of the spin tensor $Y^{\lambda\mu\nu}$ is as follows. The component $Y^{ij0}$ is a volume density of spin. This means that

$$dS^{ij} = Y^{ij0}dV \qquad (27)$$

is the spin of electromagnetic field inside the spatial element $dV$ . The component $Y^{ijk}$ is a flux density of spin flowing in the direction of the $x^k$ axis. For example,

$$dS_z / dt = dS^{xy} / dt = d\tau^{xy} = Y^{xyz}da_z \qquad (28)$$

is the *z*-component of spin flux passing through the surface element $da_z$ per unit time, i.e. the torque acting on the element.

The explicit expression for the spin tensor is [10, 11] (see also Supplement)

$$Y^{\lambda\mu\nu} = A^{[\lambda}\partial^{|\nu|}A^{\mu]} + \Pi^{[\lambda}\partial^{|\nu|}\Pi^{\mu]}, \qquad (29)$$

where $A^{\lambda}$ and $\Pi^{\lambda}$ are magnetic and electric vector potentials which satisfy

$$2\partial_{[\mu}A_{\nu]} = F_{\mu\nu}, \quad 2\partial_{[\mu}\Pi_{\nu]} = -e_{\mu\nu\alpha\beta}F^{\alpha\beta} \qquad (30)$$

where $F^{\alpha\beta} = -F^{\beta\alpha}$, $F_{\mu\nu} = F^{\alpha\beta} g_{\mu\alpha} g_{\nu\beta}$ is the field strength tensor of a free electromagnetic field.

A relation between $\Pi$ and $F$ can be readily obtained in the vector form as follows.

If $\mathsf{div}\mathbf{E} = 0$, then $\mathbf{E} = \mathsf{curl}\,\Pi$, and if $\partial\mathbf{E}/\partial t = \mathsf{curl}\mathbf{H}$ as well, then

$$\partial\Pi/\partial t = \mathbf{H}. \tag{31}$$

This reasoning is analogous to the common:

If $\mathsf{div}\mathbf{H} = 0$, then $\mathbf{H} = \mathsf{curl}\mathbf{A}$, and if $\partial\mathbf{H}/\partial t = -\mathsf{curl}\mathbf{E}$ as well, then

$$\partial\mathbf{A}/\partial t = -\mathbf{E}. \tag{32}$$

Now use the spin tensor (29) for calculating of the spin flux into Beth plate. Since the orbital term is zero, wall terms, $\partial_x$, $\partial_y$, in equations (24), (25) may be neglected, and we have for the fields

$$\mathbf{E} = 2(\mathbf{x}\cos z - \mathbf{y}\sin z)E_0\cos t, \tag{33}$$

$$\mathbf{H} = -2(\mathbf{x}\cos z - \mathbf{y}\sin z)E_0\sin t, \tag{34}$$

$$\mathbf{A} = -\int \mathbf{E}dt = -2(\mathbf{x}\cos z - \mathbf{y}\sin z)E_0\sin t, \tag{35}$$

$$\Pi = \int \mathbf{H}\,dt = 2(\mathbf{x}\cos z - \mathbf{y}\sin z)E_0\cos t. \tag{36}$$

When calculating the spin tensor (29) the signature of metric tensor must be taken into account. Because $g^{ij} = -1$, $\partial^i = -\partial_i$. Thus, the spin flux density onto the low side of the Beth plate is

$$\mathbf{Y}^{xyz} = (A^x\partial^z A^y - A^y\partial^z A^x)/2 + (\Pi^x\partial^z\Pi^y - \Pi^y\partial^z\Pi^x)/2$$
$$= 2E_0^2(\sin^2 t + \cos^2 t) = 2E_0^2. \tag{37}$$

The same calculation for the domain above the plate gives $\mathbf{Y}^{xyz} = -2E_0^2$. This means that $S^{xy}$-component of the spin moves opposite the $z$-direction, i.e. towards the plate also. It follows that *the plate receives the spin flux density, or torque density, of* $4E_0^2$ *in the absence of energy flux!* Thus, the torque is

$$\tau = 4\int E_0^2\,da, \tag{38}$$

and recalling (10), we get ( $\omega = 1$ )

$$\tau = 4\mathsf{P} \tag{39}$$

as Beth's experiment shows.

It is remarkable that volume density of spin equals zero, i.e. $Y^{xy0} = 0$. The use of (35), (36) shows this. This is natural because the beams of the same handedness, which propagate in the opposite directions, are summed up. So, the Beth's double beam contains spin flux and energy without spin and energy flux.

Other applications of the spin tensor are presented in [12][13] and at online homepages: http://www.mai.ru/projects/mai_works/, and http://www.sciprint.org (see folder user Khrapko). Absorption and reflection of a circularly polarized beam is calculated there. Radiation from a rotating electric dipole and other topics are also considered in these works.

### 6.4 Supplement: Electrodynamics spin tensor

The standard classical electrodynamics starts from the free field canonical Lagrangian [14]

$$\mathsf{L}_c = -F_{\mu\nu}F^{\mu\nu}/4\,,\quad F_{\mu\nu} = 2\partial_{[\mu}A_{\nu]}\,,\tag{40}$$

Using this Lagrangian, by the Lagrange formalism physicists obtain the canonical energy-momentum tensor

$$T_c^{\ \lambda\mu} = \partial^\lambda A_\alpha \frac{\partial \mathsf{L}_c}{\partial(\partial_\mu A_\alpha)} - g^{\lambda\mu}\mathsf{L}_c = -\partial^\lambda A_\alpha F^{\mu\alpha} + g^{\lambda\mu}F_{\alpha\beta}F^{\alpha\beta}/4\,,\tag{41}$$

and the canonical total angular momentum tensor

$$J_c^{\ \lambda\mu\nu} = 2x^{[\lambda}T_c^{\ \mu]\nu} + Y_c^{\ \lambda\mu\nu}\tag{42}$$

where

$$Y_c^{\ \lambda\mu\nu} = -2A^{[\lambda}\delta_\alpha^{\mu]}\frac{\partial \mathsf{L}_c}{\partial(\partial_\nu A_\alpha)} = -2A^{[\lambda}F^{\mu]\nu}\,,\tag{43}$$

is the canonical spin tensor.

Unfortunately, the canonical tensors are not electrodynamics tensors. True electrodynamics tensors must be in accordance with experimental facts. In particular, it should be

$$\partial_\mu T^{\lambda\mu} = -F^{\lambda\mu}j_\mu = F^{\lambda\mu}\partial^\nu F_{\mu\nu}\,.\tag{44}$$

But $T_c^{\ \lambda\mu}$ has a wrong divergence,

$$\partial_\mu T_c^{\lambda\mu} = -\partial^\lambda A^\mu j_\mu = \partial^\lambda A^\mu \partial^\nu F_{\mu\nu}, \tag{45}$$

and is asymmetric. Physicists undertook an attempt to modify these tensors. They "put in by hand" specific addends [15, 16] to the canonical tensors and arrive to the standard energy-momentum tensor $\Theta^{\lambda\mu}$, the standard total angular momentum tensor $J_{st}^{\lambda\mu\nu}$, and the standard spin tensor $Y_{st}^{\lambda\mu\nu}$, which is zero,

$$\Theta^{\lambda\mu} = T_c^{\lambda\mu} - \partial_\nu \widetilde{Y}_c^{\lambda\mu\nu}/2$$

$$= -\partial^\lambda A_\nu F^{\mu\nu} + g^{\lambda\mu} F_{\alpha\beta} F^{\alpha\beta}/4 + \partial_\nu (A^\lambda F^{\mu\nu}),$$

$$\widetilde{Y}_c^{\lambda\mu\nu} \stackrel{def}{=} Y_c^{\lambda\mu\nu} - Y_c^{\mu\nu\lambda} + Y_c^{\nu\lambda\mu} = -2A^\lambda F^{\mu\nu}, \tag{46}$$

$$J_{st}^{\lambda\mu\nu} = J_c^{\lambda\mu\nu} - \partial_\kappa (x^{[\lambda} \widetilde{Y}_c^{\mu]\nu\kappa}), \tag{47}$$

$$Y_{st}^{\lambda\mu\nu} = J_{st}^{\lambda\mu\nu} - 2x^{[\lambda} \Theta^{\mu]\nu} = Y_c^{\lambda\mu\nu} - \widetilde{Y}_c^{[\lambda\mu]\nu} = 0. \tag{48}$$

But the standard tensors are not true electrodynamics tensors either:

1. $\Theta^{\lambda\mu}$ obviously contradicts experiments. It is asymmetric and has wrong divergence as well

$$\partial_\mu \Theta^{\lambda\mu} = \partial_\mu T_c^{\lambda\mu} = \partial^\lambda A^\mu \partial^\nu F_{\mu\nu} \tag{49}$$

Tensor $\Theta$ is never used. The Maxwell tensor

$$T^{\lambda\mu} = -F^\lambda{}_\alpha F^{\mu\alpha} + g^{\lambda\mu} F_{\alpha\beta} F^{\alpha\beta}/4 \tag{50}$$

is used in the electrodynamics instead of $\Theta^{\lambda\mu}$.

2. The main defect is the absence of spin, $Y_{st}^{\lambda\mu\nu} = 0$. In contrast to the canonical pair, $T_c^{\lambda\mu}, Y_c^{\lambda\mu\nu}$, the standard pair, $\Theta^{\lambda\mu}, Y_{st}^{\lambda\mu\nu} = 0$, is defective. *The standard energy-momentum tensor is not accompanied by a spin tensor.*

Thus the Belinfante-Rosenfeld procedure [15, 16] is not fit for determining the true electrodynamics tensors. This procedure is

$$\Theta^{\lambda\mu} = T_c^{\lambda\mu} + t_{st}^{\lambda\mu}, \qquad t_{st}^{\lambda\mu} = -\partial_\nu \widetilde{Y}_c^{\lambda\mu\nu}/2 = \partial_\nu (A^\lambda F^{\mu\nu}), \tag{51}$$

$$\underset{st}{Y}{}^{\lambda\mu\nu} = \underset{c}{Y}{}^{\lambda\mu\nu} + \underset{st}{s}{}^{\lambda\mu\nu} = 0, \quad \underset{st}{s}{}^{\lambda\mu\nu} = -\widetilde{\underset{st}{Y}}{}^{[\lambda\mu]\nu} = 2A^{[\lambda}F^{\mu]\nu}. \tag{52}$$

Another way of using the canonical pair $\underset{c}{T}{}^{\lambda\mu}, \underset{c}{Y}{}^{\lambda\mu\nu}$ is presented in [11 – 13]. Note that the Maxwell tensor (50) can be obtained by adding a term

$$t^{\lambda\mu} = T^{\lambda\mu} - \underset{c}{T}{}^{\lambda\mu} = \partial_\nu A^\lambda F^{\mu\nu} \tag{53}$$

to the canonical energy-momentum tensor $\underset{c}{T}{}^{\lambda\mu}$. Here a question arises, what term $s^{\lambda\mu\nu}$, instead of $\underset{st}{s}{}^{\lambda\mu\nu}$, must be added to the canonical spin tensor $\underset{c}{Y}{}^{\lambda\mu\nu} = -2A^{[\lambda}F^{\mu]\nu}$ for changing it from the canonical spin tensor to an unknown electrodynamics spin tensor $Y^{\lambda\mu\nu} = \underset{c}{Y}{}^{\lambda\mu\nu} + s^{\lambda\mu\nu}$? Our answer is [11 – 13] that the addends $t^{\lambda\mu}$, $s^{\lambda\mu\nu}$ must satisfy the relationship

$$\partial_\nu s^{\lambda\mu\nu} - 2t^{[\lambda\mu]} = 0, \quad \text{i.e.} \ \partial_\nu s^{\lambda\mu\nu} - 2\partial_\alpha A^{[\lambda}F^{\mu]\alpha} = 0. \tag{54}$$

The simple expression

$$s^{\lambda\mu\nu} = 2A^{[\lambda}\partial^{\mu]}A^\nu \tag{55}$$

satisfies Eq. (54). So, the suggested electrodynamics spin tensor is

$$2\underset{e}{Y}{}^{\lambda\mu\nu} = \underset{c}{Y}{}^{\lambda\mu\nu} + s^{\lambda\mu\nu} = -2A^{[\lambda}F^{\mu]\nu} + 2A^{[\lambda}\partial^{\mu]}A^\nu.$$

$$= 2A^{[\lambda}\partial^{|\nu|}A^{\mu]} \tag{56}$$

The expression (56) was obtained heuristically. It is not a final one. Spin tensor (56) is obviously not symmetric in the sense of electric - magnetic symmetry. It represents only the electric field, $\mathbf{E}, \ \mathbf{A} = -\int\mathbf{E}dt$. A true spin tensor of electromagnetic waves must depend symmetrically on the magnetic vector potential $A_\alpha$ and on an electric vector potential $\Pi_\alpha$ (30). So the spin tensor of electromagnetic waves has the form (29).

# 6    Unsolved Problems in other areas of Science

> Any sufficiently advanced technology is indistinguishable from magic.
> --Arthur C. Clarke

It seems worthwhile to consider here a number of other questions which may be of value to ponder.

Not all of these questions are very serious.

## 6.1  Relativity theory

Here are some questions pertaining to relativity theory.

Questions:

i.    What is the geometric shape of a relativistic rotating plane (Unruh experiment)? Is it Minkowskian, Euclidean, or Riemannian?

ii.   A disc rotating at high speed will exert out-of-plane forces resembling an accelerating field. Is the principle of equivalence also applicable for this process?

iii.  Will someone inside an elevator in free-fall and rotating around its vertical centre, feel a gravitational force? Or will he feel a gravitational force larger than what equivalence principle requires? Does the equivalence principle remain applicable here?

iv.   An aeroplane flies at an altitude of 1 km. The co-pilot drops an elevator-room <u>without</u> a passenger inside it. After one second has elapsed, the co-pilot drops four grenades in the direction of the freely-falling elevator's path. The question: Will the grenades reach the elevator before it reaches the ground? If no, why? If yes, which grenade?

v.    What is the effect of space temperature on the spacetime curvature? Or does temperature affect the metric of spacetime? (Ref. D. Colladay, et al., arXiv:hep-ph/0602071)  This question reso-nates Einstein words (something like): 'If you stay near flaming fire, a minute will be felt like a century.'

## 6.2 Questions related to Bell's theorem

Here are some questions pertaining to Bell's theorem in the measurement theory of Quantum Mechanics.

Question: Is it possible to modify Bell's theorem to include multivalued photon-pair position? (see Ref: R. Ogden, arXiv: cs.IT/0507032).

### 6.2.1 Introduction

It is generally accepted that Bell's theorem [81] is quite exact for describing the linear hidden-variable interpretation of quantum reality. Therefore null result of this proposition implies that no hidden-variable theory could provide a sound explanation of quantum reality.

Nonetheless, after further thought we can see that Bell's theorem is nothing more than another kind of abstraction of quantum observation based on a set of assumptions and propositions [87]. Therefore, one should be careful before making further generalizations on the null result from experiments which are 'supposed' to verify Bell's theorem. For example, the most blatant assumption of Bell's theorem is that it takes into consideration only the classical statistical problem of chance of outcome A or outcome B. In other words, it is nothing more than an adaptation of the 'coin toss' problem into the complicated quantum reality, whereas simultaneous appearance of a photon pair could occur at different places, and hence there is a small chance of the effect of multivalued logic.

Therefore in the present paper we will extend this Bell's theorem into a modified version which takes into consideration this multivalued outcome, in particular using the information fusion theory of Dezert-Smarandache [82][83][84]. We suppose that in quantum reality the outcome of $P(A \cup B)$ and also $P(A \cap B)$ shall also be taken into consideration. This is where DSmT theory could be found useful. [82]

It could be expected that such a modified version would be useful for describing quantum reality in a more precise way. Further experiments are of course recommended in order to verify or refute this proposition.

### 6.2.2    Bell's theorem and its inherent assumption

Despite widespread belief of its ability to describe hidden-variables of quantum reality [81], it will be noted that Bell's theorem starts with a set of assumptions inherent in its formulation. It is assumed that each of pair of particles possesses a particular value of $\lambda$, and we define the quantity $p(\lambda)$ so that the probability of a pair being produced with a value of $\lambda$ between $\lambda$ and $\lambda + d\lambda$ is $p(\lambda)d\lambda$. It is also assumed that this is normalized so that:

$$\int p(\lambda)d\lambda = 1 \tag{1}$$

Further analysis shows that the integral that measures the correlation between two spin components that are at an angle of $(\delta - \phi)$ to each other, is therefore equal to $C''(\delta - \phi)$, where C'' represents *average* chance after some number of observations. We can therefore write:

$$\left| C''(\phi) - C''(\delta) \right| - C''(\delta - \phi) \leq 1 \tag{2}$$

which is known as Bell's theorem, and it represents any local hidden-variable theorem. But it will be noted that his theorem cannot actually be tested completely because it assumes that all particle pairs have been detected. In other words, we find that a hidden assumption behind the Bell's theorem is that it uses a classical probability assertion, and therefore it neglects the possibility of including $P(A \cup B)$ and $P(A \cap B)$, which may be useful in describing quantum reality.

### 6.2.3    Set theoretic extension of Bell's theorem, DSmT

In the context of a physical interpretation of information [88, p.378], Barrett has noted that "there ought to be a set theoretic language which applies directly to all quantum interactions." This is because the idea of a bit is itself straight out of *classical set theory*, the definitive and unambiguous assignment of an element of the set {0,1}, and so the assignment of an information content of the photon itself is fraught with the same difficulties [88, p.378].

We know that for quantum reality, the photon could appear (sometimes) in two different places at the same time [85]. Therefore it could be useful to introduce a modified version of an information assertion which is capable of representing such multivalued positions. In other words we should find an extension to the *standard* proposition in statistical theory [88, p.388]:

$$P(AB|C) = P(A|BC)P(B|C) = P(B|AC)P(A|C) \qquad (3)$$

$$P(A|B) + P(\overline{A}|B) = 1 \qquad (4)$$

Such an extension is already known in the area of information fusion [82], known as Dempster-Shafer theory:

$$m(A) + m(B) + m(A \cup B) = 1 \qquad (5)$$

Furthermore, Dezert and Smarandache [82] introduced a further refinement of Dempster-Shafer theory by taking into consideration chance of observe the intersection between A and B:

$$m(A \cup B) + m(A) + m(B) + m(A \cap B) = 1 \qquad (6)$$

Therefore, introducing this extension from equation (6) into equation (2), one finds a modified version of Bell's theorem in the form:

$$\left| C''(\phi) - C''(\delta) \right| - C''(\delta - \phi) + C''(\delta \cup \phi) + C''(\delta \cap \phi) \leq 1 \qquad (7)$$

which could be called the modified Bell's theorem according to Dezert-Smarandache theory [82]. Of course, further experiment is recommended in order to verify and to find various implications of this new proposition.

Interestingly, one could consider further refinement of equation (6) in the context of Unification of Fusion Theories (UFT) as recently proposed by Smarandache, by considering chance to observe negation of A and B, and also by considering chance to observe '*none of these*' (i.e. intersection of negation of A and negation of B), which can be written as follows:

$$m(A \cup B) + m(A) + m(B) + m(A \cap B) +$$
$$m(\neg A \cup \neg B) + m(\neg A) + m(\neg B) + m(\neg A \cap \neg B) = 1 \qquad (6a)$$

Where $\neg X$ represents negation of X. Deriving an extension of Bell's theorem using this new proposition (UFT) remains open, but we leave this problem for the readers.

### 6.2.4    Alternative interpretation using Kholevo's theorem

Alternatively, we can offer another interpretation of the above modification of Bell's theorem from the viewpoint of Kholevo's theorem as described by Schumacher [88, p.31].

Let's suppose the information content of a message is given by the information function [88, p.31]:

$$H(X) = -\sum_i p(x_i) \log p(x_i) \qquad (8)$$

where H(X) is expressed in 'bits'. If the channel is 'noisy', then the receiver may have a non-zero degree of uncertainty of the message X – on average—as follows:

$$H(X|Y) = \sum_k p(y_k)H(X|y_k) = H(X,Y) - H(Y) \qquad (9)$$

Then the receiver has gained an amount of information:

$$H(X:Y) = H(X) - H(X|Y) = H(X) + H(Y) - H(X,Y) \qquad (10)$$

which is equivalent to Bell's theorem given that :

$$H(X:Y) \leq 1 \qquad (11)$$

Therefore we can write [88, p.31]:

$$H(X) + H(Y) - H(X,Y) \leq 1 \qquad (12)$$

Now let's suppose that the ensemble average of any quantum observable is described by the density operator [88, p.31]:

$$\rho = \sum_i p(x_i).\rho(x_i) \qquad (13)$$

For example, the average signal energy is $< E >= Tr\rho H$ . The entropy of the signal ensemble is defined by:

$$S[\rho] = Tr\rho.\log\rho \qquad (14)$$

which is equivalent to the information function H(X).

Kholevo sets a bound on H(X), the amount of information transmitted by the quantum channel Q [88, p.31]:

$$H(X:A) \leq S[\rho] - \sum_i \rho(x_i)S[\rho(x_i)] \qquad (15)$$

Since the 'subtracted' term from the right side is non-negative, it trivially follows that

$$H(X:A) \leq S[\rho] \qquad (16)$$

that is, the quantum channel Q can deliver an amount of information no greater than the entropy of the ensemble.

A further plausible extension of the above proposition is introduced by Maesser & Uffinck [88, p.31] which describes for any mass X coded into Q:

$$H(X:A) + H(X:B) = H(A) + H(B) - \\ \left[H(A|X) + H(B|X)\right] \leq 2\log N - C \qquad (17)$$

For a spin-1/2 system Q, C=log2, and so we get [8, p.31]:

$$H(XY : \sigma_x) + H(XY : \sigma_y) \le 1bit \qquad (18)$$

which can be considered as an alternative interpretation of Bell's theorem in the context of information theory.

Now by using Dezert-Smarandache theory [82] to extend the standard statistical assumption as described in the preceding section, we get:

$$H(A) + H(B) - \left[ H(A|X) + H(B|X) \right] +$$
$$H(A \cup B) + H(A \cap B) \le 2 \log N - C \qquad (18a)$$

And for a spin-1/2 system Q with C=log 2, we get a modified version of equation (18) as follows:

$$H(\sigma_z) + H(\sigma_x) - \left[ H(\sigma|XY) + H(\sigma_x|XY) \right] +$$
$$H(\sigma_z \cup \sigma_x) + H(\sigma_z \cap \sigma_x) \le 1bit \qquad (18b)$$

which can be viewed as the Kholevo-Dezert-Smarandache interpretation of Bell's theorem.

In the above section we discussed a plausible modified version of Bell's theorem which could take into consideration the chance to observe outcomes beyond classical statistical theory, in particular using the information fusion theory of Dezert-Smarandache. It is recommended to conduct further experiments in order to verify and also to explore various implications of this new proposition, including perhaps for quantum computation theory [85].

## 6.3  Questions related to Mind-Matter interaction, hidden mystery of water

Here are some questions pertaining to the mind-matter interaction, in particular with links to the mystery of water structure.

Questions:
- Is there a valid theoretical relation to describe the Mind-Matter interaction in the context of quantum physics?
- Is there experimental observation supporting the above linkage between Mind and Matter?

- Could this observation be used to explain other known phenomena in quantum physics?
- What is the plausible linkage between the Mind–Matter interaction and Cymatics observation?
- Are there new observations predicted by the Mind-Matter interaction paradigm, in particular concerning the hidden structure of water?

Background information:

- A Japanese quantum physicist [i], performed a series of experiments on water crystals and revealed the fact that water is receptive to external messages. The formation of water crystals is correlated to exposure of the water to messages from human language, music, and printed characters. In the meantime, other researchers have also conjectured that the special structure of water could appear as a degree of coherence [ii].

- Similar experiments were performed by G. Thomas etc. (http://www.dallasinstitute.org/Programs/Previous/Fall%202001/talks/gthomaswater.htm)

- Actually these works are not really new. In 1967, Hans Jenny, a Swiss doctor, artist, and researcher, published the bilingual book Kymatik - Wellen und Schwingungen mit ihrer Struktur und Dynamik/ Cymatics - The Structure and Dynamics of Waves and Vibrations. The tonoscope was constructed to make the human voice visible without any electronic apparatus as an intermediate link. This yielded the amazing possibility of being able to see the physical image of the vowel, tone or song a human being produced directly.

 (*http://www.mysticalsun.com/cymatics/cymatics.html* )

- Is it real or merely an untested hypothesis?

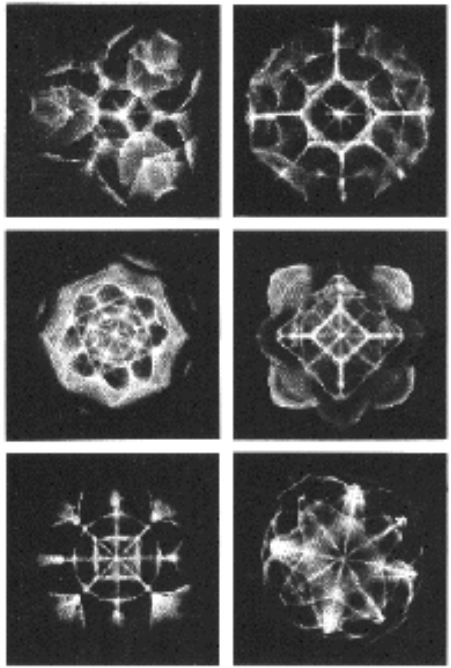

Figure 5. Swinging water drop (by Hans Jenny) Ref. *http://www.mysticalsun.com/cymatics/cymatics.html*

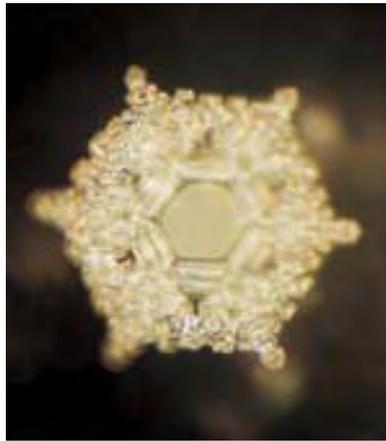

Figure 6. Crystal Photo of Water That Had Been Shown Label "Love and Thanks" (http://www.dallasinstitute.org/Programs/Previous/Fall%202001/talks/gthomaswater.htm)

Ref.    (i)    http://www.pureinsight.org/pi/index.php?news=1626,    (ii) www.0disease.com/0waterheal.html

## 6.4  Quaternionic wave interpretation of superconductors

Here are some questions pertaining to quaternion number and superconductivity.

Questions:

- Is there a valid theoretical link between quaternion number, Maxwell's equations and superconductivity phenomena?
- Is there experimental observation supporting the above linkage between quaternion number and electrodynamics of superconductors?
- Could this observation be used to explain other known phenomena such as the Meissner effect or the Josephson junction?
- What is the plausible linkage between this quaternionic electrodynamics view of superconductivity and BCS theory?
- Are there new observations predicted by a quaternionic electrodynamics viewpoint that has not been predicted before by BCS theory?

Background theory:

It is often recognized that there are some aspects of Maxwell's electromagnetic equations which are not replaced by quantum mechanics. In the light of some recent development in this apparently 'blurred area' between classical electrodynamics and quantum mechanics [89], we would like to argue that the Ginzburg-Landau-Schrödinger (Gross-Pitaevskii) wave equation, which is often used in the context of superconductivity, could also be extended into the form of nonlinear electromagnetic equations [93], in particular via the use of biquaternion numbers.

In this section, we begin with a plausible relation between the Ginzburg-Landau-Schrödinger (Gross-Pitaevskii) wave equation and Schrödinger equation, and from Schrödinger equation to biquaternion number. In effect, we argue that the Ginzburg-Landau-Schrödinger (Gross-Pitaevski) wave equation could be represented in the form of double biquaternion wave equations. Because it is known that biquaternion wave equations could represent both Maxwell's equations and Dirac's equation [94], then it seems that the double biquaternion wave equations could be interpreted as nonlinear versions of the Maxwell-Dirac equations, i.e. to include new phenomena which cannot be predicted by ordinary electrodynamics theory. From this viewpoint, we suggest verification of such phenomena, either numerically or via experiments.

### From Ginzburg-Landau to (double-wave) Schrödinger equation

Consider the well-known Gross-Pitaevskii equation in the context of superfluidity or superconductivity [90]:

$$i\hbar \frac{\partial \Psi}{\partial t} = -\frac{\hbar^2}{2m} \Delta \Psi + (V(x) - \gamma |\Psi|^{p-1})\Psi, \tag{19}$$

where p<2N/(N-2) if N≥3. In physical problems, the equation for p=3 is known as the Gross-Pitaevskii equation. Equation (19) has a standing wave solution quite similar to Schrödinger's equation, in the form:

$$\Psi(x,t) = e^{-iEt/\hbar}.u(x) \tag{20}$$

Substituting equation (20) into equation (19) yields:

$$-\frac{\hbar^2}{2m}\Delta u + (V(x) - E)u = |u|^{p-1}u, \tag{21}$$

which is nothing more than the time-independent linear form of Schrödinger equation, except for the term $|u|^{p-1}$ [90].

Using Maclaurin series expansion, we get for (20):

$$\Psi(x,t) = \left(1 - iEt/\hbar + \frac{(Et/\hbar)^2}{2!} + \frac{(-iEt/\hbar)^3}{3!} + \dots\right)u(x) \qquad (22)$$

Therefore we can say that the standing wave solution of the Gross-Pitaevskii equation (20) is similar to the standing wave solution of Schrödinger equation (u), except for a nonlinear term which comes from the Maclaurin series expansion (22). By neglecting third and fourth terms of equation (22), one gets:

$$\Psi(x,t) = \left(1 - iEt/\hbar.\right)u(x) \qquad (23)$$

Note that this equation (23) is very near to the hyperbolic form $z = x + iy$. Therefore one could conclude that the standing wave solution of the Gross-Pitaevskii equation is merely an extension of the ordinary solution of Schrödinger equation into a Cauchy (imaginary) plane. We will use this result in the following section, but first we consider how to derive a biquaternion from Schrödinger equation.

**From Schrödinger's equation to biquaternion electrodynamics**

It is known that solutions of Riccati's equation are logarithmic derivatives of solutions of Schrödinger's equation, and *vice versa* [91]:

$$u'' + vu = 0 \qquad (24)$$

The biquaternion of differentiable function of x=(x$_1$,x$_2$,x$_3$) is defined as [91]:

$$Dq = -div(q) + grad(q_o) + rot(q) \qquad (25)$$

By using alternative representation of Schrodinger equation [91]:

$$\left(-\Delta + u\right)f = 0, \qquad (26)$$

where f is twice differentiable, and introducing the quaternion equation:

$$Dq + q^2 = -u. \qquad (27)$$

Then we could find q, where q is a purely vectorial differentiable biquaternion valued function [91].

We note here that solutions of (26) are related to solutions of (25) as follows:
➔ For any nonvanishing solution f of (26), its logarithmic derivative:

$$q = \frac{Df}{f}, \qquad (28)$$

is a solution of equation (27), and *vice versa*. [91]

Furthermore, we also note that for an arbitrary twice-differentiable scalar function f, the following equality is permitted [91]:

$$\left(-\Delta + u\right)f = \left(D + M^h\right)\left(D - M^h\right)f \ , \tag{29}$$

provided h is solution of equation (27).

Therefore in summary, *given a particular solution of Schrödinger equation (26), the general solution reduces to the first order equation* [91, p.9]:

$$\left(D + M^h\right)F = 0, \tag{30}$$

where

$$h = \frac{D\sqrt{\varepsilon}}{\varepsilon} \ . \tag{31}$$

Interestingly, equation (30) is equivalent to Maxwell's equations.[91]

Now we can generalize our result from the preceding section, in the form of the following conjecture:

<u>Conjecture 1</u>: Given *a particular solution of Schrödinger equation (26), then the approximate solution of the Gross-Pitaevskii equation (19) reduces to the first order equation:*

$$\left(1 - iEt / \hbar.\right)\left(D + M^h\right)F = 0 \tag{32}$$

## 6.5 Solar dynamics

With regards to Solar dynamics, there are apparently some unsolved problems (read for instance http://thermalphysics.org/Sun.evidence.1.pdf). Glen Deen <glen.deen@gte.net> contributed the following problems:

Questions:

-   Could some of these neutrons be captured by atomic nuclei, raising their Z numbers, and the rest experience beta decay and become ionized hydrogen?
-   Could these nuclear reactions and beta decay events be sufficient for all the Sun's thermal energy and hence luminosity? (There might be no need for the fusion of protons into helium nuclei.)

- Couldn't a proton capture a neutron and become a deuterium nucleus?
- Couldn't a deuterium nucleus capture a neutron and become a tritium nucleus?
- Couldn't a tritium nucleus experience beta decay and become a helium-3 nucleus?
- Couldn't a helium-3 nucleus capture a neutron and become a helium-4 nucleus? And so on up the periodic table to iron and beyond?
- What if the real problem is that the Sun's thermal energy does not come from 4p -> 4He + 2e+ + 2nu_e as the standard model claims?
- What if all the thermal power in the Sun comes from the neutron capture reactions and beta decays due to the neutron flux from the core?
- Concerning the abundance of all the isotopes in the photosphere: Why can't we compute the thermal energy released by each neutron capture up the periodic table and tune the neutron flux to make the total power produced by neutron capture and decay equal to the known power of the Sun and add up the number of neutrinos produced by each reaction to see if it matches observations?

Background argument:

Eliminating thermonuclear reactions of protons fusing into helium nuclei in the core as the Sun's source of thermal power means that the core's temperature can be much lower than the million Kelvins required by the gaseous model. The neutron flux energy source idea is more compatible with the condensed matter model. Perhaps the neutrons at great depths are moving too fast to be captured so there is very little heat generated deep inside the Sun. That would mean the temperature might rise as you move up radially from the core to the surface because the neutrons slow down on their way up and their capture becomes more likely as they approach the surface. What is interesting about this idea is that Population I stars like the Sun would produce their own heavy elements by neutron capture.

Furthermore, there is perhaps link to the Neutrino problem, because researchers found some high-energy muon and tauon neutrinos, but not the number they expected to account for the 2/3 shortfall of electron neutrinos only. Even so, he says that the standard model of a massless neutrino needs to be revised to permit the oscillation between the three species of neutrinos.

Other Ref.: (i)
http://www.space.com/scienceastronomy/060619_mystery_monday.html; (ii) Robitaille, "The Solar Photosphere: Evidence for Condensed Matter," *Progress in Physics* Vol. 2 No. 2 (2006), www.ptep-online.com/index-files/2006/PP-05-04.pdf.

## 6.6 The Science of Conservation of Energy and Modified Newton-Coulomb theory

This is a Contribution by Fu Yuhua (fuyh@cnooc.com.cn).

Background argument: Introduction to the Science of Conservation of Energy

**The Science of Conservation of Energy** is presented by taking the law of Energy Conservation as the foundation and central factor. For all the problems concerned with energy, the law of Energy Conservation is the *only* truth; other laws will be derived by the law of Energy Conservation, or verified by it, or proved wrong. In this section, some questions are discussed.

Firstly, the relationship between force, mass, and velocity is reconsidered according to the law of Energy Conservation, and the general expression given by F=f(m,v,x,y,z,t), as for the standard form, should be derived by the law of Energy Conservation.

Secondly, other laws such as the law of gravity and Coulomb Law, are re-derived by the law of Energy Conservation. In passing, the rule for changing the gravitational coefficient (the so-called gravitational constant, G) is given.

Thirdly, other laws will be verified and yet others proved to be erroneous, such as the law of angular momentum conservation and the law of momentum conservation (the results given by them contradict the law of Energy Conservation).

Fourthly, changing an old subject into a new subject; for example, changing Newton's mechanics into New Newtonian mechanics in which the law of Energy Conservation is taken as the source law. From the source law, the law of gravity and Newton's second law can be derived. New Newtonian mechanics can be used for partly replacing relativity and solving problems that cannot be solved by (special) relativity theory. The science of conservation of energy may be widely used in physics, mechanics, engineering, chemistry, medicine, biology and the like. We also discuss some unsolved problems, such as the dimension of space and the dimension of time.

Some related Questions:

- Is it possible to come up with a form of Coulomb law which has repulsion effect near R = 0?

- Can we re-derive Newton's second law from first principles –like the principle of Energy Conservation? And does it imply a modified version of Newton's second law?

- Is it possible to come up with a flat-metric gravitation theory based on a modification of Newton's second law, which can explain the standard tests of General Relativity (see also Synge's theory)?

- Note that the Newton's second law is written as : d[mv]/dt = v.[dm/dt] + m.[dv/dt]. Therefore when [dm/dt]=0, we get F=m.[dv/dt] = m.a. Is there a (Haussdorf) fractal version of this formula and if so what is it?

- Is there a fundamental motivation to come up with a modification of Newton's second law?

- Is there observational data supporting a fractal version of Newton's second law?

<u>A new proposition:</u>

To discuss the possibility of deriving the Coulomb law and Newton's second law theoretically, according to the law of Energy Conservation, the variable dimension fractal method is developed, and used to improved Newton's second law and the Coulomb law in an example (a small charge ball moves down along a long incline within the electric field due to an electrified globe). The results from this example with constant dimension fractal form are as follows: the improved Coulomb law (inverse non-square Coulomb law), $f = kq_1q_2 / r^{1.99989}$; the improved Newton's second law $F = ma^{1.01458}$.

The standard form of Coulomb law ascertained from experimental results, it reads

$$f = \frac{kq_1q_2}{r^2} \qquad (33)$$

The Coulomb law and the law of gravity presented by Newton are inverse square laws, while the inverse non-square gravitational law was presented in reference [104] according to the improved Newton's formula of universal gravitation presented in reference [105]. The main results are as follows: the inverse non-square gravitational law with the form of a variable dimension fractal reads, $F=-GMm/r^D$, where $D=f(r)$ instead of $D=2$. The values of D are different for different problems.

Similar to the inverse non-square law of gravitation, is two charged bodies in relative motion (in this case the Newton's second law must be considered); the force between the two bodies will agree with the inverse non-square Coulomb law. Newton's second law in accordance with experimental results, it reads

$$F = ma \qquad (34)$$

Can these two laws be derived theoretically? It is possible in the case that there is a more extensive law. The law of Energy Conservation can be used for this important task. The reason for this is that the Coulomb law and Newton's second law can be used for handling macrocosmic physical phenomena only, while the law of Energy Conservation can be used for handling the macrocosmic and the microcosmic physical phenomena.

To discuss the possibility of deriving the Coulomb's law and Newton's second law theoretically, according to the law of Energy Conservation, the variable dimension fractal method is developed, and used to improved Newton's second law and the Coulomb law in an example (a small charged ball moves down along a long incline within the electric field due to an electrified globe). Since the analytic process is complicated, the results suitable for this example with a constant dimension fractal form will be given.

### 7.6.1. Variational principle for deriving the Coulomb law and Newton's second law simultaenously

The law of Energy Conservation is a basic one in natural science. Its main content can be stated briefly as follows, in a closed system, the total sys-temic energy is equal to a constant. Now the variational principle established by the law of Energy Conservation can be given by the *least-square method*. Supposing that the initial total energy of a closed system equals $W(0)$ andfor time $t$ the total energy equals $W(t)$, then according to the law of Energy Conservation:

$$W(0) = W(t) \tag{35}$$

This can be written as:

$$R_{\text{w}} = \frac{W(t)}{W(0)} - 1 = 0 \tag{36}$$

According to the *least-square method*, for the interval $[t_1, t_2]$ we can write the following variational principle:

$$\Pi = \int_{t_1}^{t_2} R_{\text{w}}^2 \mathrm{d}t = \min{}_0 \tag{37}$$

where $\min_0$ denotes the minimum value of the functional $\Pi$ and it should be equal to zero [106]. Besides the time coordinate, another one can also be used. For example, for the interval $[x_1, x_2]$, the following variational principle can be given according to the law of Energy Conservation

$$\Pi = \int_{x_1}^{x_2} R_{\text{w}}^2 \mathrm{d}x = \min{}_0 \tag{38}$$

The above-mentioned principle is established by using the law of Energy Conservation directly. Sometimes, a certain principle should be established by using the law of Energy Conservation indirectly. For example, a special physical quantity $Q$ may be of interest; not only can it be calculated by using the law of Energy Conservation, but it can also be calculated by using other laws (for this discussion they are the Coulomb law and Newton's second law). For distinguishing the values, let's denote the value given by other laws as $Q$, while $Q'$ denotes the value given by the law of Energy Conservation, then the value of $R_{\text{w}}$ can be redefined as follows

$$R_{\text{w}} = \frac{Q}{Q'} - 1 = 0 \tag{39}$$

Substituting Eq.(39) into Eqs.(37) and (36) as $Q'$ is the result calculated with the law of Energy Conservation, yields the variational principle established by using the law of Energy Conservation indirectly. Otherwise, it is clear that the extent of the value of $Q$ accords with $Q$'.

### 7.6.2. <u>Improved Newton's Second Law and Coulomb Law with the Form of a Variable Dimension Fractal and the Like</u>

In Newton's mechanics, the law of gravity has been improved. For example, in reference [105], the following improved formula was obtained

$$F = -\frac{GMm}{r^2} - \frac{3G^2M^2mp}{c^2r^4} \qquad (40)$$

where G — is gravitational constant, M and m — the masses of the two bodies, r — the distance between the two bodies, c — the velocity of light, p — the half normal chord for the body m moving around the body M along with a curve, and the value of p is given by:

$\quad$ p = a(1-e$^2$)$\quad$ (for ellipse)
$\quad$ p = a(e$^2$-1)$\quad$ (for hyperbola)
$\quad$ p = y$^2$/2x$\quad$ (for parabola)

For the problem of gravitational deflection of a photon orbit around the Sun and the problem of the advance of Mercury's perihelion, by using the improved formula of universal gravitation, the same results as given by general relativity can be obtained.

Referring to Eq.(40) the general form of the improved law of gravity can be written as follows

$$F = -\frac{GMm}{r^2}(1 + \frac{a_1}{r^2} + \frac{a_2}{r^4} + \cdots) \qquad (41)$$

Similarly, besides the static case, the general form of the improved Coulomb law can be written as follows:

$$f = \frac{kq_1q_2}{r^2}(1 + \frac{a_1}{r^2} + \frac{a_2}{r^4} + \cdots) \qquad (42)$$

In addition, the fractal method has yielded excellent results in many fields recently. The fractal distribution reads [107]:

$\quad$ N = $\dfrac{C}{r^D}$

Where r— represents characteristic scale such as length, time and the like; N— represents a quantity related to r, such as temperature, force and the like; C— represents a constant to be determined; and D—is (Haussdorf) fractal dimension.

For the case of D is a constant, this kind of fractal can be called a *constant dimension fractal*. For the case of D is not a constant, this kind of fractal can be called a *variable dimension fractal* [108-110].

The general form of the improved law of Coulomb with the form of variable dimension fractal can be written as follows

$$f = \frac{kq_1 q_2}{r^D} \tag{43}$$

where, $D = f(r)$. For example, it may take the form

$$D = a_1 + a_2 r + a_3 r^2 + \cdots \tag{44}$$

In this discussion, only the form of the *constant dimension fractal* will be taken, i.e., $D = \text{const}$.

In Newton's mechanics, the Newton's second law can also be improved, the general form may be written as

$$F = ma + k_1 a^2 + k_2 a^3 + \cdots \tag{45}$$

Similarly, the general form of the Newton's second law with the form of a *variable dimension fractal* may be written as

$$F = ma^{D'} \tag{46}$$

where, $D' = f(r)$, for example, it may be written as

$$D' = k_1 + k_2 a + k_3 a^2 + \cdots \tag{47}$$

In this discussion, only the form of the *constant dimension fractal* will be taken, i.e., $D' = \text{const}$. For the sake of convenience, it may be written as

$$F = ma^{1+\varepsilon} \tag{48}$$

where $\varepsilon = \text{const}$.

### 7.6.3. <u>A method for deriving Coulomb Law and Newton's second law</u>

Substituting Eq.(42) or Eq.(43) and the related quantities calculated by Eq.(45) or Eq.(46) into Eq.(37) or (38), then the equations derived by the condition for an extremum can be written as follows

$$\frac{\partial \Pi}{\partial a_i} = \frac{\partial \Pi}{\partial k_i} = 0 \tag{49}$$

After solving these equations, the improved form of Coulomb law and Newton's second law can be reached at once. According to the value of $\Pi$, the effect of the solution can be interpreted, therefore *the nearer the value of $\Pi$ is to zero, the better the effect of the solution.*

Obviously, the laws derived in this discussion are not depending to any experimental result. Whether or not the improved forms of Coulomb law and Newton's second law derived in this way will be suitable for other cases is a separate issue.

It should be noted that besides of solving equations, optimum-seeking methods could also be used for finding the minimum and the constant to be determined. In fact, the optimum seeking method will be used in this discussion.

### 7.6.4. Example for Deriving Coulomb Law and Newton's Second law

As shown in Fig.7, let a circle $O'$ denotes an electrified globe, $m$ denotes the mass of a small charged ball (treated as a point mass P), the electric charges of the globe and ball be $q_1$ and $q_2$ of positive and nega-tive polarities respectively, and there are positive and negative separately. We will assume that the small ball rolls along a long incline from A to B. Its initial velocity is zero and gravity and friction are neglected. Supposing that O'A is a plumb line, coordinate x is orthogonal to O'A, coordinate y is orthogonal to coordinate x (parallel to O'A), BC is orthogonal to O'A. The lengths of OA, OB, BC, and AC are all equal to H, and O'C equals the radius R of the globe.

In this example, the value of $v_P^2$ which is the square of the velocity for the ball located at point $P$ is investigated. To distinguish the quanti-ties, denote the value given by the improved forms of Coulomb law and Newton's second law as $v_P^2$, whilst $v'^2_P$ denotes the value given by the law of energy conservation, then Eq.(38) can be written as

$$\Pi = \int_{-H}^{0} \left( \frac{v_P^2}{v'^2_P} - 1 \right)^2 \mathrm{d}x = \min_0 \tag{50}$$

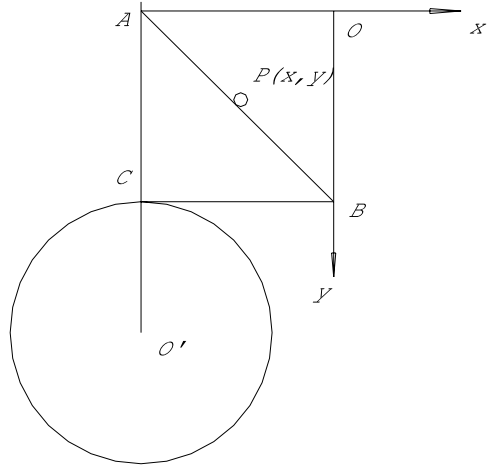

Fig.7  A small charged ball rolls from $A$ to $B$

Similar to the gravitational potential energy, from Eq.(43), the potential energy of the charged ball caused by electricity on the ball located at point $P$ is

$$V = -\frac{kq_1q_2}{(D-1)r_{O'P}^{D-1}} \tag{51}$$

According to the law of Energy Conservation, we can get

$$-\frac{kq_1q_2}{(D-1)r_{O'A}^{D-1}} = \frac{1}{2}mv_P'^2 - \frac{kq_1q_2}{(D-1)r_{O'P}^{D-1}} \tag{52}$$

And therefore

$$v_P'^2 = \frac{2kq_1q_2}{m(D-1)}[\frac{1}{r_{O'P}^{D-1}} - \frac{1}{(R+H)^{D-1}}] \tag{53}$$

Considering the general case, the rolling curve is

$$y = y(x)$$

(54)

For the ball located at point P,

$$\mathrm{d}v / \mathrm{d}t = a$$

(55)

because

$$\mathrm{d}t = \frac{\mathrm{d}s}{v} = \frac{\sqrt{1 + y'^2}\, \mathrm{d}x}{v}$$

Therefore

$$\mathrm{d}v = a\mathrm{d}t = a\,\frac{\sqrt{1 + y'^2}\, \mathrm{d}x}{v}$$

which gives

$$v\mathrm{d}v = a\sqrt{1 + y'^2}\, \mathrm{d}x$$

(56)

According to the improved form of Coulomb law, for point P, the attracted force acting on the ball is

$$F_{\mathrm{P}} = \frac{kq_1 q_2}{r_{\mathrm{O'P}}^{D}}$$

The force along to the tangent is

$$F_{\mathrm{a}} = \frac{kq_1 q_2}{r_{\mathrm{O'P}}^{D}}\,\frac{y'}{\sqrt{1 + y'^2}}$$

(57)

According to the improved form of Newton's second law, for point P, the acceleration along to the tangent is

$$a = \left(\frac{F_{\mathrm{a}}}{m}\right)^{1/1+\varepsilon} = \left(\frac{kq_1 q_2 y'}{mr_{\mathrm{O'P}}^{D}\sqrt{1 + y'^2}}\right)^{1/1+\varepsilon}$$

(58)

And from Eq.(56), it gives

$$v\mathrm{d}v = \left\{\frac{kq_1 q_2 y'}{m[(H+x)^2 + (R+H-y)^2]^{D/2}\sqrt{1+y'^2}}\right\}^{1/1+\varepsilon}\sqrt{1 + y'^2}\, \mathrm{d}x$$

(59)

For the two sides, we run the integral operation from A to P, for

$$v_{\mathrm{P}}^2 = 2\int_{-H}^{x_{\mathrm{P}}} \left\{\frac{kq_1 q_2 y'}{m[(H+x)^2 + (R+H-y)^2]^{D/2}\sqrt{1+y'^2}}\right\}^{1/1+\varepsilon}\sqrt{1 + y'^2}\, \mathrm{d}x$$

(60)

Considering the simplest case, the straight line between A and B is

$$y = H + x$$

(60a)

Substituting Eq.(60a) into Eq.(60), and setting $x = -z$ gives

$$v_{\mathrm{P}}^2 = 2 \int_{-x_{\mathrm{P}}}^{H} \{ \frac{kq_1 q_2}{m[(H-z)^2 + (R+z)^2]^{D/2}} \}^{1/1+\varepsilon} (\sqrt{2})^{\varepsilon/1+\varepsilon} \mathrm{d}z \qquad (60\mathrm{b})$$

then the value can be calculated by a method of numerical integral.

**Example 1. T**he given data is assumed to be: for the charged globe and ball, $\dfrac{kq_1 q_2}{m}$ =3.99×10$^{14}$m$^3$/s$^2$; the radius of the globe $R$ =6.37×10$^6$m, $H = R$ /10. Now try to solve the problem shown in Fig. 1, finding the solution for the value of $v_{\mathrm{B}}^2$, and derive the improved form of Coulomb law and the improved form of Newton's second law at once.

Firstly, according to the original form of Coulomb law and the original form of Newton's second law (i.e., let $D$ =2 in Eq.(43) and $\varepsilon$ =0 in Eq.(48)) and the law of Energy Conservation, all the related quantities can be calculated. Then substituting them into Eq.(50), gives

$\Pi_0$ =571.4215

Here, according to the law of Energy Conservation, it gives $v_{\mathrm{B}}^2$ =1.0767×10$^{7,}$ while according to the original form of Coulomb law and the original form of Newton's second law, it gives $v_{\mathrm{B}}^2$ =1.1351×10$^7$. The difference is about 5.4 %.

Since $\Pi_0$ is not equal to zero, then the values of $D$ and $\varepsilon$ can be determined by the optimum seeking method.

The optimum seeking methods can be divided into two kinds: (a) the problem is independent of the initial values – this is complicated, (b) ) the problem requires initial values -- this program is simple.

For type (b) approach, i.e., the searching method [106] will be used. Firstly, the value of $D$ is fixed so let $D$ =2, then search using $\varepsilon$ as $\varepsilon$ =0.0146. The value of $\Pi$ reaches the minimum of 139.3429, then the value of $\varepsilon$ is fixed. Search using $D$ =1.99989, the value of $\Pi$ reaches the minimum 137.3238. Then the value of $D$ is fixed. Search using $\varepsilon$ =0.01458, the value of $\Pi$ reaches the minimum of 137.3231. Because the last two results are very close to the first result, the search can be stopped. Then the result is as follows

$D$ =1.99989; $\varepsilon$ =0.01458; $\Pi$ =137.3231

Here the value of $\Pi$ is only 24% of $\Pi_0$. According to the law of Energy Conservation, it gives $v_B^2$ =1.0785×10$^7$, and according to the improved Coulomb law and the improved Newton's second law, it gives $v_B^2$ =1.1073×10$^7$. The difference is about 2.7 % only.

The results suitable for this example, with the constant dimension fractal form, are as follows:

- the improved Coulomb law, $f = \dfrac{kq_1q_2}{r^{1.99989}}$ ;

- the improved Newton's second law, $F = ma^{1.01458}$ .

Finally we discuss the dimension(unit) of the improved Coulomb Law and the improved Newton's second law. Two precepts can be given.

(i): To fix the dimensions of $a^{1+\varepsilon}$ and $r^{2-\varepsilon}$ use the same for $a^1$ and $r^2$ separately.

(ii): To handle the dimensions for each formula, multiply the right side by a factor, for example, the improved Newton's second law can be written as $F = K'ma^{1+\varepsilon}$, where the value of $K'$ is equal to 1, while the dimension of $K'$ should be chosen to make the dimensions of the left side and right side identical.

The first precept is used in this discussion for the advantage that the form of the formula cannot be changed, while for the second one the form of the formula will be changed. Of course, other precept also may be discussed further.

Therefore we conclude that the Coulomb law and Newton's second law could be modified, based upon experimental results. The example given herein shows that these original two laws should be improved. In order to derive these two laws theoretically, the law of Energy Conservation can be used. By way of an example (a small charged ball moves down along a long incline within the electric field due to a charged globe), the variable dimension fractal method was developed, and used to derive the improved Newton's second law and the improved form of Coulomb law at once.

Other Ref. (i) Newton's Second Law, the Law of Gravity and the Law of Coulomb with the Form of Variable Dimension Fractal. [111] (ii) Science of Conservation of Energy. [112]

## 6.7  Do fundamental constants in Nature vary with time?

Contributed by Fu Yuhua. (E-mail: fuyh@cnooc.com.cn)

One of interesting questions in recent years is whether the fundamental questions of Nature vary with time? Some theoretical physicists are convinced of this issue, while plenty others are not so sure. In this section we discuss how gravitational constant (G) may also vary. While it is not so conclusive yet, perhaps this section could inspire other approaches.

**Ref.**

http://news.yahoo.com/s/space/20060711/sc_space/scientistsquestionnatures fundamentallaws

Abstract: For the static state, Newton's law of gravity is correct, here the gravitational coefficient should be taken as a constant $G_0$ ( i.e., the gravitational constant in common meaning, and equals $6.67\cdots\times10^{-11}$ N•m²/kg² ) . For other cases, Newton's law of gravity will not be correct, for the condition that requiring the form of the law of inverse square for gravitation, then the changing rule for gravitational coefficient $G$ must be considered. In this paper, the expressions of $G$ for problem of advance of Mercury's perihelion, the problem of gravitational deflection of photon orbit around the sun and the example that a small ball moves down along a long incline are given. For problem of advance of Mercury's perihelion, $(1+5.038109\times10^{-8})G_0 \leq G \leq (1+1.162308\times10^{-7})G_0$; for problem of gravitational deflection of photon orbit around the sun, $G_0 \leq G \leq 2.5G_0$ ; for the example that a small ball moves down along a long incline on earth surface, the result given by constant dimension fractal is as follows: $1.001725G_0 \leq G \leq 1.001735G_0$ ; the more accurate result given by variable dimension fractal is as follows: $G_0 \leq G \leq 1.000012 G_0$ .

In reference [105], the following improved formula of universal gravitation was presented (according to the requirement of this paper, the original letter $G$ already rewrites as $G_0$ )

$$F = -\frac{G_0 Mm}{r^2} - \frac{3G_0^2 M^2 mp}{c^2 r^4} \qquad (61)$$

where : $G_0$ is the gravitational constant in common meaning, and equals $6.67 \cdots \times 10^{-11}$ N•m²/kg² ; M and m — masses of the two bodies ; r — the distance between the two bodies ; c — velocity of light ; p — half normal chord for body m moving around the body M with a curve, and the value of p reads

p = a(1-e²)    (for ellipse)
p = a(e²-1)    (for hyperbola)
p = y²/2x    (for parabola)

For the problem of gravitational deflection of photon orbit around the sun and the problem of advance of planet perihelion, by using the improved formula of universal gravitation, the same results as given by general relativity can be reached.

From Eq.(61), for the condition that requiring the form of the law of inverse square for gravitation, then we have

$$-\frac{GMm}{r^2} = -\frac{G_0 Mm}{r^2} - \frac{3G_0^2 M^2 mp}{c^2 r^4}$$

It gives the changing rule for gravitational coefficient $G$ as follows

$$G = G_0 (1 + \frac{3G_0 Mp}{c^2 r^2}) \qquad (62)$$

Therefore for problem of advance of Mercury's perihelion, we have

$(1 + 5.038109 \times 10^{-8}) G_0 \leq G \leq (1 + 1.162308 \times 10^{-7}) G_0$

For problem of gravitational deflection of photon orbit around the sun, we have

$G_0 \leq G \leq 2.5 G_0$

In reference [113], for an example of a small ball moves down along a long incline on earth surface, the results suitable for this example with the constant dimension fractal form were as follows

Improved law of gravity:

$$F = -\frac{G_0 Mm}{r^{1.99989}} \tag{63}$$

Improved Newton's second law:

$$F = ma^{1.01458} \tag{64}$$

From Eq.(63), for the form of the law of inverse square for gravitation, then we have

$$-\frac{GMm}{r^2} = -\frac{G_0 Mm}{r^{1.99989}}$$

It gives

$$G = G_0 r^{0.00011} \tag{65}$$

The result given by constant dimension fractal is as follows

$$1.001725 G_0 \leq G \leq 1.001735 G_0$$

In ref [114], considering that for the small ball at the beginning of static state, the Newton's law of gravity and Newton's second law are correct, then the more accurate results with the form of variable dimension fractal are as follows:

Improved law of gravity:

$$F = -G_0 Mm / r^{2-\delta} \tag{66}$$

where $\delta = 1.206 \times 10^{-12} x$, $x$ equals the horizontal distance of the small ball movement（for the beginning static state, $x = 0$）.

Improved Newton's second law:

$$F = ma^{1+\varepsilon} \tag{67}$$

where $\varepsilon = 8.779 \times 10^{-8} x$.

From Eq.(66), for the form of the law of inverse square for gravitation, then we have

$$-\frac{GMm}{r^2} = -G_0 Mm / r^{2-\delta}$$

It gives

$$G = G_0 r^{\delta} \tag{68}$$

Then the more accurate result given by variable dimension fractal is as follows

$G_0 \leq G \leq 1.000012 G_0$

## 6.8  Scale  invariance principle and coherent picture between microscales and macroscales

Provided we accept the 'fractal explanation' described in the preceding section, then does this 'scale-invariance principle' could be used to explain and reconcile the observed, great difference in the properties/behaviours of:

(i).  a few grains of sand vs. the vast terrain of the Sahara;

(ii).  individual particle vs transport phenomena;

(iii). micromolecule (e.g. $H_e$ or $H_z$) vs. macromolecule (such as DNA);

(iv). paternless or simple pattern of a small patch of something vs complicated --- albeit repeated --- pattern of beautiful fractals?

## 6.9  Does coral reef data support slowing-Earth-day hypothesis?

Some geophysicists argue that the oldest coral reef is 4200 years old, and from this viewpoint they put forward the argument that the age of the Earth cannot be much larger than 4200 years. Questions:

- What is the link (and how) between the age of a coral reef and the age of the Earth?
- Is there experimental observation supporting the hypothesis that coral reel data may imply that Earth-day may be less than 24 hours in the past? (see Ref. ii, iii).
- If yes, then what can be shown to validate this experimental observation with other related data? (see Ref. v, vi)
- Provided such data comparison exists, then what is the 'aggregate conclusions' of various observed data from coral reef and the age of the Earth?
- Is there any implication of this new finding of the age of the Earth to Earth's climatic changes?

## 6.10 Plausible linkage between Planckian quantization and quantized information

Here are some questions concerning plausible relations between Planckian quantization and quantized information. (See also R. Ogden, arXiv: cs.IT/0507032).

Questions:

- Is there valid theoretical relation between Planckian quantization and quantization of information?
- Does this relation provide a sound theoretical basis for quantum computation?
- Is there also relation between quantization of information and Maxwell's equations and superconductivity phenomena?
- Is there experimental observation supporting the above linkage between Planckian quantization and quantization of information?

Background argument:

It is known that a quantum liquid may exhibit quantum computation phenomena. Whether the quantum liquid has a valid theoretical link with quantized information, however, has not been explored adequately in the literature.

We begin with Landauer principle [85] that energy to erase a bit of information can be expressed as follows:

$$E_{erase} \approx kT.\ln(2) \tag{69}$$

which also has relation with Shannon's entropy theorem. It is argued that the same amount of energy is required to create a bit [85].

Now it is interesting to note here that a similar relation between Boltzmann's constant and the quantization of energy is known in the context of superconductivity model of the universe [86]:

$$\frac{1}{2}\hbar\omega \approx kT_c = kT_c.\ln(e) \qquad (70)$$

Equating both equations, now we get:

$$E_{erase} \approx kT.\ln(2) = kT_c.\ln(e) \qquad (71)$$

or

$$\frac{T_c}{T} = \frac{\ln(2)}{\ln(e)} \qquad (72)$$

Therefore one could conclude that there is a plausible relation between minimum temperature for a system to exhibit quantum computation phenomena and critical temperature of a superconductor [86].

Using the above argument, one could also hypothesize that the amount of energy required for creating or erasing a bit of information in the quantum (computation) universe will be quantized according to equation [88]. And because this quantized information has relation to Planckian quantization, then one could argue that this quantized energy to create/erase a bit of information should be carried via Planck mass.

Now it is also interesting to remark here, that there is a plausible link between fundamental four forces and Planck mass*, via [87]:

$$\frac{37}{57}\left(\frac{g_r.g_w}{\alpha G.10^6}\right) \approx m_p \qquad (73)$$

which seems to support our argument that the linkage between a quantum universe model and quantum computation may not be merely a hypothetical question.

We now return to equation (70):

$$E = \frac{1}{2}\hbar\omega \approx kT_c = kT_c.\ln(e) \qquad (74)$$

which can be written as:

$$\hbar \approx \frac{2}{\omega} k T_c . \ln(e)$$

(75)

By substitution of equation (72) into equation (75), one gets another plausible linkage:

$$\hbar \approx \frac{2}{\omega} k T . \ln(2) \tag{76}$$

Summarizing, we may conjecture a new hypothesis that Planck's constant may be *temperature-dependent*, in particular near the critical-temperature of superconductor [86]. Of course, this hypothesis should be verified by experiments.

<u>Note</u>:
*However, Crothers and Dunning-Davies have recently argued that Planck particle does not exist. [87b]

# 7   Postscript: A description of anomalous electromagnetic phenomena known as the Hutchison effect

> Mathematics is not yet ready for such problems..
> --Paul Erdos, *The American Mathematical Monthly*

In the preceding chapters, we have discussed some unsolved problems in various areas of science.

Now in this Postscript we will review an array of phenomena which is almost unknown in the present body of science. This article was adapted and rewritten (when necessary) from the homepage www.hutchisoneffect.biz, with the kind permission of John Hutchison (heffect@infinet.net). He is well-known for his research on the 'Hutchison effect', which will be described here. Some of his more recent research has been filmed and can be found at www.bluebookfilms.com. This postscript was included here to stimulate further research concerning an anomalous effect which is unknown to formal science.

(a) **Definition:**

The 'Hutchison effect' refers to a collection of phenomena discovered by inventor John Hutchison in 1979. Electromagnetic influences developed by a combination of electric power equipment, including Tesla coils, have produced levitation of heavy objects (including a 60-pound cannon ball), fusion of dissimilar materials such as metal and wood, anomalous heating of metals without burning of adjacent material, spontaneous fracturing of metals, and changes in the crystalline structure and physical properties of metals. The effects have been well documented on film and videotape, and witnessed many times by credentialed scientists and engineers, but are difficult to reproduce consistently.

Some phenomena which has been witnessed include: a super-strong molybdenum rod was bent into an S-shape as if it were soft metal; a length of high-carbon steel shredded at one end and transmuted into lead the other; a piece of PVC plastic disappeared into thin air; bits of wood became embed-ded in the middle of pieces of aluminium; and all sorts of objects levitated. (http://www.hutchisoneffect.biz/Research/pdf/ESJAug201997.pdf)

Another known phenomena related to the Hutchison Effect is the Mini Ark. The following chapter is the description of this experiment. Recently he was able to build a full-scale Ark (http://www.hutchisoneffectonline.com/). In this system the large Ark was loaded with an external Tesla Coil. This produced great power and in 2006 they levitated 1800 pounds of machinery for National Geographic TV.

### (b) Hutchison's mini Ark

It is known that all of the most ancient recorded writings throughout the cultures have a made reference to a certain box. This was a box of great magical power; this was a box that was used to transmute the morning dew into seed to feed an army that went about the planet.

John Hutchison had a theory about the magical box, and it goes something like this... The box itself is a capacitor, the poles antenna, the cheribum a spark gap, and the carpeting shrouding it a medium capturing ambient AC current, storing it, then slowly charging it back through the antenna. In other words it is a self charging capacitor static scavenging unit. The Discovery Channel decided to put John's theory to the test and hired Bruce Burges of Bluebook Films to make a documentary. The title of this documentary is The Ark of the Covenant Revealed. (www.bluebookfilms.com)

 What is quite amazing in this story is that The Ark worked. But that's not all, there were orbs and entities that became apparent, very apparent, in person and on camera. These entities have been described as "*Angelic Beings of Light*" flying into and out of a portal that ignited out of the plasma ark between the cheribum.

In this section we discuss how it could be possible and what is actually observed. The voltage will be measured and collected by rectennaes, noble gas will be generated, and captured, a camera mounted in the box, a scale checking its weight to see if it is constant, home made vacuum tubes with the captured noble gas, powered by rectennaes, are but a few examples of the theories of Hutchison.

We will discuss some basic sketches we are working from, then explain and clarify them. They started with 7 sheets of 3/4 " spruce plywood. Two 30" high Tesla coils wound in opposite directions, a turnable type transiever with a few extra turns, 6" pvc white plastic pipe, black has carbon in it... 1/2" copper pipe coils, bendable water pipe (we bought 2 rolls of 60' each), 6lbs of #22ga wire (we bought 10lbs), 1 gallon of shellac urethane, 6 spray cans, 6 cans of gold paint, 4 litres of polyester resin and two metres of cloth, and several sheets of sand paper.

The entire interior of the box is lined with copper, silver plated is better... 1/64" thick, First use the resin to seal the inside of the box with fibre glass cloth around the bottom corners. Allow the resin to fill the grooves in the rough plywood. Now it's water or (grapefruit juice) proof. Line the inside of the box with copper of any thickness, minimum 1/64". Stop 2.5" from the top inside and out to prevent voltage leaking. The bottom must be at least 2" off the ground, so glue and pin with wooden dowel some non conductive feet onto the bottom of the box.

The mercy seat interlocks into the top of the box, then 3 more pieces so the numbered pieces from bottom to top go as follows #5, #3, #7, #4. Cap that with a piece of granite 30 1/4" x 50 3/4". The quartz in the slab has a resonating effect. The measurements should be staggered to allow for your choice of non-conductive moulding. The bottom of the seat needs to be lined with copper, as well as the top. At the cheribum, or griffins, or eagles or whatever you use as a spark gap, leave a 2.5" space all around the base of each, if you use the Tesla Coil system, the griffins will sit in the top of each coil, making a spring loaded contact with the inside coil.

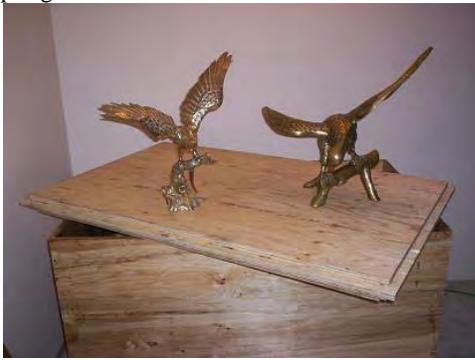

Figure 8. The Ark model of Hutchison (www.hutchisoneffect.biz)

There should be 10 turns in the secondary coils 1" apart, probably a couple extra but we want it there for experimentation. The suggested caps should be 12" lengths of 6" white pvc pipe with the end caps, lined inside and out with tin foil. The variation in the schematic, and possible adjustable spark gap - moveable or replaceable cheribum, is to increase or decrease the distance. This version can be juiced with a neon sign transformer to charge things up and activate the metals so they are more sensitive to ions.

When making the coils start 1" up from the bottom, drill a small hole. pull through a foot or two, shoe goop the hole closed, wind it perfectly until you get 1" from the top, then pull an extra foot or two through a small hole 1" from the top, shoe goop that closed too. Use urethane shellac to coat it liberally. The cheribum will make contact with the top wire inside the tube with a copper spring loaded contraption that compresses inside each tube when the lid (mercy seat) is lowered onto the pipe caps. Use springs brass studs, and silver dollars with holes drilled in them.

We settled on how the poles should be connected, one should be insulated from the outside of the box, and make a connection with the inside (positive); the other should be making contact with the outside (negative) --you wouldn't want to pick it up when charged... might not be the best way to do it, but we want to see if it collects more that way, or if both should be wired to the inside. Time will tell, for those that want to make an Ark, make it in such a way that it can be varied, rewired, and experimented with.

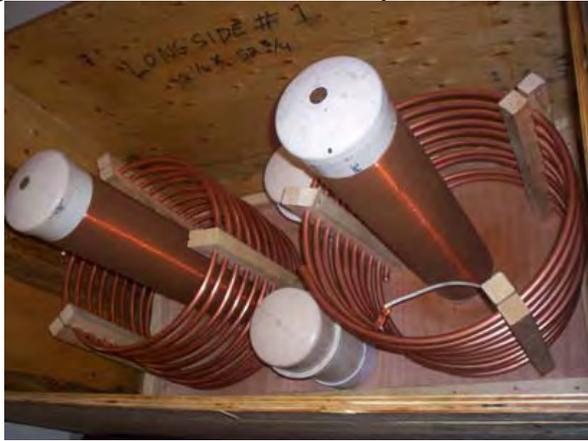

Figure 9. The poles inside the box.

When you have finished winding things in the opposite directions, mount the apparatus on a piece of plywood about 21"x51"x3/4" thick. This way nothing changes positions between experiments... The "Warp Core" can be ejected if you want to change up or maintain it. Screws hold it together until you can drill and glue in 5/16" dowels in their place. Heavy string works to bind the secondary coils, a piece of 2x4 in the centre of the pipe, pinned and glued holds the centre where it is. Once everything is glued, take a gallon of urethane and coat the whole thing evenly, especially the wood, and 22ga wire, this will keep everything intact. John uses 50 coats for durability, 6 minimum, you decide where to stop.

A warning... the box, or Ark it's self itself reactive, meaning it will store and collect energy, lots of it, which can be gathered slowly just from your carpet, a passing storm... once you put the lining in it, it is active unless you short it out! Lethal folks, that's how and why people who make these things die, they think it needs a charge to be lethal, wrong. The coils themselves are receivers and react similarly, especially in storms. As you assemble this keep everything shorted out until you get it in a place where it is safe to play with a few million volts or you won't live very long... even during travel. This has not yet been tested in proximity to gasoline either, so be careful...

You can build a smaller version, apparently John fired it up in his lab and little orbs and entities became visible. As for grounding, just get some wire and hook alligator clips to it and short it out, grounding won't help, short the outside to the inside of the box or it will collect, if you ground the outside to your water pipe for instance you just helped it charge, the outside is negative, inside positive. Keep the poles out of it when not in use, as they are like antennae, so they collect too. Short out all of the secondaries to the cheribum too. Look up salt water capacitor, it may help you understand what the box is. Also look up home made Tesla coils too.

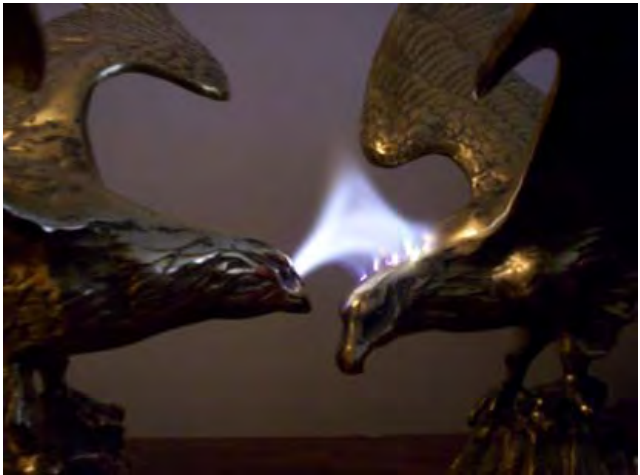

Figure 10. This is what they observed: the Hutchison Effect

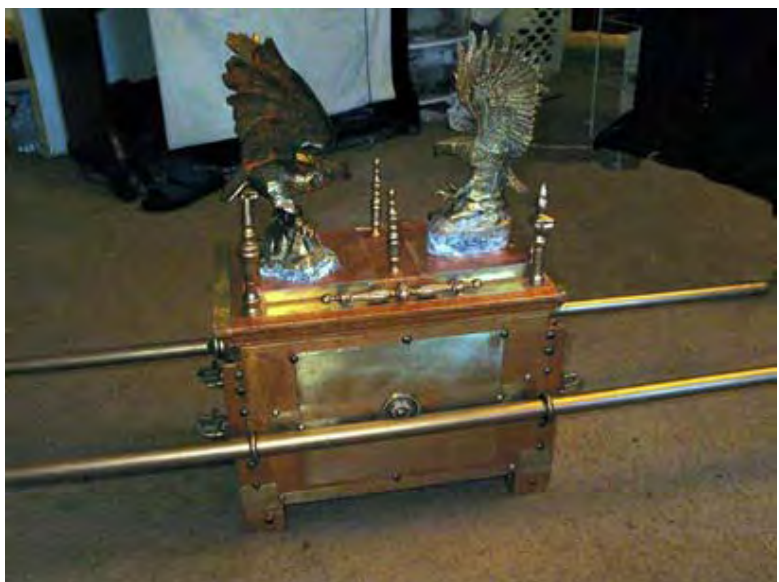

Figure 11. After all those steps, you get your Mini Ark

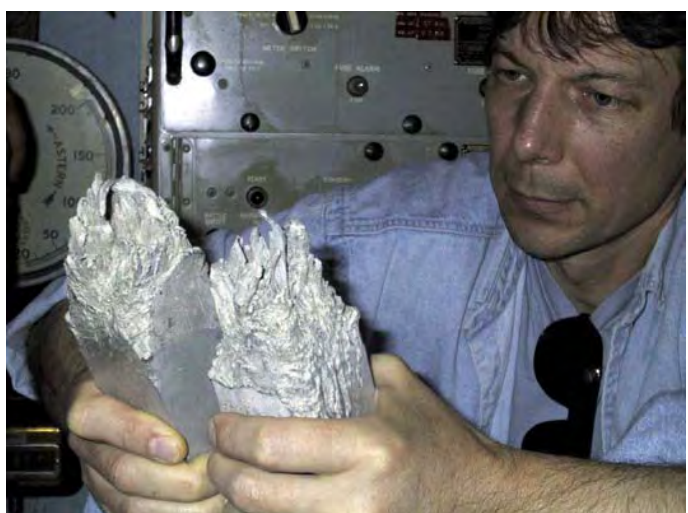

Figure 12. Allumunium bar jellified by the Hutchison effect.

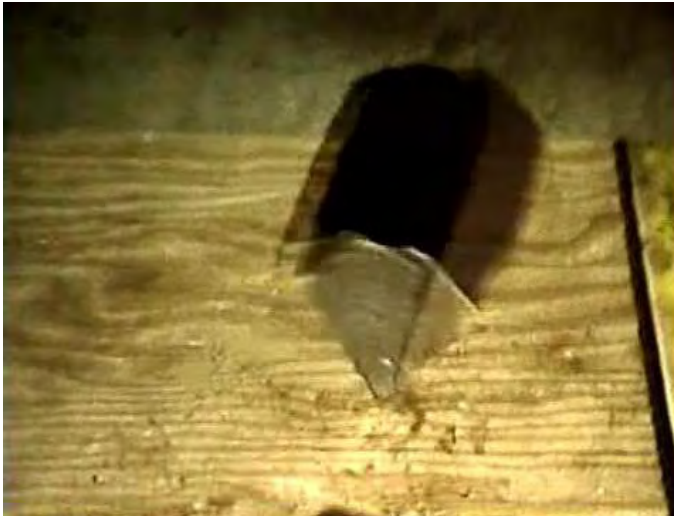

Figure 13. Very first footage of the Hutchison Effect. The sample wobbling in and out of existence. (http://www.hutchisoneffect.biz)

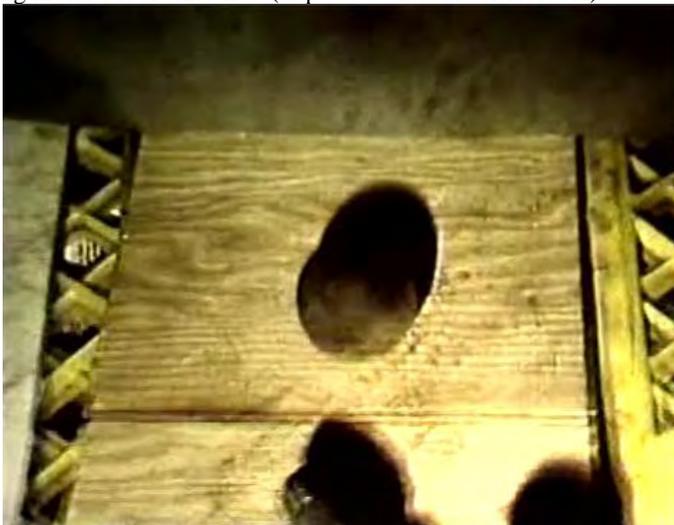

Figure 14. The disappearing cannonball.

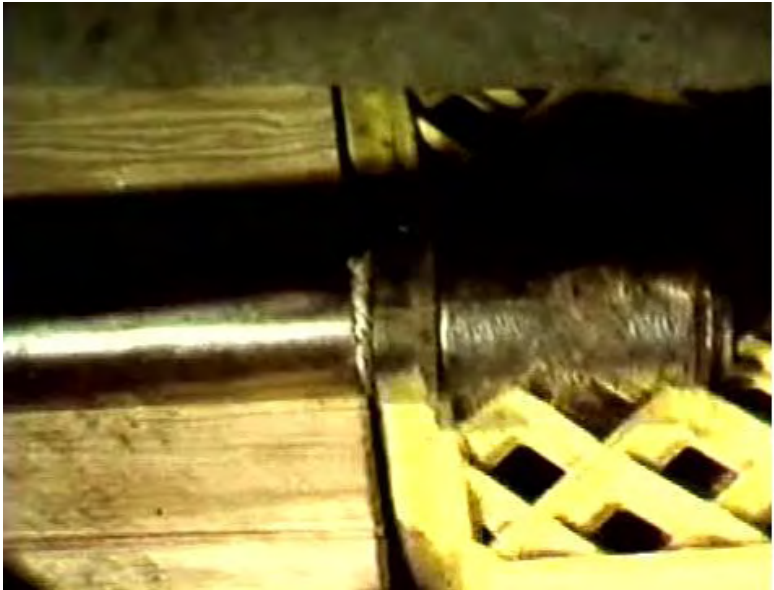

Figure 15. A cast iron sample phasing out of existence.

We have introduced here a simplified description of the Hutchison effect as reported in his homepage (www.hutchisoneffect.biz). Of course, those who want to make a replica must be very cautious as to its effects.

Hopefully, some of these examples could be interesting enough to attract further research on these anomalous effects.

WARNING: The re-creation of certain aspects of these experiments can cause all many as yet unknown and uncharted scalar effects such as *temporal distortion*. Exercise the utmost caution if you intend to rebuild any of these experiments.

Ref.:

# Epilogue

> It isn't that they can't see the solution. It is that they can't see the problem.
>                                                        --G. K. Chesterton

Throughout this book, we discuss some unsolved problems in various branches of science, including mathematics, theoretical physics, etc.

It is our hope that some of the problems discussed in this book will find their place either in theoretical exploration or further experiments, while some problems may be useful only for scholarly stimulation.

FS, VC, FY, RK, JH
August 28th, 2006

*First version: 18th April, 2006. 1st revision: 8th May 2006.2nd revision: 16th June 2006; 3rd version: 28th June 2006, 4th version: 28th Aug. 2006*

# Acknowledgement

The lead author would like to thank Dmitri Rabounski and Dr. Larissa Borissova for their valuable comments on this book.

The co-author would like to thank to Profs. C. Castro, M. Pitkänen, and A. Kaivarainen, for insightful discussions and remarks; and also a number of journal editors and reviewers. This book is also dedicated to all physicist fellows at www.archivefreedom.org and www.sciprint.org, who contributed via discussion and in numerous other ways.

The authors would also like to thank to Prof T. Love and Steven Crothers for proof-reading the draft and various suggestions for improvement. Of course, any errors are solely attributed to the authors.

Florentin Smarandache
University of New Mexico
Gallup, NM 87301, USA
**smarand@unm.edu**

V. Christianto
Jakarta, Indonesia
http://www.sciprint.org, email: admin@sciprint.org

Fu Yuhua – China Offshore Oil Research Centre, P.O. Box 4728, Beijing, 100027, China (email: fuyh@cnooc.com.cn)

R. Khrapko -- Moscow Aviation Institute

John Hutchison – Independent Inventor, **www.hutchisoneffect.biz,** (**heffect@infinet.net**)

# Appendix A:

# Observation of anomalous potential electric energy from distilled water under solar heating

In order to stimulate both further theoretical and experimental research pertaining to new alternative energy, we discuss in this Appendix section a very simple experiment with distilled water. It can be shown that such an experiment will exhibit anomalous potential electric energy. Whereas the result is less impressive compared to the common LENR/CANR experiments, it is recommended to carry out further research along this direction.

### Introduction

There has been a somewhat regained interest for the alternative energy technologies based on *low-energy chemical-aided reaction* [1]. This process includes various different methods ranging from the well-known gas discharge process until the exotic processes such as microwave-induced reaction.[2][3] Some theoretical explanation has also been proposed in recent years. [4][5]

Nonetheless, from the viewpoint that our Earth is presently seeking a rapid change to alternative energy, one could imagine that it is required to find a 'less-exotic' energy source, which can be generated with minimum preparation. Therefore, the definition of 'energy input' term shall also include the energy amount needed to make preparation for the source and also for the equipment.

In this regard, we re-visit a well-known process of finding excess electrical energy out of 'distilled water.' It can be shown via experiment, that with very minimum preparation one can obtain anomalous excess electrical energy from distilled water, in particular under solar (photon) exposure. The result is summarized in Table 1.

In the last section we will discuss a few alternative approaches to explain this observed anomalous effect, for instance using the concept of 'zero point energy' of the phion-fluid condensate medium. [6]

Nonetheless, further experiment is recommended in order to verify or refute our proposition as described herein.

## Experimental preparation and result

The basic idea of this experiment comes from reading various papers related to *chemical aided reaction* [1][2]. There is also an abstract requirement for minimum preparation energy, so that it would be easier for rapid implementation (if chance permits).

Therefore we come to analogue to dc battery: *a used battery will re-gain part of its electric energy once it is put under exposure to the Sun light for a few hours*. This analogy leads us to hypothesize that the Sun light emits photon flux with sufficient 'zero point energy' which could trigger chemical reaction in the electrolyte. Then the re-gained electric energy of the used battery will last for a few days more.

Possible implication for this experiment could include usage of distilled water as an efficient method for battery charger, while possible future use in transportation etc. remains open.

So, in this simple experiment we consider a few alternative scenarios:

(i)     ordinary water without exposure to Sun light or to external dc potential (as control for this experiment);

(ii)    ordinary water with exposure to Sun light;

(iii)   distilled water without exposure to Sun light or to external dc potential;

(iv)    distilled water with exposure to Sun light;

(v)     distilled water with exposure to carbon alkali (chemical inside battery);

(vi)    distilled water with exposure to external dc potential;

(vii)   distilled water with exposure to Sun light and carbon alkali (chemical inside battery);

(viii)  distilled water with exposure to carbon alkali and to external dc potential.

Distilled water is used in this experiment instead of heavy-water (deuterium) which is commonly used in LENR experiment [1][2], with simple reason that it is easier to obtain almost anywhere. Therefore no excessive preparation for such water is needed. Of course, for better result it is recommended to conduct this experiment with heavy-water. (For instance, Belyaev *et al*. already conducted various experiments with heavy-water.)

The preparation for this experiment is described as follows.

We use 20 mm-diameter aluminium tube and fill it with ordinary water for control, then we measure its electrical resistance and also its electrical voltage (Type iA experiment). Then we put this tube under the exposure of Solar daylight (high noon), and using a 60mm-diameter magnifying lens at its focal distance in order to focus the Solar's photon flux into our tube. Then we measure again the electrical resistance and also its electrical voltage. (Type iB experiment)

We use another 20 mm-diameter aluminium tube and fill it with distilled water, then we put these tubes under the exposure of Solar daylight (Type iiB). Thereafter we repeat the procedure once again after introducing an external 1.5V DC potential into the electrolytes. Then we measure again the electrical resistance and also its electrical voltage. (Type iiC) After around 5-10 minutes, we release the external potential (1.5 DC volt) and put the tube again under solar light exposure. (Type iiD)

Then, we repeat the procedure after filling the tube with carbon alkali from used-batteries 1.5V DC. Then we measure again the electrical resistance and also its electrical voltage. (Type iiiA) Thereafter we repeat the procedure once again after introducing an external 1.5V DC potential into the electrolytes. Then we measure again the electrical resistance and also its electrical voltage. (Type iiiC) After around 5-10 minutes, we release the external potential (1.5 DC volt) and put the tube again under solar light exposure. (Type iiiD)

The experimental configuration is shown in the following diagrams, both with and without external 1.5Volt DC potential.

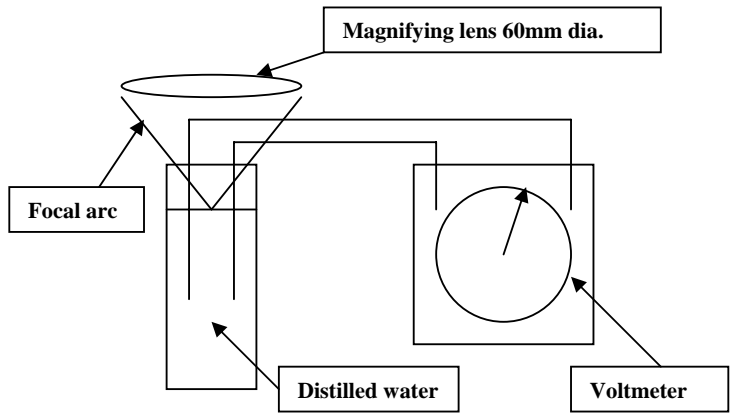

Diagram A1. Experiment with distilled water and no external DC (Type iiA)

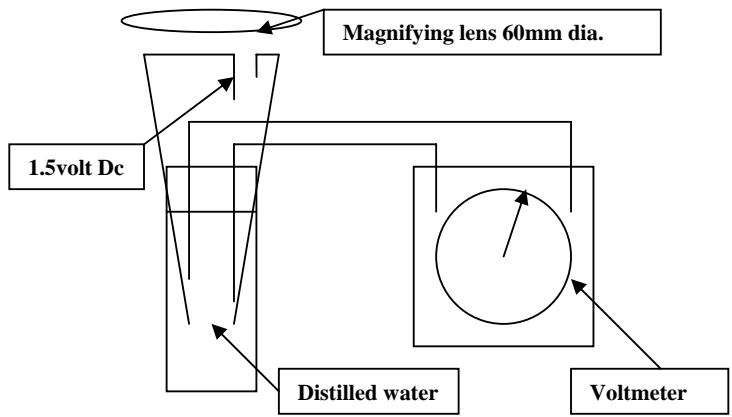

Diagram A2. Experiment with distilled water and external 1.5V DC (iiC)

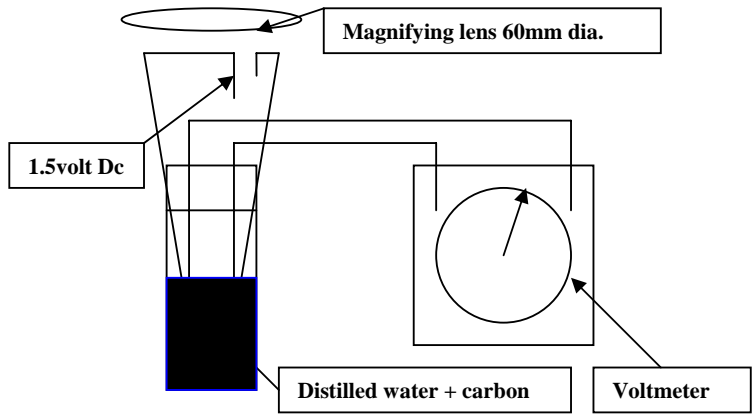

Diagram A3. Experiment with distilled water with carbon alkali and external 1.5Volt DC (Type iiiC + iiiD)

In simple words, in this experiment we want to know whether the effect of Solar heating (photon flux) is similar compared to the resulting effect from introducing carbon alkali material or introducing 1.5V DC potential into the electrolytes. As shown in Table A.1 below, it turns out that both photon flux and external 1.5V DC potential could induce significant impact to the observed anomalous potential, while carbon alkali almost has no further effect (at least to the experimental configuration as described herein).

The experiment was conducted in the backyard, around 21$^{st}$ Aug. 2006.

**Table A.1. Observation result with distilled water**

| Description | Without solar exposure | With solar exposure (magnifying lens) | Before external 1.5V DC. Without solar exposure | After external 1.5V DC. With solar exposure (magnifying lens) |
|---|---|---|---|---|
| | **A** | **B** | **C** | **D** |
| Ordinary water [i] | V=0 Volt; R>>1000Ω | V=0 Volt; R>>1000Ω | | |
| Distilled water [ii] | V=0 Volt; R>>1000Ω | V=0.2 Volt; R=600Ω ~1000Ω | V=0.8-1.0 Volt; R=600Ω ~1000Ω | V=0.6-0.8 Volt; R=100Ω ~600Ω |
| Distilled water with carbon alkali material [iii] | V=0.2 Volt; R>>1000Ω | V=0.6 Volt; R=600Ω ~1000Ω | V=0.6-0.8 Volt; R=600Ω ~1000Ω | V=0.6-0.8 Volt; R=100Ω ~600Ω |

From Table A.1 we can observe a few interesting results, as follows:

(i)     That *within bounds of experimental precision limits* we observe that there is anomalous potential energy in distilled water as much as 0.6-0.8 Volt (DC) after sufficient exposure to solar light, and after a few minutes introducing external 1.5Volt (DC) potential into the electrolytes. (Type iiC)

(ii)    Using carbon alkali material will add no further effect into this anomalous observed potential energy (Type iiiC). The exact source of this observed anomalous potential energy remains unknown.

(iii)   Furthermore, it is also interesting to note here that after around two hours (the external 1.5Volt DC potential has been released), measurement reading for configuration [iiD] remains showing anomalous potential electric energy ~ 0.4-0.6 Volt and resistance R=~100Ω.

(iv)     After around 24 hours (the next day), measurement reading for configuration [iiD] remains showing anomalous potential electric energy ~ 0.1-0.2 Volt and resistance R=~100Ω.

(v)     Therefore we can conclude to summarize this experimentation, that a small DC potential and photon flux (Solar light) could play significant role in the LENR/CANR-type processes which so far this effect has been almost neglected in reported LENR/CANR experiments.[1][2]

For clarity, we draw diagram showing observed anomalous potential energy (the lower bound value) in experiment type iiA, iiB, iiC, iiD for the first 24 hours of this experiment (Table A.2 and diagram A.4). It is clear here that the peak of anomalous potential energy was observed after introducing external 1.5Volt DC potential, and its impact not last yet after around 24 hours.

**Table A.2. Observation result in each step of experiment Type ii**

| Step | Hours | Observed potential (volt) |
|---|---|---|
| Without solar light | 0 | 0 |
| After solar light | 0.2 | 0.2 |
| With external 1.5Volt | 0.4 | 0.8 |
| Without external 1.5Volt, after solar light | 0.5 | 0.6 |
| After 2 hours | 2.5 | 0.4 |
| After ~24 hours | 24 | 0.1 |

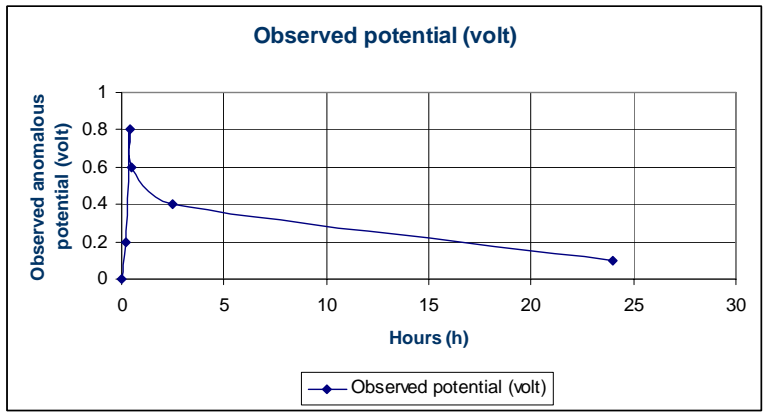

**Diagram A.4. Observation result in each step of experiment Type ii**

In our opinion, it is very likely that this *photon flux* could trigger effect just like in '*photo-synthesis' process* which is known in various biological forms of life. However, this proposition requires further theoretical considerations.

If this proposition corresponds to the facts (concerning the role of photo-synthesis), then perhaps this experiment does not belong to typical LENR-CANR experiments [1][[2], instead it is perhaps more convenient to call it PSCR (*PhotoSynthesis-catalyzed Chemical Reaction*).

Nonetheless, it should also be noted here that there is shortcoming of this experimentation, for instance we don't exactly measure how much carbon alkali material has been introduced into the electrolyte, nor how long the solar light exposure shall be maintained (it could take 5-10 minutes). It is because this experiment is merely to assess the viability of the idea, instead of becoming a rigorous experiment. Further experiments are of course suggested to verify this proposition with better precision.

**A few alternative interpretations of the above anomalous effect**

In order to explain the above anomalous potential energy, we consider a few possible alternative interpretations, as follows:

- photon magnified energy;
- photon Hall effect;
- photon condensate's zero point energy;
- phion condensate's Gross-Pitaevskii energy.

The rationale for each of these alternatives is discussed as follows:

(a) Photon Magnified Energy. It can be shown by the use of special relativity that the energy momentum relation actually also depends on the 'scale' of the frame of reference. Therefore the use of magnifying lens that focuses photon energy in the electrolyte will be not the same again with E=p.c for the area of magnifying lens, but:

$$E_{focused} = n^2 . E_{photon-flux} \tag{A.1}$$

Where n represents scaling factor, similar to refractive index.

(b) Photon Hall effect. It is known that photon takes the form of boson [13][13a]. Now it is possible also to assume that the photon condensate will induce Hall effect [12a], therefore we could use total particle momentum expression instead of conventional momentum [12a]:

$$p = mv + m\Omega \times r + qA \tag{A.2}$$

Therefore the energy-momentum relation becomes:

$$E = pc = (mv + m\Omega \times r + qA).c \tag{A.3}$$

If we neglect the first term (assuming photon is massless), then:

$$E = pc = (qA).c \tag{A.4}$$

We shall note here that Vigier and others suggested photon has mass.

(c) Photon condensate's Zero Point Energy. Starting with the assumption that photon is bosonic, then we could also use zero point energy of bose condensate for photon [13]. It is also known that zero point energy could play significant role in LENR experiments [2]. The zpe for bose condensate could be expressed as follows [13, p.13]:

$$\in = \frac{1}{v} \left\langle \hat{H}_{QFT} \right\rangle_{vac} \tag{A.5}$$

Nonetheless it is not yet clear, how zpe could trigger anomalous effect. This zpe could have linkage with interpretation of Dirac's negative energy [5].

(d) Phion condensate's Gross-Pitaevskii energy. We could also start with assumption that there exists phion fluid medium which is unobserved [6][14]. Recent paper by Moffat [6a] has shown that phion condensate model is at good agreement with CMBR temperature and also with galaxies rotation curve data. It could also be shown that using Gross-Pitaevskii equation one could derive Schrödinger equation, also planetary quantization.[11] Using the mechanism of photon-photon interaction [6], the solar's photon flux interacts with the surrounding phion condensate medium. And therefore the energy collected by the magnifying lens is not only its own 'photon flux' energy but also includes the energy of the phion condensate medium. This energy then triggers chemical reaction in the electrolyte. It is known that Ginzburg-Landau (Gross-Pitaevskii) equations has free energy term due to its nonlinear effect, therefore it perhaps could explain why the effect on the electrolyte remains quite significant (more than 0.2volt) after a few hours.

Further experiments are of course recommended in order to verify or refute these alternative explanations.

**Concluding remarks**

We have described here an experiment which could exhibit anomalous electrical energy in distilled water with *very minimum preparation energy*. While this observed excess energy here is less impressive than [1][2] and the material used is also far less exotic than common LENR/CANR experiments, from the viewpoint of minimum preparation requirement –and therefore less barrier for rapid implementation--, it seems that further experiments could be recommended in order to verify and also to explore various implications of this new proposition.

Practical implications of this experiment could include possibility for using distilled water+carbon alkali for *battery charger*, as an alternative to *polymer electrolyte charger* (PEFC) method introduced by DoCoMo by July this year (2006).

We shall note here that perhaps this experiment does not belong to 'standard' LENR-CANR experiments [1][[2], instead it is perhaps more convenient to call it PSCR (*PhotoSynthesis-catalyzed Chemical* Reaction). Nonetheless, the present simple experiment was reported merely to encourage further experiments along similar line of thought.

# Appendix B:

## On the origin of macroquantization in astrophysics and celestial motion

ABSTRACT. Despite the use of Bohr radius formula to predict celestial quantization has led to numerous verified observations, the cosmological origin of this macroquantization remains an open question. In this article various plausible approaches are discussed. Further observation to verify or refute this proposition is recommended, in particular for exoplanets.

RÉSUMÉ:  En dépit de l'utilisation de la formule de rayon de Bohr de prévoir la quantification céleste a mené aux nombreuses observations véri-fiées, l'origine cosmologique de ce macroquantization est une question en suspens. En cet article de diverses approches plausibles sont discutées. Promouvez l'observation pour vérifier ou réfuter cette proposition est re-commandée, en particulier pour des exoplanets.

**Introduction**

It is known that the use of Bohr radius formula [1] to predict celestial quantization has led to numerous verified observations [2][3]. This approach was based on Bohr-Sommerfeld quantization rules [4][5]. Some implications of this quantum-like approach include exoplanets prediction, which has become a rapidly developing subject in recent years [6][7]. While this kind of approach is not widely accepted yet, this could be related to a recent sug-gestion to reconsider Sommerfeld's conjectures in Quantum Mechanics [8].

While this notion of macroquantization seems making sense at least in the formation era of such celestial objects, i.e. "*all structures in the Universe, from superclusters to planets, had a quantum mechanical origin in its earliest moments*" [9], a question arises as to how to describe the physical origin of wave mechanics of such large-scale structures [5].

A plausible definition of the problem of quantization has been given by Grigorescu [10]: "select an infinite, discrete number of quantum possible real motions, from the continuous manifold of all mechanically possible motions." While this quantization method has been generally acceptable to describe physical objects at molecular scale, there is not much agreement why shall we also invoke the same notion to describe macrophenomena, such as celestial orbits. Nonetheless, there are plenty efforts in the literature in attempt to predict planetary orbits in terms of wave mechanics, including a generalisation of Keplerian classical orbits [11].

In this article we discuss some plausible approaches available in the literature to describe such macroquantization in astrophysics, in particular to predict celestial motion:

  a.  Bohr-Sommerfeld's conjecture;
  b.  Macroquantum condensate, superfluid vortices;
  c.  Cosmic turbulence and logarithmic-type interaction.

While these arguments could be expected to make the notion of macro-quantization a bit reasonable, it is beyond the scope of this article to conclude which of the above arguments is the most consistent with the observed data. There is perhaps some linkage between all of these plausible arguments. It is therefore recommended to conduct further research to measure the reliability of these arguments, which seems to be worthwhile in our attempt to construct more precise cosmological theories.

**Bohr-Sommerfeld's quantization rules**

In an attempt to describe atomic orbits of electron, Bohr proposed a conjecture of quantization of orbits using analogy with planetary motion. From this viewpoint, the notion of macroquantization could be considered as returning Bohr's argument back to the celestial orbits. In the meantime it is not so obvious from literature why Bohr himself was so convinced with this idea of planetary quantization [12], despite such a conviction could be brought back to Titius-Bode law, which suggests that celestial orbits can be described using simple series. In fact, Titius-Bode were also not the first one who proposed this kind of simple series [13], Gregory-Bonnet started it in 1702.

In order to obtain planetary orbit prediction from this hypothesis we could begin with the Bohr-Sommerfeld's conjecture of quantization of angular momentum. As we know, for the wavefunction to be well defined and unique, the momenta must satisfy Bohr-Sommerfeld's quantization condition [14]:

$$\oint_{\Gamma} p.dx = 2\pi.n\hbar \tag{1}$$

for any closed classical orbit $\Gamma$. For the free particle of unit mass on the unit sphere the left-hand side is

$$\int_0^T v^2 .d\tau = \omega^2 .T = 2\pi.\omega \qquad (2)$$

where T=$2\pi/\omega$ is the period of the orbit. Hence the quantization rule amounts to quantization of the rotation frequency (the angular momentum): $\omega = n\hbar$. Then we can write the force balance relation of Newton's equation of motion:

$$GMm / r^2 = mv^2 / r \qquad (3)$$

Using Bohr-Sommerfeld's hypothesis of quantization of angular momentum (2), a new constant g was introduced:

$$mvr = ng / 2\pi \qquad (4)$$

Just like in the elementary Bohr theory (before Schrödinger), this pair of equations yields a known simple solution for the orbit radius for any quantum number of the form:

$$r = n^2 .g^2 /(4\pi^2 .GM .m^2) \qquad (5)$$

or

$$r = n^2 .GM / v_o^2 \qquad (6)$$

where r, n, G, M, $v_o$ represents orbit radii (semimajor axes), quantum number (n=1,2,3,…), Newton gravitation constant, and mass of the nucleus of orbit, and specific velocity, respectively. In this equation (6), we denote

$$v_o = (2\pi / g) .GMm \qquad (7)$$

The value of m is an adjustable parameter (similar to g).

Nottale [1] extends further this Bohr-Sommerfeld quantization conjecture to a gravitational-Schrödinger equation by arguing that the equation of motion for celestial bodies could be expressed in terms of a scale-relativistic Euler-Newton equation. For a Kepler potential and in the *time independent* case, this equation reads (in Ref [1c] p. 380):

$$2D^2 \Delta\Psi + (E / m + GM / r).\Psi = 0 \qquad (8)$$

Solving this equation, he obtained that planetary orbits are quantized according to the law:

$$a_n = GMn^2 / v_o^2 \qquad (9)$$

where $a_n$,G,M,n,$v_o$ each represents orbit radius for given n, Newton gravitation constant, mass of the Sun, quantum number, and specific velocity ($v_o$=144 km/sec for Solar system and also exoplanet systems), respectively.

These equations (8)-(9) form the basis of Nottale's Scale Relativity prediction of planetary orbits [1]; and equation (9) corresponds exactly with equation (6) because both were derived using the same Bohr-Sommerfeld's quantization conjecture. Another known type of observed quantization in astronomy is Tifft's 72 km/sec quantization [13].

**Macroquantum condensate, superfluid vortices**

Provided the above Bohr-Sommerfeld description of macroquantization corresponds to the facts, then we could ask further what kind of physical object could cause such orbital quantization. Thereafter we could come to the macroquantum condensate argument. In this regard, astrophysical objects could be seen as results of vacuum condensation [15][16]. For instance Ilyanok & Timoshenko [17] took a further step by hypothesizing that the universe resembles a large Bose Einstein condensate, so that the distribution of all celestial bodies must also be quantized. This conjecture may originate from the fact that according to BCS theory, superconductivity can exhibit macroquantum phenomena [18]. There is also a known suggestion that the vacua consist of hypercrystalline: *classical spacetime coordinate and fields are parameters of coherent states* [19].

It is perhaps interesting to remark here that Ilyanok & Timoshenko do not invoke argument of *non-differentiability* of spacetime, as Nottale did [1]. In a macroquantum condensate context, this approach appears reasonable because Bose-Einstein condensate with Hausdorff dimension $D_H \sim 2$ could exhibit fractality [20], implying that non-differentiability of spacetime conjecture is not required. The same fractality property has been observed in various phenomena in astrophysics [21], which in turn may also correspond to an explanation of the origin of multifractal spectrum as described by Gorski [22]. In this regard, Antoniadis *et al.* have discussed CMBR temperature (2.73° K) from the viewpoint of conformal invariance [23], which argument then could be related to Winterberg's hypothesis of superfluid Planckian phonon-roton aether [24].

Based on previous known analogy and recent research suggesting that there is neat linkage between gravitation and condensed matter physics [25][26], we could also hypothesize that planetary quantization is related to quantized vortex. In principle, this hypothesis starts with observation that in quantum fluid systems like superfluidity, it is known that such vortexes are subject to quantization condition of integer multiples of $2\pi$, or $\oint v_s.dl = 2\pi.n\hbar / m_4$. Furthermore, such quantized vortexes are distributed in equal distance, which phenomenon is known as vorticity [4]. In large superfluid system, usually we use Landau two-fluid model, with normal and superfluid component. The normal fluid component always possesses some non-vanishing amount of viscosity and mutual friction. Similar approach with this proposed model has been considered in the context of neutron stars [27], and this quantized vortex model could also be related to Wolter's vortex [28].

### Cosmic turbulence and logarithmic type interaction

Another plausible approach to explain the origin of quantization in astronomy is using turbulence framework. Turbulence is observed in various astrophysical phenomena [21], and it is known that such turbulence could exhibit a kind of self-organization, including quantization.

Despite such known relations, explanation of how turbulence could exhibit orbital quantization is not yet clear. If *and only if* we can describe such a flow using Navier-Stokes equation [29], then we can use R.M. Kiehn's suggestion that there is exact mapping from Schrödinger equation to Navier-Stokes equation, using the notion of quantum vorticity [30]. But for fluid which cannot be described using Navier-Stokes equation, such exact mapping would not be applicable anymore. In fact, according to Kiehn the Kolmogorov theory of turbulence is based on assumption that the turbulent state consists of "vortices" of all "scales" with random intensities, but it is not based on Navier-Stokes equation explicitly, in fact "*the creation of the turbulent state must involve discontinuous solutions of Navier-Stokes equations.*" [31] However, there is article suggesting that under certain conditions, solutions of 3D Navier-Stokes equation could exhibit characteristic known as Kolmogorov length [32]. In this kind of hydrodynamics approach, macroquantization could be obtained from solution of diffusion equation [33].

In order to make this reasoning of turbulence in astrophysics more consistent with the known analogy between superfluidity and cosmology phenomena [26], we could also consider turbulence effect in quantum liquid. Therefore it seems reasonable to consider *superfluid turbulence* hypothesis, as proposed for instance by Kaivarainen [34]. There are also known relations such as discrete scale invariant turbulence [35], superstatistics for turbulence [36], and conformal turbulence. Furthermore, such a turbulence hypothesis could lead to logarithmic interaction similar to Kolmogorov-type interaction across all scales [28].

Another way to put such statistical considerations into quantum mechanical framework is perhaps using Boltzmann kinetic gas approach. It is known that quantum mechanics era began during Halle conference in 1891, when Boltzmann made a remark: "*I see no reason why energy shouldn't also be regarded as divided atomically.*" Due to this reason Planck subsequently called the quantity $2\pi\hbar$ after Boltzmann – 'Boltzmann constant.' Using the same logic, Mishinov *et al.* [37] have derived Newton equation from TDGL:

$$m * d_t V_p(t) = e * .E - m * V_p(t) / \tau_p \qquad (10)$$

This TDGL (time-dependent Ginzburg-Landau) equation is an adequate tool to represent the *low-frequency* fluctuations near $T_c$, and it can be considered as more universal than GPE (Gross-Pitaevskii equation).

### Concluding remarks

In this article, some plausible approaches to describe the origin of macro-quantization in astrophysics and also celestial motion are discussed. While all of these arguments are interesting, it seems that further research is required to verify which arguments are the most plausible, corresponding to the observed astrophysics data.

After all, the present article is not intended to rule out the existing methods in the literature to predict quantization of celestial motion, but instead to argue that perhaps this macroquantization effect in various astronomy phenomena requires a new kind of theory to describe its origin.